\documentclass[11pt]{article}
\usepackage{amssymb,amsmath,amsthm,tikz,multirow,nccrules,graphicx,subfig,overpic,indentfirst}
\usetikzlibrary{arrows,calc}

\title{Tilings of the sphere by congruent quadrilaterals III: edge combination $a^3b$ with general angles}
\author{Yixi Liao, Pinren Qian, Erxiao Wang\thanks{Corresponding author (wang.eric@zjnu.edu.cn). Research was supported by Key projects of Zhejiang Natural Science Foundation No. LZ22A010003 and ZJNU Shuang-Long Distinguished Professorship Fund No. YS304319159.}, Yingyun Xu \\
	Zhejiang Normal University}

\usepackage[hidelinks, hyperindex]{hyperref}

\hyperref[sec:function]{}

\newcommand{\sub}{\subset}

\usepackage{amsmath}  
\usepackage{mathdots}  
\newcommand\bbb{\beta}
\newcommand\ccc{\gamma}
\newcommand\aaa{\alpha}
\newcommand\ddd{\delta}
\newcommand\bn{\boldsymbol{n}}
\newcommand\vn{(n_1\,n_2\,n_3\,n_4)}

\newcommand{\thin}{\hspace{0.1em}\rule{0.7pt}{0.8em}\hspace{0.1em}}
\newcommand{\thick}{\hspace{0.1em}\rule{1.5pt}{0.8em}\hspace{0.1em}}

\newtheorem{theorem}{Theorem}
\newtheorem{lemma}[theorem]{Lemma}
\newtheorem{remark}[theorem]{Remark}

\newtheorem{proposition}[theorem]{Proposition}
\newtheorem*{theorem*}{Theorem}

\theoremstyle{definition}
\newtheorem*{definition*}{Definition}

\numberwithin{equation}{section}

\begin{document}

\date{}
\maketitle

\begin{abstract}
	Edge-to-edge tilings of the sphere by congruent quadrilaterals are completely classified in a series of three papers. This last one classifies the case of $a^3b$-quadrilaterals with some irrational angle: there are a sequence of $1$-parameter families of quadrilaterals admitting $2$-layer earth map tilings together with their basic flip modifications under extra condition, and $5$ sporadic quadrilaterals each admitting a special tiling. 
	
	A summary of the full classification is presented in the end.  
	
	{\it Keywords}: 
	spherical tiling, quadrilateral, classification, earth map tiling, irrational angle. 
	
	{\it2000 MR Subject Classification} 52C20, 05B45
\end{abstract}

	\section{Introduction}
 	In an edge-to-edge tiling of the sphere by congruent quadrilaterals, the tile can only have four edge arrangements \cite{ua2,lpwx}:  $a^2bc$, $a^2b^2$, $a^3b$, $a^4$. Akama and Sakano classified tilings for $a^2b^2$ and $a^4$ in  \cite{sa} via the list of triangular tilings in  \cite{ua}. Tilings for $a^2bc$ are classified in the first paper \cite{lpwx} of this series via the methods in \cite{wy1,wy2,awy,wy3} developed for pentagonal tilings. Tilings for $a^3b$ with rational angles are classified in the second paper \cite{lw} of the series by solving some trigonometric Diophantine equations. In this last paper of the series, we classify tilings for  $a^3b$ with some irrational angle, and therefore obtain the full classification. 
 
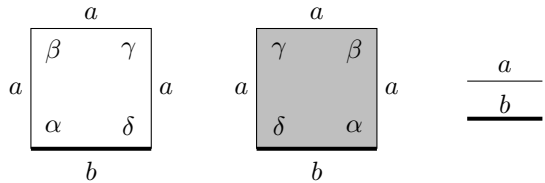
\begin{figure}[htp]
	\centering
	\begin{tikzpicture}
		
		\fill[gray!50]  (-0.8+3,-0.8) -- (-0.8+3,0.8) -- (0.8+3,0.8) -- (0.8+3,-0.8);
		
		\foreach \a in {0,1}
		{
			\begin{scope}[xshift=3*\a cm]
				\draw
				(-0.8,-0.8) -- (-0.8,0.8) -- (0.8,0.8) -- (0.8,-0.8);

				\draw[line width=1.5]
				(-0.8,-0.8) -- (0.8,-0.8);	    	\end{scope}}   
		
		\draw  (4.5+0.5,0.1)--(5.5+0.5,0.1);

		\draw[line width=1.5] (4.5+0.5,-0.4)--(5.5+0.5,-0.4);
		
		\node at (-0.5,0.5) {\small $\bbb$};
		\node at (0.5,0.5) {\small $\ccc$};
		\node at (-0.5,-0.5) {\small $\aaa$};
		\node at (0.5,-0.5) {\small $\ddd$};
		
		\foreach \a in {0,1,2}
		\node at (90*\a:1) {\small $a$};
		\node at (0,-1.1) {\small $b$};
		
		\node at (3,1) {\small $a$};\node at (2,0) {\small $a$};\node at (4,0) {\small $a$};
		\node at (3,-1.1) {\small $b$};
		
		\node at (2+3.5,0.3) {\small $a$};
		
		\node at (2+3.5,-0.2) {\small $b$};

		\node at (-0.5+3,0.5) {\small $\ccc$};
		\node at (0.5+3,0.5) {\small $\bbb$};
		\node at (-0.5+3,-0.5) {\small $\ddd$};
		\node at (0.5+3,-0.5) {\small $\aaa$};
		
	\end{tikzpicture} 
	\caption{Quadrilateral of type $a^3b$.}
	\label{quad}
\end{figure}

    Type $a^3b$ quadrilateral is given by the same Fig.\,\ref{quad} as in the second paper \cite{lw} of our series: one quadrilateral is the mirror image or flip of the other, and the shading can help present the tiling without  indicating all angles. Throughout this paper, an $a^3b$-tiling is always an edge-to-edge tiling of the sphere by congruent simple quadrilaterals in Fig.\,\ref{quad}, such that all vertices have degree $\ge 3$.

\begin{figure}[htp]
	\centering
	\includegraphics[scale=0.168]{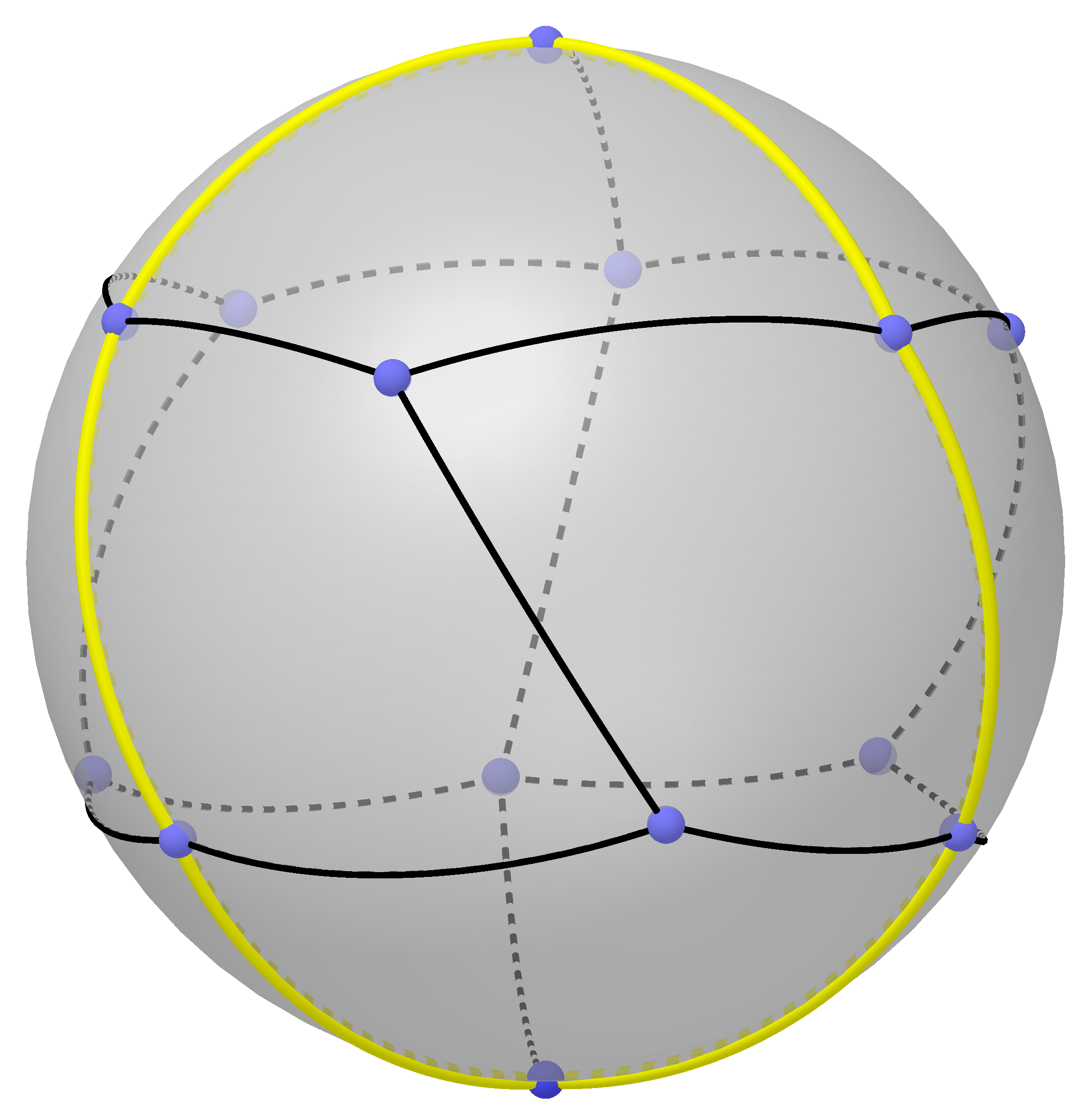}
	\includegraphics[scale=0.146]{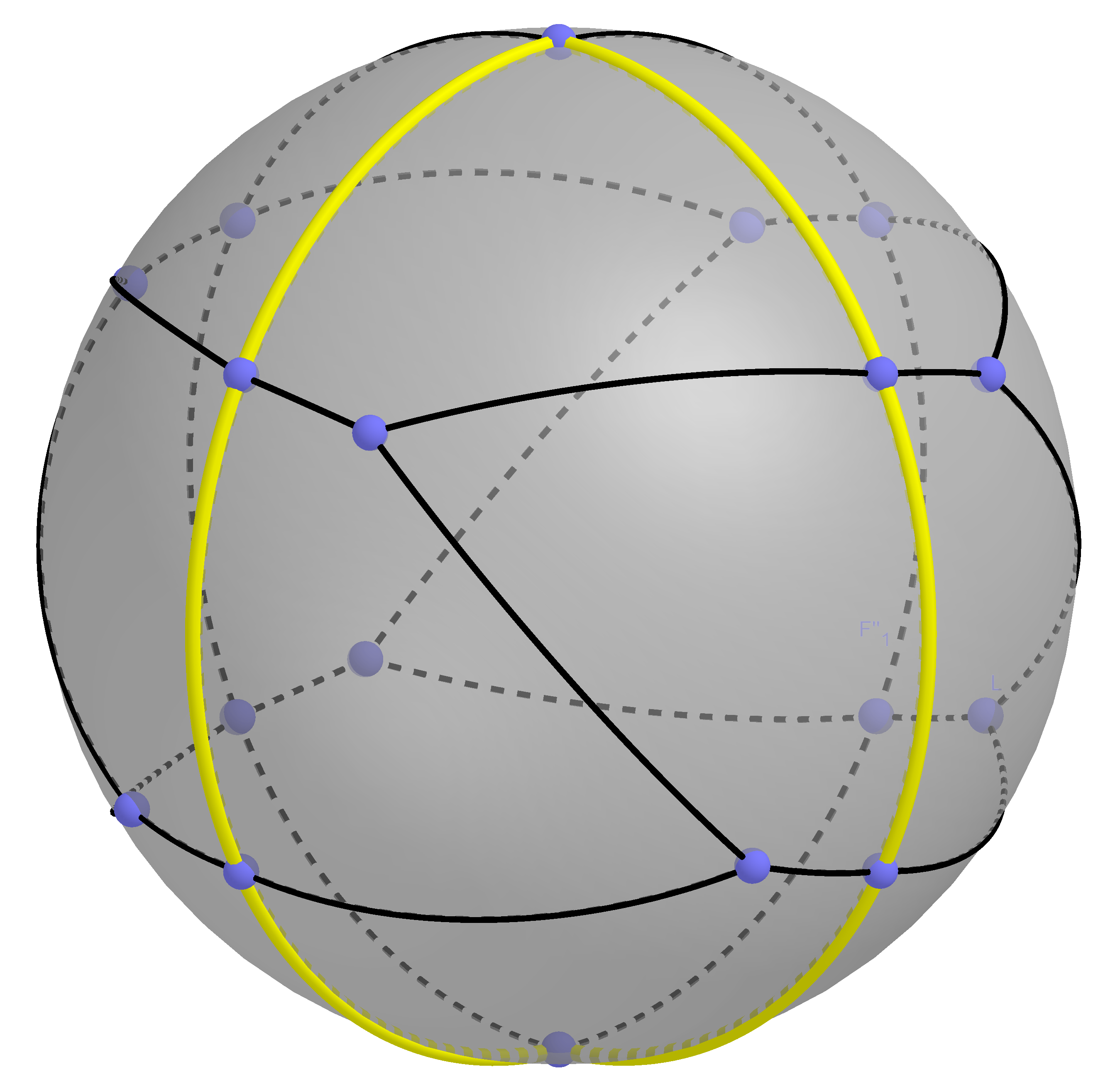}				
	\includegraphics[scale=0.15]{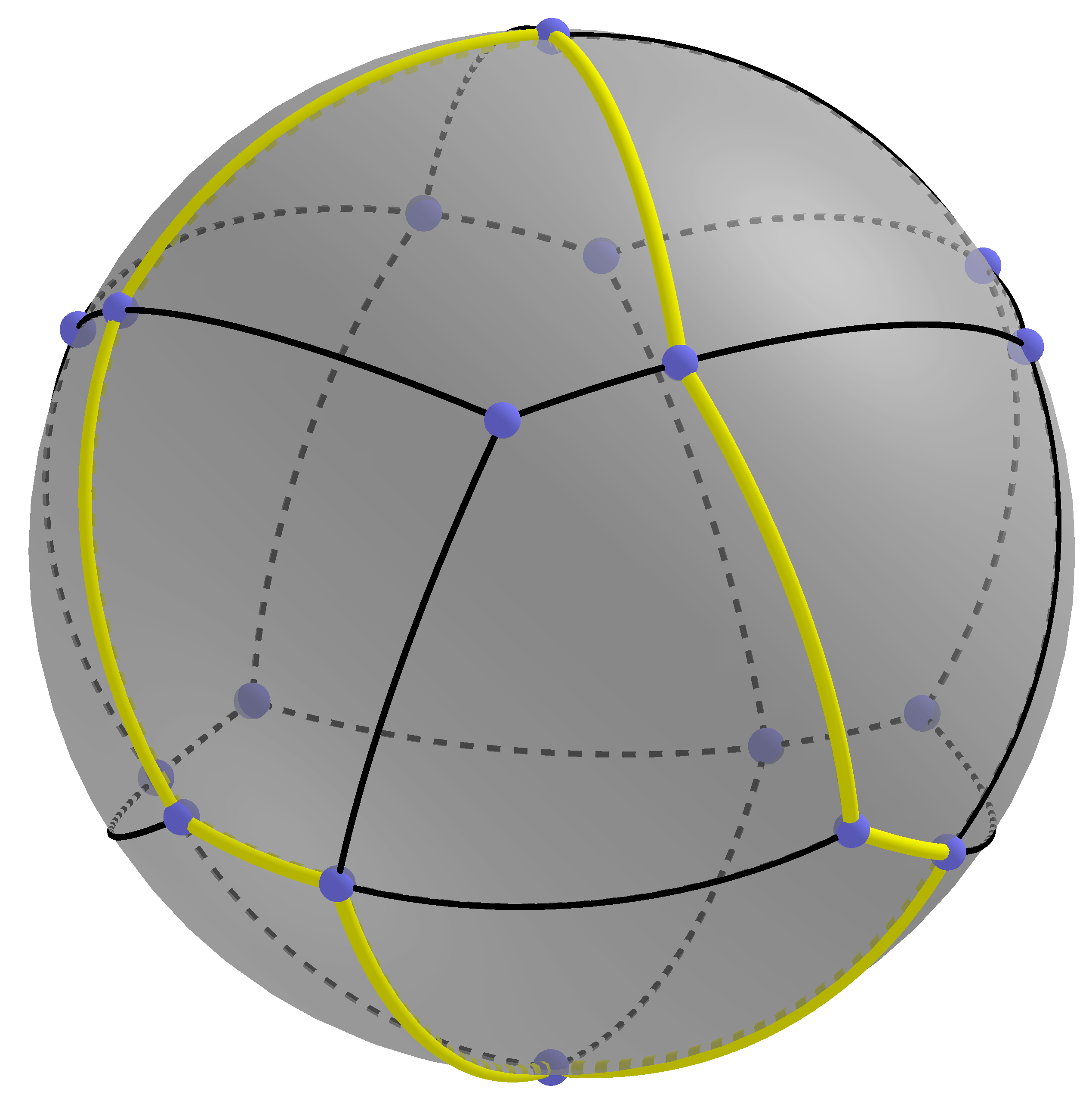}  
	\includegraphics[scale=0.15]{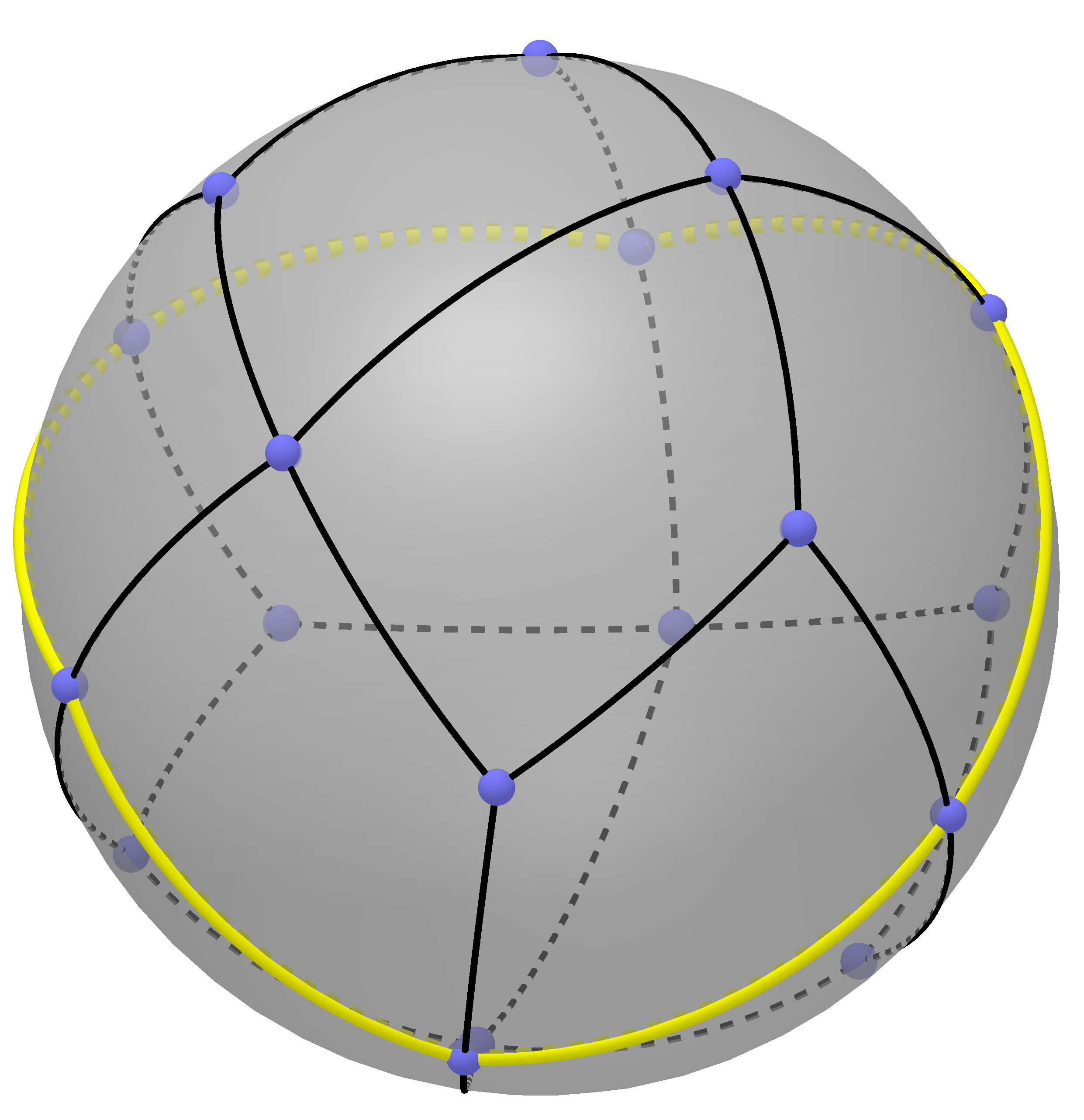}

	\caption{Four sporadic tilings with $f=12,16,16,16$.} 
	\label{a2}	
\end{figure}

\begin{theorem*}
	All $a^3b$-tilings with some irrational angle are the following:
	\begin{enumerate}
		\item Five sporadic tilings in Table \ref{tab-1}  (each has a unique quadrilateral), also shown in Fig.\,\ref{a2} and Fig.\,\ref{case a3 d4}.
		\item A sequence of $1$-parameter families of quadrilaterals admitting $2$-layer earth map tilings $T(f\aaa\bbb\ddd,2\ccc^{\frac{f}{2}})$ for every even integer $f\ge6$ in Table \ref{tab-2}, together with their basic flip modifications in Table \ref{tab-3} when $\bbb$ is an integer multiple of $\ccc$. All quadrilaterals in these families have some irrational angles except the $12$ sporadic and $3$ infinite sequences of cases listed in \cite{lw}. 		
	\end{enumerate}
\end{theorem*} 

 \begin{table}[htp]
	\begin{center}
	
		\bgroup
		\def\arraystretch{1.5}
		  \begin{tabular}[htp]{  c|c|c }

				$f$ & $(\aaa,\bbb,\ccc,\ddd),a,b$ & all vertices  \\ \hline
				
				$12$ &  \parbox[c][3.0cm][c]{8cm}{$(\aaa,2-2\aaa,\frac23,-\frac13+\aaa)\\ \aaa=1-\arcsin{\frac{\sqrt6}{4}}\approx 0.7902$
				\\$a=\arccos{\frac{2\sqrt5-3}{3}}\approx0.3367$
				\\$b=\arccos{(3\sqrt5-6)}\approx0.2495$} & $T(6\aaa^2\bbb,2\ccc^3,6\bbb\ccc\ddd^2)$  \\
				\hline
				$16$& \parbox[c][3.0cm][c]{8cm}{$(\aaa,2-2\aaa,\frac12,-\frac14+\aaa)\\ \aaa=1-\arcsin{\frac{\sqrt{3\sqrt{10}-3\sqrt5-3\sqrt2+15}}{6}}\approx 0.7898$\\
					$a=\arccos{\frac{\sqrt{10}+\sqrt{5}-\sqrt2-3}{2}}\approx 0.3362$\\
					$b=\arccos{\frac{(27\sqrt5-43)\sqrt2+23\sqrt5-27}{(196\sqrt5-420)\sqrt2-267\sqrt5+623}}\approx 0.1052$}  &  $T(8\aaa^2\bbb,8\bbb\ccc\ddd^2,2\ccc^4)$  \\
				\hline
				$16$ & \parbox[c][3.0cm][c]{8cm}{$(\aaa,\frac12,1-\aaa,\frac34)\\ \aaa=1-\arcsin{\frac{\sqrt{2\sqrt2+1}}{2}}\approx 0.5664$\\
					$a=\arccos{\frac{1}{\sqrt{2\sqrt2+1}}}\approx 0.3292$\\
					$b=\arccos{\frac{7}{(2\sqrt2+1)^{\frac32}}}\approx 0.1158$} 
				&  $T(8\bbb\ddd^2,8\aaa^2\ccc^2,2\bbb^4)$\\
				\hline
				$16$ & \parbox[c][3.0cm][c]{8cm}{$(\aaa,\frac12,\frac34,1-\aaa)\\ \aaa=1-\arcsin{\tfrac{\sqrt{10+4\sqrt2}}{\sqrt{17}}}\approx 0.5906$\\
					$a=\tfrac14$\\
					$b=\arccos{\tfrac{2\sqrt 2-1}{4}}\approx 0.3488$} &  $T(8\bbb\ccc^2,6\aaa^2\ddd^2,4\aaa\bbb^2\ddd)$\\
				\hline
				$24$	 & \parbox[c][3.0cm][c]{8cm}{$(\aaa,\frac23,1-\aaa,\frac12)\\ \aaa=\arcsin{\tfrac{\sqrt{4+\sqrt3}}{\sqrt6}}\approx 0.4322$\\
					$a=\arcsin{\tfrac{\sqrt2}{\sqrt{4+\sqrt3}}}\approx 0.2011$\\
					$b=\tfrac12-a\approx 0.2988$} 
				&  $T(8\bbb^3,12\aaa^2\ccc^2,6\ddd^4)$\\	
				\hline
				
		\end{tabular}
		\egroup
			\caption{Five sporadic quadrilaterals and their unique tilings.}\label{tab-1}
	\end{center}
\end{table}

\begin{table}[htp]
	\begin{center}
		\bgroup
		\resizebox{\textwidth}{21mm}{  \begin{tabular}[htp]{ccccc}				
				\multicolumn{4}{c|}{$\aaa,\bbb,\ccc,\ddd,a,b$} & vertices and tilings \\ 
				\hline
				\multicolumn{4}{c|}{ $\,\,\,\,\,\,f=6: \, \forall\,\bbb\in(\frac13,\frac32)\backslash\{\frac23\}$; $\,\,\,\,\,\,\,f\ge8:\, \forall\,\bbb\in(\frac12,\frac32)\backslash\{1-\frac{2}{f}\} \,\,\,\,\,\,$   }& \multirow{7}{*}{$T(f\aaa\bbb\ddd,2\ccc^{\frac{f}{2}})$}\\
				\multicolumn{4}{c|}{$\cos a=-\frac{\cos \bbb}{1-\cos \bbb}$,$\quad\ccc=\frac 4f$, $\quad\ddd=2-\aaa-\bbb$ }&\multirow{7}{*}{($\forall$ even $f\ge6$)}\\
				\multicolumn{4}{c|}{$\cos b= \cos^3 a(1-\cos \bbb)(1-\cos \frac{4\pi}{f})-\cos^2 a\sin \bbb\sin \frac{4\pi}{f}+$}&\\
				\multicolumn{4}{c|}{$\,\,\,\,\,\quad\quad\,\,\cos a(\cos \bbb+\cos \frac{4\pi}{f}-\cos \bbb\cos \frac{4\pi}{f}) +\sin\bbb\sin\frac{4\pi}{f}$}&\\
				\multicolumn{4}{c|}{$\tan \aaa=\frac{\cos \! \bbb \left(\left(\cos \! \bbb -1\right) \sin \! \frac{4 \pi}{f}+2 \sin \! \bbb  \sin^2 \! (\frac{2 \pi}{f})\right)}{\left(\sin \! \frac{4 \pi}{f} \sin \! \bbb -1\right) \left(\cos \! \bbb -1\right)-2 \cos^2 \! (\bbb) \sin^{2} \! (\frac{2 \pi}{f})}$, $\cos\aaa=\frac{\cos\frac{2\pi}{f}-\cos\frac a2 \cos\frac b2}{\sin\frac a2 \sin\frac b2}$}& \\
				\multicolumn{4}{l|}{$f=6$: $\frac13\le a<1$. If $\frac13<\bbb<\frac43$, then $b<1$; If $\frac43\le\bbb<\frac32$, then $1\le b<\frac76$.}&\\
				\multicolumn{4}{l|}{$f\ge8$: $\frac13\le a<\frac12$, $b<\frac12+\frac{4}{f}$.  $\lim_{f\to\infty}b=a,\,\lim_{f\to\infty}\aaa=0,\,\lim_{f\to\infty}\ddd=2-\bbb$}&\\
				\hline		
		\end{tabular}}
		\egroup
		\caption{A sequence of $1$-parameter families of $2$-layer earth map $a^3b$-tilings.}\label{tab-2}
	\end{center}
\end{table}
\begin{table}[htp]
	\begin{center}
		
		\bgroup
		\resizebox{\textwidth}{38mm}{  \begin{tabular}[htp]{c|c}				
                $\aaa+\ddd=m\ccc$& \multirow{2}{*}{Flip modifications only if $f\ge8$}\\
                $\bbb=(\frac{f}{2}-m)\ccc$&\\
                \hline
		    	\multirow{3}{*}{$ f/8<m\le f/6$}		    	
				&$(f-2)\aaa\bbb\ddd,2\aaa\ccc^{\frac{f}{2}-m}\ddd,2\bbb\ccc^m$\\				
				&$(f-4)\aaa\bbb\ddd,2\aaa^2\ccc^{\frac{f}{2}-2m}\ddd^2,4\bbb\ccc^m$: $\lfloor \frac{f-4m+4}{4} \rfloor$ tilings\\
				&$(f-6)\aaa\bbb\ddd,2\aaa^3\ccc^{\frac{f}{2}-3m}\ddd^3,6\bbb\ccc^m$: $\lfloor \frac{f-6m+4}{4} \rfloor +\langle \frac{(f-6m)^2}{48}\rangle$ tilings\\
				\hline
				\multirow{2}{*}{$ f/6<m< f/4$} 
				&$(f-2)\aaa\bbb\ddd,2\aaa\ccc^{\frac{f}{2}-m}\ddd,2\bbb\ccc^m$\\
				&$(f-4)\aaa\bbb\ddd,2\aaa^2\ccc^{\frac{f}{2}-2m}\ddd^2,4\bbb\ccc^m$: $\lfloor \frac{f-4m+4}{4}\rfloor$ tilings\\
				\hline
				\multirow{2}{*}{$m= f/4$}
				&$(f-2)\aaa\bbb\ddd,2\aaa\ccc^{m}\ddd,2\bbb\ccc^{m}$\\
				\multirow{2}{*}{$(\bbb=1)$}&$(f-4)\aaa\bbb\ddd,2\aaa^2\ddd^2,4\bbb\ccc^{m}$\\
				&If $m=(f+2)/4$, then $a=b$ and it is of type $a^4$ in Section \ref{conclusion}\\
				\hline
				\multirow{2}{*}{$(f+4)/4\le m<f/3$}
				&$(f-2)\aaa\bbb\ddd,2\aaa\ccc^{\frac{f}{2}-m}\ddd,2\bbb\ccc^m$\\
				\multirow{2}{*}{$(f>12)$}&$(f-4)\aaa\bbb\ddd,4\aaa\ccc^{\frac{f}{2}-m}\ddd,2\bbb^2\ccc^{2m-\frac{f}{2}}$: $\lfloor \frac{4m-f+4}{4}\rfloor$ tilings\\
				&If $m=(f+4)/4$, then $(\aaa,\bbb,\ccc,\ddd)=(\frac4f,1-\frac4f,\frac4f,1)$\\
				\hline
				\multirow{3}{*}{$f/3\le m< 3f/8$} 
				&$(f-2)\aaa\bbb\ddd,2\aaa\ccc^{\frac{f}{2}-m}\ddd,2\bbb\ccc^m$\\
				&$(f-4)\aaa\bbb\ddd,4\aaa\ccc^{\frac{f}{2}-m}\ddd,2\bbb^2\ccc^{2m-\frac{f}{2}}$: $\lfloor \frac{4m-f+4}{4}\rfloor$ tilings\\
				&$(f-6)\aaa\bbb\ddd,6\aaa\ccc^{\frac{f}{2}-m}\ddd,2\bbb^3\ccc^{3m-f}$: $\lfloor \frac{3m-f+2}{2}\rfloor+\langle  \frac{(3m-f)^2}{12}\rangle$ tilings\\
				\hline			
		\end{tabular}}
		\egroup
		\caption{Basic flip modifications if $\bbb$ is an integer multiple of $\ccc$ in Table \ref{tab-2}.} \label{tab-3}
	\end{center}
\end{table}
\begin{table*}[htp]          
	\centering     

	~\\ 
	\resizebox{\textwidth}{9mm}{\begin{tabular}{c|ccccccccccc}	 
			
			$f$ &8&$24k-16$&$24k-14$&$24k-12,24k-10$&$24k-8$&$18$&$24k-6$&$24k-4$&$24k-2,24k,24k+2$&$24k+4$&$24k+6$ \\
			\hline
			$k$&&$\ge2$&$\ge1$&$\ge1$&$\ge1$&&$\ge2$&$\ge1$&$\ge1$&$\ge1$&$\ge1$\\
			\hline
			$\mathcal{Q}_{1}(f)$&1&$6k-5$&$6k-5$&$6k-3$&$6k-3$&$3$&$6k-3$&$6k-1$&$6k-1$&$6k+1$&$6k+1$  \\
			\hline			
			$\mathcal{Q}_{2}(f)$&0&$2$&$1$&$1$&$2$&1&$0$&$2$&$1$&$2$&$0$ \\	
			\hline	
			$\mathcal{Q}_{3}(f)$&1&$6k-7$&$6k-6$&$6k-4$&$6k-5$&$2$&$6k-3$&$6k-3$&$6k-2$&$6k-1$&$6k+1$  \\
			\hline	
	\end{tabular}}
	\caption{Counting $a^3b$-quadrilaterals (rational and general) admitting flips.}\label{Tab-1.3}        
\end{table*}

 The numerical data and exact formulas are listed in Table \ref{tab-1} and \ref{tab-2}, where the angles and edge lengths are expressed in units of $\pi$,  and the last column counts all vertices and also all tilings when they are not uniquely determined by the vertices. A rational fraction, such as $\ccc=\frac{2}{3}$, means the precise value $\frac{2\pi}{3}$. A decimal expression, such as $b\approx0.2495$, means an approximate value  $0.2495\pi < b < 0.2496\pi$. We put $\pi$ back in any trigonometric functions to avoid confusion.  

\begin{figure}[htp]
	\centering
	\includegraphics[scale=0.16]{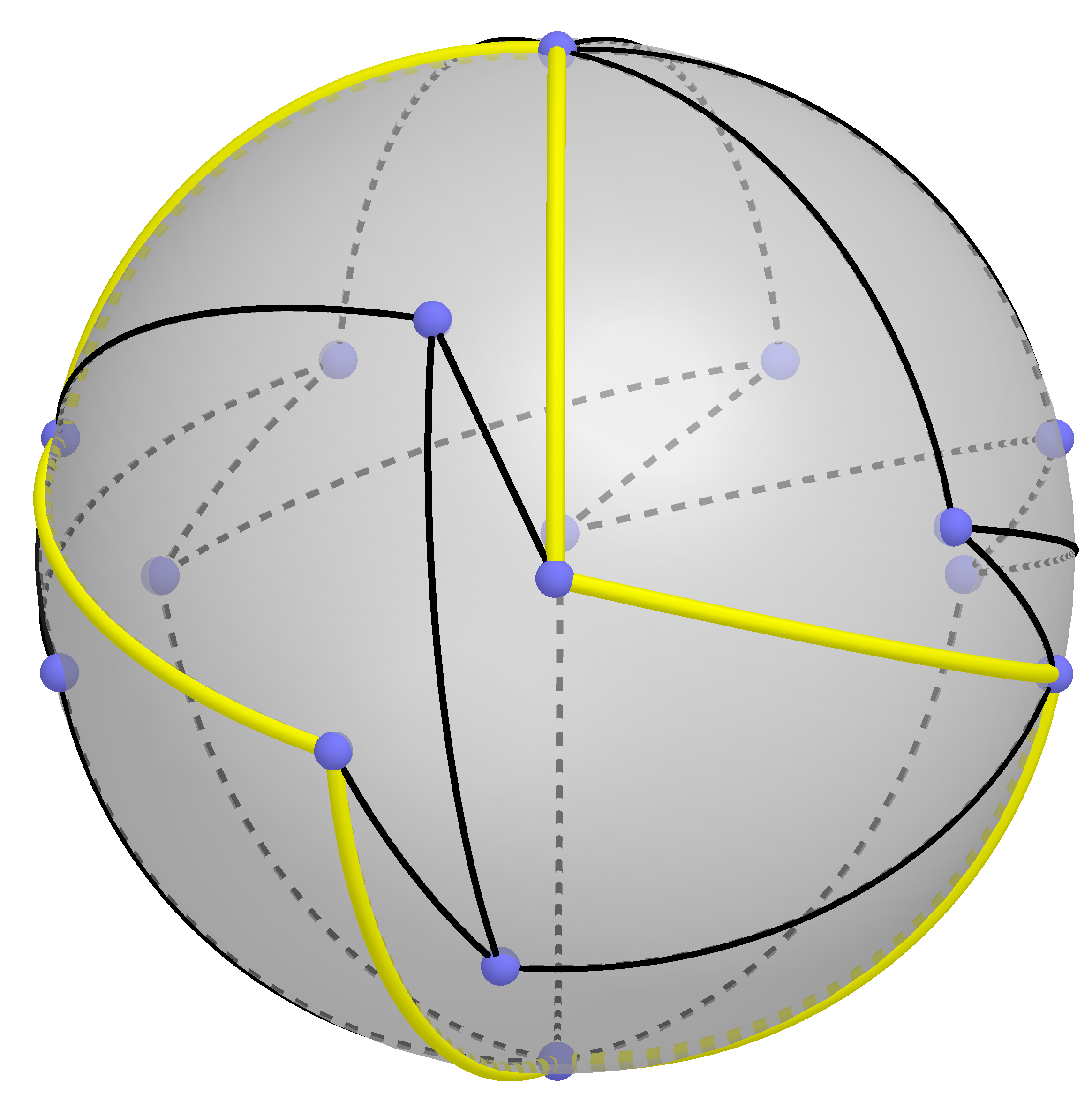}				
	\includegraphics[scale=0.155]{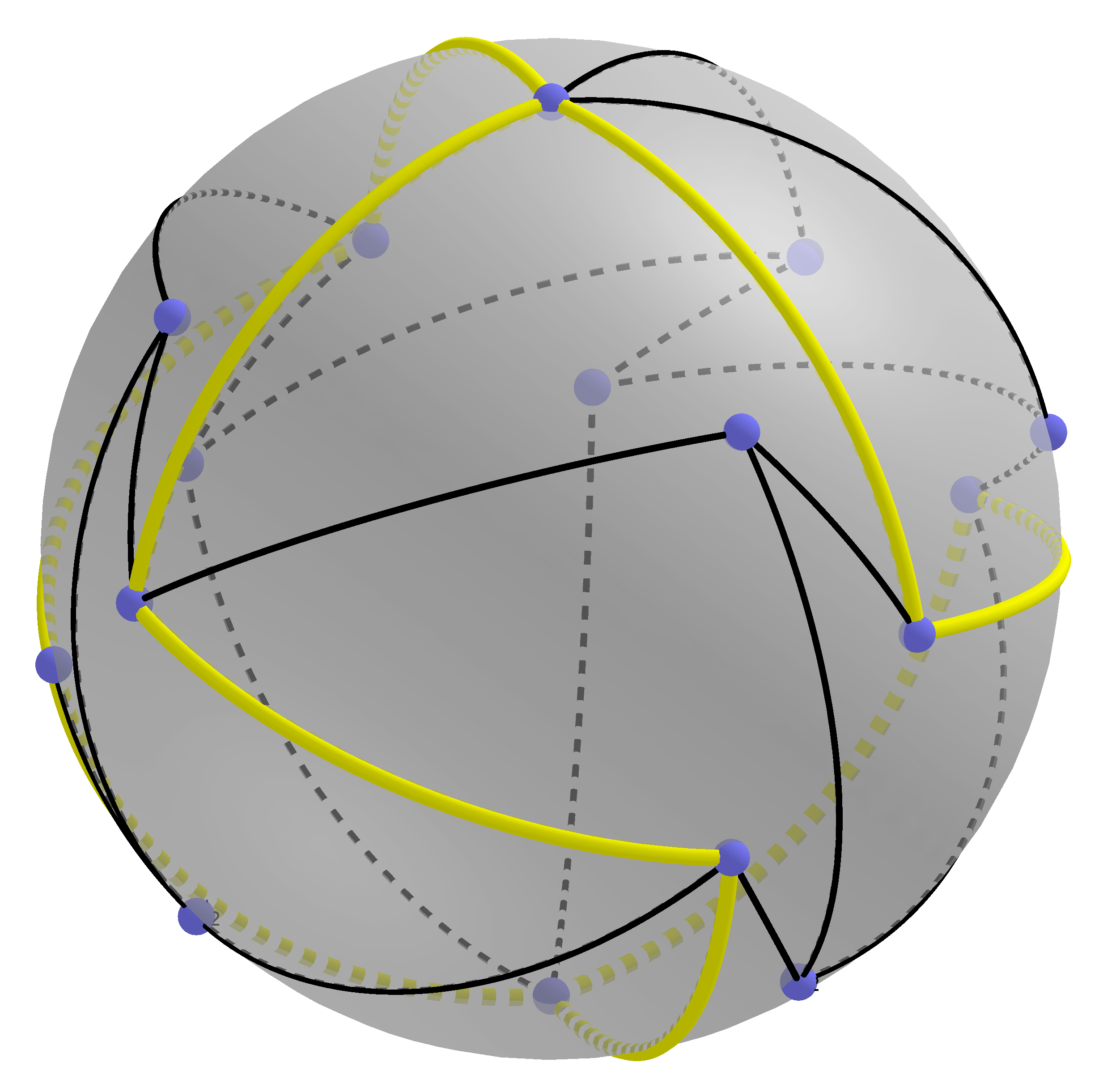}  
	\includegraphics[scale=0.16]{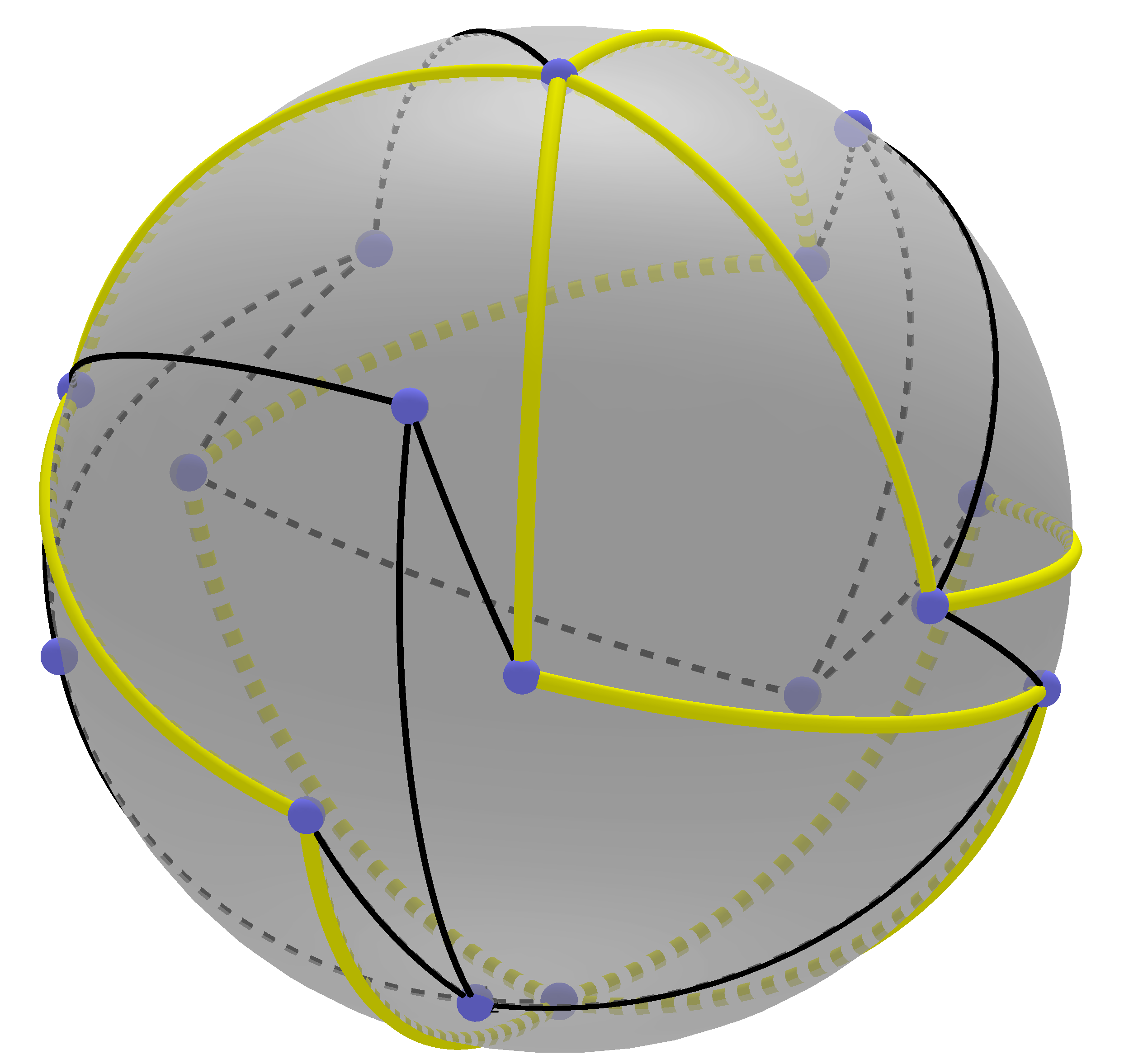}
	\includegraphics[scale=0.16]{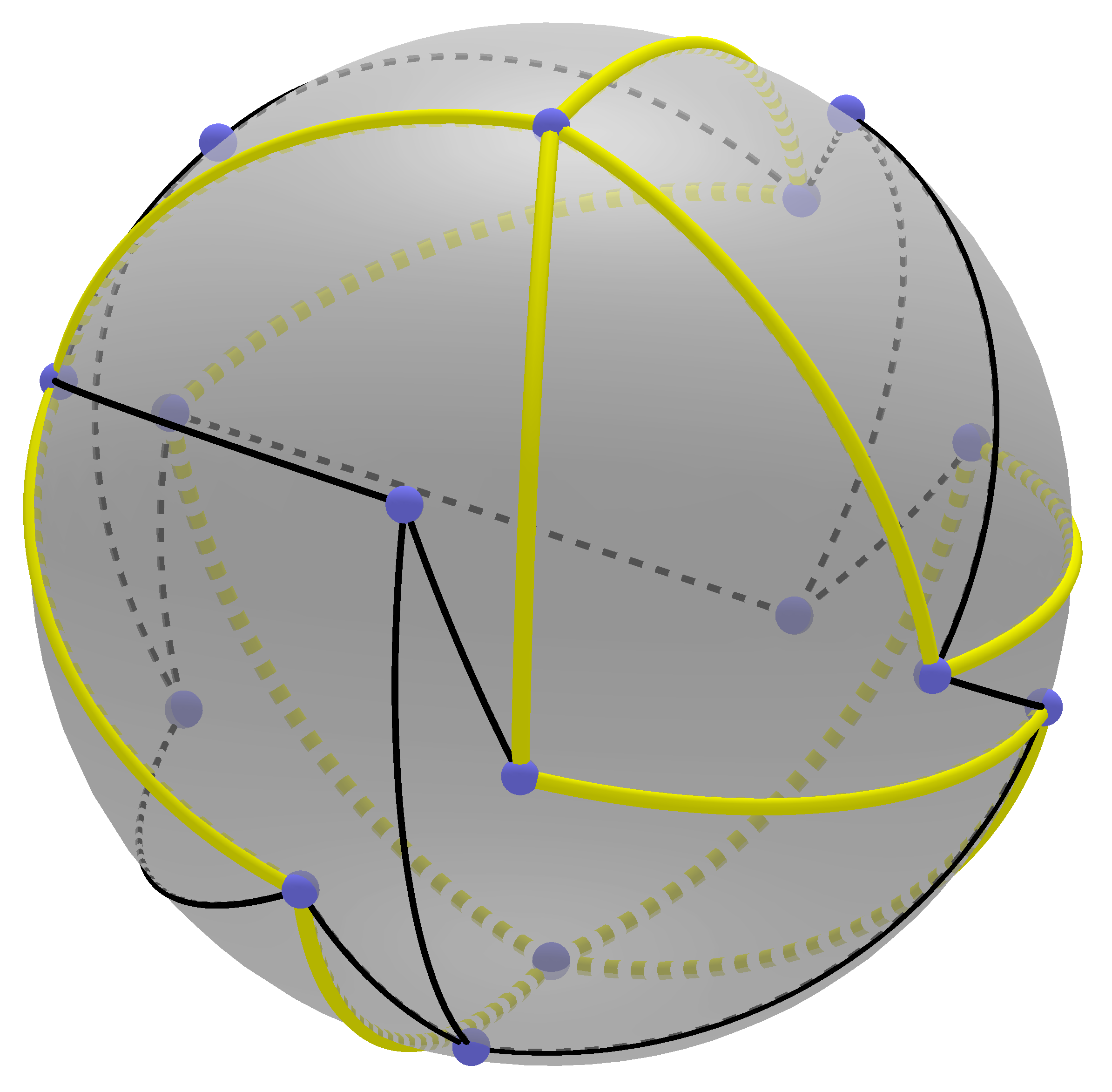}
	
	\caption{Four flips for $f=14$ in the last case of Table \ref{tab-3}.} 
	\label{a3}	
\end{figure}

 There are two types of basic flip modifications in Fig.\,\ref{flip1}, depending on $\bbb<1$ or $\bbb\ge1$. 
 Four pictures of Fig.\,\ref{a3} illustrate that $1,2$ or $3$ basic flips can be applied simultaneously on special $2$-layer earth map $a^3b$-tilings in many different ways depending on the space between these flips. 
 
 Table \ref{Tab-1.3} counts the total number $\mathcal{Q}_{1}$ of $a^3b$-quadrilaterals admitting any flips in Table \ref{tab-3} ($\mathcal{Q}_{2}$ for rational case, $\mathcal{Q}_{3}=\mathcal{Q}_{1}-\mathcal{Q}_{2}$ for general case).  
 In the last column of Table \ref{tab-3}, the total number of different tilings are counted using unordered integer partitions, where the function $\langle x \rangle$ gives the closest integer to $x$. The details for Table \ref{Tab-1.3} and \ref{tab-3} are discussed in the end of Section \ref{sec-111}.

 \subsubsection*{Outline of the paper}
 The second of our series \cite{lw} shows that most rational $a^3b$-tilings are $2$-layer earth map tilings and their rearrangements. This phenomena turns out to be true in general. After reviewing the basic results from \cite{lpwx,lw} in Section \ref{basic_facts}, we characterize the quadrilaterals in $2$-layer earth map $a^3b$-tilings by some special degree $3$ vertex, and find all their tilings in Section \ref{sec-111}. Sections \ref{sec-two3} and \ref{sec-one3} handle all other possibilities for degree $3$ vertices and find five sporadic tilings. A summary of all edge-to-edge spherical tilings by congruent quadrilaterals is presented in the last section.

\section{Basic Facts}
\label{basic_facts}

We will always express angles in $\pi$ radians for simplicity. Therefore, the sum of all angles ({\em angle sum}) at a vertex is $2$. Most of the following lemmas have been proved in \cite{lpwx,lw}.

Let $v,e,f$ be the numbers of vertices, edges, and tiles in a quadrilateral tiling. Let $v_k$ be the number of vertices of degree $k$. Euler's formula $v-e+f=2$ implies (see \cite{lpwx})
\begin{align}
	f&=6+ \sum_{k=4}^{\infty}(k-3)v_k
	=6+v_4+2v_5+3v_6+\cdots, \label{vcountf} \\
	v_3 &=8+\sum_{k=5}^{\infty}(k-4)v_k=8+v_5+2v_6+3v_7+\cdots. \label{vcountv}
\end{align}

Therefore, $f\ge 6$ and $v_3 \ge 8$.

\begin{lemma}[{\cite[Lemma 2]{lpwx}}] \label{anglesum} 
	If all tiles in a tiling of the sphere by $f$ quadrilaterals have the same four angles $\aaa,\bbb,\ccc,\ddd$, then 
	\[
	\aaa+\bbb+\ccc+\ddd=2+\frac{4}{f}, 
	\]
	ranging in $(2,\frac83]$. In particular no vertex contains all four angles.  
\end{lemma}

\begin{lemma}[{\cite[Lemma 3]{wy2}}]\label{geometry1}
	If the quadrilateral in Fig.\,\ref{quad} is simple, then  $\bbb<\ccc$ is equivalent to $\aaa>\ddd$.
\end{lemma}

\begin{lemma}[{\cite[Lemma 3]{lw}}]\label{geometry3}
	If the quadrilateral in Fig.\,\ref{quad} is simple, then $\bbb=\ddd$ if and only if $\aaa=1$. Furthermore, if it is convex with all angles $<1$, then $\bbb>\ddd$ is equivalent to $\aaa<\ccc$, and  $\bbb<\ddd$ is equivalent to $\aaa>\ccc$.    
\end{lemma}

\begin{lemma}[{\cite[Lemma 4]{lw}}]\label{geometry4}
	If the quadrilateral in Fig.\,\ref{quad} is simple, and $\ddd\le 1$, then $2\aaa+\bbb>1$ and $\bbb+2\ccc>1$ .
\end{lemma}

\begin{lemma}[Parity Lemma, {\cite[Lemma 10]{wy2}}] \label{parity}
	In an $a^3b$-tiling, the total number of $ab$-angles $\aaa$ and $\ddd$ at any vertex is even.
\end{lemma}

In a tiling of the sphere by $f$ congruent tiles, each angle of the tile appears $f$ times in total. If one vertex has more $\aaa$ than $\bbb$, there must exist another vertex with more $\bbb$ than $\aaa$.
Such global counting induces many interesting and useful results. 

\begin{lemma}[Balance Lemma, {\cite[Lemma 6]{lw}} ]\label{balance}
	If either $\aaa^2\cdots$ or $\ddd^2\cdots$ is not a vertex, then any vertex either has no $\aaa,\ddd$, or is of the form $\aaa\ddd\cdots$ with no more $\aaa,\ddd$ in the remainder.
\end{lemma}

The very useful tool \textit{adjacent angle deduction} (abbreviated as AAD) has been introduced in \cite[Section 2.5]{wy1}. The following is \cite[Lemma 10]{wy1}.

\begin{lemma} \label{aadlemma}
	Suppose $\lambda$ and $\mu$ are the two angles adjacent to $\theta$ in a quadrilateral.
	\begin{itemize}
		\item If $\lambda\thin\lambda\cdots$ is not a vertex, then $\theta^n$ has the unique AAD $\thin^{\lambda}\theta^{\mu}\thin^{\lambda}\theta^{\mu}\thin^{\lambda}\theta^{\mu}\thin\cdots$.
		\item If $n$ is odd, then we have the AAD $\thin^{\lambda}\theta^{\mu}\thin^{\lambda}\theta^{\mu}\thin$ at $\theta^n$.
	\end{itemize}
\end{lemma}

\begin{lemma}[{\cite[Lemma 8]{lw}}]\label{proposition-1}
	There is no tiling of the sphere by congruent quadrilaterals with two angles $\ge1$. 
\end{lemma}

\begin{lemma}[{\cite[Proposition 3]{ua2}}]\label{relation}
	Any convex spherical $a^3b$-quadrangle with all angles $<1$ satisfies 
	\[ 	\aaa+\ddd<1+\bbb, \quad	\aaa+\ddd<1+\ccc, \quad	\aaa+\bbb<1+\ddd, \quad	\ccc+\ddd<1+\aaa.\]
\end{lemma}

\begin{lemma}[{\cite[Lemma 13]{lw}}]\label{calculate-1}	
	For $a^3b$-quadrilaterals, the following equations hold:  
\end{lemma}
\begin{equation}
	\begin{split}
		\cos b=& \cos^3 a(1-\cos \bbb)(1-\cos \ccc)-\cos^2 a\sin \bbb\sin \ccc+\\ &\cos a(\cos \bbb+\cos \ccc-\cos \bbb\cos \ccc) +\sin\bbb\sin\ccc;  \label{4}
	\end{split}
\end{equation} 
\begin{equation}
	\cos a=\frac{\sin \aaa+\cos\ddd\sin\ccc}{2\sin\ddd\sin^2\frac{\ccc}{2}}=\frac{\sin \ddd+\cos\aaa\sin\bbb}{2\sin\aaa\sin^2\frac{\bbb}{2}} \quad (\aaa,\ddd\neq1);   \label{4-2}
\end{equation}
\begin{align}
	\sin(\aaa-\frac{\ccc}{2})\sin\frac{\bbb}{2}=\sin\frac{\ccc}{2}\sin(\ddd-\frac{\bbb}{2}), \label{4-7}
\end{align}
\begin{align}
	\text{or}\,\,\,	\sin(\aaa+\frac{\ccc}{2})\sin\frac{\bbb}{2}=-\sin\frac{\ccc}{2}\sin(\ddd+\frac{\bbb}{2}).  \label{4-8}
\end{align} 

When $\aaa,\bbb,\ccc,\ddd<1$, the equations \eqref{4-2}, \eqref{4-7} were proved in \cite[Theorem 2.1]{coolsaet} while the equation \eqref{4-8} was dismissed.

\begin{lemma}\label{geometry2}
	In an $a^3b$-tiling, then we cannot have three different degree $3$ vertices with only one $b$-edge in a tiling by the quadrilateral. Moreover, the following are the only possible combinations of two degree $3$ vertices with only one $b$-edge:
	\begin{enumerate}
		\item $\aaa\bbb\ddd,\ccc\ddd^2$.
		\item $\aaa^2\bbb,\aaa\ccc\ddd$.
		\item $\aaa^2\bbb,\ccc\ddd^2$.
	\end{enumerate}
\end{lemma}
The proof is exactly the same as Lemma $5$ in \cite{wy2}.

\vspace{9pt}

Following Rao \cite{rao}, a vertex $\aaa^{n_1}\bbb^{n_2}\ccc^{n_3}\ddd^{n_4}$ can be efficiently represented by its vector type $\bn=(n_1\, n_2 \, n_3 \, n_4)\in \mathbb{N}^4$. Lemma \ref{anglesum} implies that one of $n_j$ must be $0$. We will use both representations interchangeably afterwards for convenience. 

Since the tilings by rational $a^3b$-quadrilaterals have been classified in \cite{lw}, this paper handles the case with some irrational angle, which imposes strong constraints on all possible vertices. 

\begin{lemma}[Irrational Angle Lemma]\label{irrational}
	
	Given two different vertices $\boldsymbol{l},\boldsymbol{m}\in \mathbb{N}^4$ in a spherical quadrilateral tiling with some irrational angle and $f$ tiles, if $\boldsymbol{l},\boldsymbol{m}$ and $\boldsymbol{u}:=(1\, 1 \, 1 \, 1)$ are linearly independent, then any other vertex $\bn$ are
	\[
	\bn\in\mathbb{N}^4\cap\mathrm{Plane}\left\{\boldsymbol{l},\boldsymbol{m},\frac{f}{f+2}\boldsymbol{u}\right\}. 
	\]
	In particular, the square matrix of $\boldsymbol{l},\boldsymbol{m}, \boldsymbol{u},\bn$ has determinant $0$. 
\end{lemma}
\begin{proof} 
	By Lemma \ref{anglesum} and the vertices $\boldsymbol{l},\boldsymbol{m}, \boldsymbol{n}$, the angles satisfy a linear system of equations with the augmented matrix 
	\[C=\left(\begin{matrix}
		1 & 1 & 1 & 1 & 2+\tfrac4f \\
		l_1 & l_2 & l_3 & l_4 & 2\\
		m_1 & m_2 & m_3 & m_4 & 2\\
		n_1 & n_2 & n_3 & n_4 & 2
	\end{matrix}\right).
	\]
	If $\bn$ is not a linear combination of $\boldsymbol{l},\boldsymbol{m}, \boldsymbol{u}$, then the system has a unique rational solution, a contradiction. Therefore, the last row of $C$ must be some linear combination of the other $3$ rows, which implies that $\bn$ lies in the plane containing the three points $\boldsymbol{l},\boldsymbol{m},\frac{f}{f+2}\boldsymbol{u}$. 
\end{proof}

\section{Tilings with some vertex being $\aaa\bbb\ddd$ or $\aaa\ccc\ddd$}
\label{sec-111}

All rational $a^3b$-tilings have been
discussed in \cite{lw}. We emphasize that all propositions afterwards assume the quadrilateral has some irrational angle. 

\begin{lemma}\label{3-abd}
	In an $a^3b$-tiling, if $\aaa\bbb\ddd$ is a vertex, one of $\aaa^2\ccc^{n-2}$, $\aaa\ccc^{n-2}\ddd$, $\bbb^2\ccc^{n-2}$, $\bbb\ccc^{n-1}$, $\ccc^{n}$ or $\ccc^{n-2}\ddd^2$ must appear for some $n\ge3$. 
\end{lemma}

\begin{proof}
	If the above vertices do not appear for any $n\ge3$, then $\ccc$ appears $\le n-3$ times in any degree $n$ vertex ($\ccc$ never appears in any degree $3$ vertex). Therefore, the total number of $\ccc$ in the tiling $\#\ccc\le v_4+2v_5+3v_6+\cdots$. However, \eqref{vcountf} says $\#\ccc=f=6+ v_4+2v_5+3v_6+\cdots$, a contradiction.  
\end{proof}

\begin{proposition}\label{111}
	In an $a^3b$-tiling  with some irrational angle, if $\aaa\bbb\ddd$ is a vertex, then all tilings are $2$-layer earth map tilings in Table \ref{tab-2},  and their flip modifications in Table \ref{tab-3} when $\bbb$ is an integer multiple of $\ccc$. 
\end{proposition}

\begin{proof}
	By $\aaa\bbb\ddd$ and Lemma \ref{anglesum}, $\ccc=\tfrac{4}{f}$ and $\ccc^{\frac{f}{2}}$ could be a vertex.  We have $\aaa\neq\ddd$ and $\bbb\neq\ccc$ by Proposition $9$ in \cite{lw}. We will discuss the following cases according to Lemma \ref{3-abd}.
	
\subsubsection*{Case $\aaa^2\ccc^n$ or $\ccc^n\ddd^2$ appears}

If $\aaa^2\ccc^n$ appears, then $\aaa<1$ and $\bbb+\ddd>1$. By Irrational Angle Lemma, any vertex $\bn=\vn$ must satisfy $n_2=n_4$.  By Parity Lemma, $\bbb\cdots=\bbb\ddd\cdots=\aaa\bbb\ddd$. By Balance Lemma, all vertices are $\{f\aaa\bbb\ddd,2\ccc^{\frac{f}{2}}\}$, a contradiction. If $\ccc^n\ddd^2$ appears, we get a similar contradiction.	
	
\subsubsection*{Case $\ccc^{\frac{f}{2}}$ appears } 
If  $\ccc^{\frac{f}{2}}=\thin^{\ddd}\ccc^{\bbb}\thin^{\bbb}\ccc^{\ddd}\thin\cdots$, then this gives a vertex $\bbb^2\cdots$. By Lemma \ref{aadlemma}, $\ddd^2\cdots$ is a vertex. By Balance Lemma,  $\aaa^2\cdots$ is a vertex. 
Therefore, $\text{all angles}<1$. By $\aaa\bbb\ddd$, we have $\aaa+\ddd>1$.

If $\bbb<\ccc$, then $\aaa>\ddd$ by Lemma \ref{geometry1}. By $R(\aaa^2\cdots)<\bbb,\ccc,2\ddd$ and Parity Lemma, there is no $\aaa^2\cdots$, a contradiction.

If $\bbb>\ccc$, then $\aaa<\ddd$ by Lemma $\ref{geometry1}'$. By $R(\ddd^2\cdots)<2\aaa,\bbb$ and Parity Lemma, we get $\ddd^2\cdots=\ccc^n\ddd^2$,  a contradiction by the previous case.

\vspace{9pt}
\noindent Therefore,  $\ccc^{\frac{f}{2}}=\thin^{\bbb}\ccc^{\ddd}\thin^{\bbb}\ccc^{\ddd}\thin^{\bbb}\ccc^{\ddd}\thin\cdots$ by Lemma \ref{aadlemma}. 
If $\aaa<\ddd$, then $\bbb\ddd\cdots=\aaa\bbb\ddd$ by Parity Lemma. 
If $\aaa>\ddd$, then $\bbb<\ccc=\frac4f\le\frac23$ by Lemma \ref{geometry1}. Then  $\aaa+\ddd>1$ by $\aaa\bbb\ddd$. By $R(\aaa^2\cdots)<\bbb,\ccc,2\ddd$ and Parity Lemma, there is no $\aaa^2\cdots$. By Balance Lemma, $\bbb\ddd\cdots=\aaa\bbb\ddd$.

In Fig.\,\ref{abd-1}, $\ccc^{\frac{f}{2}}=\thin^{\bbb}\ccc_1^{\ddd}\thin^{\bbb}\ccc_2^{\ddd}\thin^{\bbb}\ccc_3^{\ddd}\thin\cdots$ determines $T_1,T_2,T_3$. Then $\bbb_2\ddd_1\cdots=\aaa_4\bbb_2\ddd_1$ determines $T_4$; $\bbb_3\ddd_2\cdots=\aaa_5\bbb_3\ddd_2$ and $\aaa_2\bbb_4\ddd_5\cdots=\aaa_2\bbb_4\ddd_5$ determine $T_5$. The argument started at $\aaa_4\bbb_2\ddd_1$ can be repeated at $\aaa_5\bbb_3\ddd_2$. More repetitions give the
unique tiling of $f$ tiles with $2\ccc^{\frac{f}{2}}$ and $f\aaa\bbb\ddd$.
\begin{figure}[htp]
	\centering
	\begin{tikzpicture}[>=latex,scale=0.6] 
	\foreach \a in {0,1,2}
	{
		\begin{scope}[xshift=2*\a cm] 
		\draw (0,0)--(0,-2)
		(2,0)--(2,-2)--(3.5,-3)--(3.5,-5)
		(1.5,-3)--(1.5,-5)
		(0,-2)--(1.5,-3);
		\draw[line width=1.5] (1.5,-3)--(2,-2);
		\node at (1.4,-2.55){\small $\aaa$};
		\node at (1.8,-3.2){\small $\ddd$};
		\node at (3.1,-3.2){\small $\bbb$};
		\node at (1,0){\small $\ccc$};
		\node at (2.5,-5){\small $\ccc$};
		
		\node at (2.1,-2.45){\small $\aaa$};
		\node at (1.7,-1.8){\small $\ddd$};
		\node at (0.3,-1.7){\small $\bbb$};
		\end{scope}
	}
	
	\fill (8,-2) circle (0.05); \fill (7.7,-2) circle (0.05);
	\fill (8.3,-2) circle (0.05);
	
	\node[draw,shape=circle, inner sep=0.5] at (1,-1) {\small $1$};
	\node[draw,shape=circle, inner sep=0.5] at (3,-1) {\small $2$};
	\node[draw,shape=circle, inner sep=0.5] at (5,-1) {\small $3$};
	\node[draw,shape=circle, inner sep=0.5] at (2.5,-4) {\small $4$};
	\node[draw,shape=circle, inner sep=0.5] at (4.5,-4) {\small $5$};
	\node[draw,shape=circle, inner sep=0.5] at (6.5,-4) {\small $6$};
	
	\end{tikzpicture}
	\caption{ $2$-layer earth map $a^3b$-tiling $T(f\aaa\bbb\ddd,2\ccc^{\frac{f}{2}})$.} \label{abd-1}
\end{figure}


\subsubsection*{Case $\bbb\ccc^m$ appears ($\ccc^{\frac{f}{2}}$ never appears)}	
By $\aaa\bbb\ddd$ and $\bbb\ccc^m$ ($m\ge2$), we have $\aaa+\ddd=m\ccc$ and $\bbb=(\frac{f}{2}-m)\ccc$.
By  Irrational Angle Lemma, any vertex $\bn$ must satisfy $n_1=n_4$. We divide the following discussion into two parts. 

$\boldsymbol{1. \,\aaa+\ddd\le1,\bbb\ge1}$. There is no $\bbb^2\cdots$. By $n_1=n_4$ and Lemma \ref{anglesum}, $\aaa\cdots=\aaa\ddd\cdots=\aaa\bbb\ddd$ or $\aaa^n\ccc^{\frac{f}{2}-mn}\ddd^n\,(n\in[1, \frac{f}{2m}]\cap\mathbb{N})$ and $\bbb\cdots=\aaa\bbb\ddd$ or $\bbb\ccc^m$. 
Therefore, 
\[\text{AVC}\sub\{\aaa\bbb\ddd,\aaa^n\ccc^{\frac{f}{2}-mn}\ddd^n,\bbb\ccc^m\,(n\in[1, \tfrac{f}{2m}]\cap\mathbb{N})\}.\]  
By Lemma $\ref{geometry4}'$, $\ccc+2\ddd>1$. By $\aaa\bbb\ddd$ and $\bbb\ccc^m$, we have $\aaa+\ddd=m\ccc\ge2\ccc$. If $n\ge5$, then $5\aaa+5\ddd>\aaa+\ddd+4\ddd\ge2\ccc+4\ddd>2$, a contradiction. Therefore, $n\le4$ (we will show $n\le3$ in the end of this section using the moduli space). 

For $n=1$, we have the unique AAD $\aaa\ccc^{\frac{f}{2}-m}\ddd=\thin^\bbb\aaa^\ddd\thick^\aaa\ddd^\ccc\thin^\bbb\ccc^\ddd\thin\cdots$.
For $n=2,3,4$, if $\aaa^n\ccc^{\frac{f}{2}-mn}\ddd^n=\thin^{\bbb}\aaa_1^{\ddd}\thick^{\ddd}\aaa_2^{\bbb}\thin^{\ddd}\ccc_3^{\bbb}\thin\cdots$
which determines $T_1,T_2,T_3$ in the left of Fig.\,\ref{abd-2}.
Then $\bbb_2\ddd_3\cdots=\aaa_4\bbb_2\ddd_3$ determines $T_4$; $\bbb_4\ccc_2\cdots=\bbb_4\ccc_2\ccc_5\cdots$ and $\ddd_1\ddd_2\cdots=\ddd_1\ddd_2\ddd_5\cdots$ determine $T_5$. By the AVC, $\ddd_1\ddd_2\ddd_5\cdots=\aaa^n\ccc^{\frac{f}{2}-mn}\ddd^n$, whose AAD implies that  $\thick\aaa^{\bbb}\thin^{\bbb}\aaa\thick\cdots$, $\thick\aaa^{\bbb}\thin^{\bbb}\ccc\thin\cdots$ or $\thin\ccc^{\bbb}\thin^{\bbb}\ccc\thin\cdots$ must appear. This gives a vertex $\bbb^2\cdots$, a contradiction.
If $\aaa^n\ccc^{\frac{f}{2}-mn}\ddd^n=\thin^{\bbb}\aaa_1^{\ddd}\thick^{\ddd}\aaa_2^{\bbb}\thin^{\ccc}\ddd_3^{\aaa}\thick\cdots$ which determines $T_1,T_2,T_3$ in the right of Fig.\,\ref{abd-2}.
Then $\bbb_2\ccc_3\cdots=\thin^{\ccc}\bbb_2^{\aaa}\thin^{\ddd}\ccc_3^{\bbb}\thin\cdots\thin^{\ddd}\ccc_4^{\bbb}\thin$ determines $T_4$; $\bbb_4\ccc_2\cdots=\bbb_4\ccc_2\ccc_5\cdots$ and $\ddd_1\ddd_2\cdots\\=\ddd_1\ddd_2\ddd_5\cdots$ determine $T_5$. Similarly, we get a contradiction. Then we get its unique AAD (here all $\ccc$'s may be distributed freely between these ``$\aaa\thick\ddd$'') \[\aaa^n\ccc^{\frac{f}{2}-mn}\ddd^n=\thin^\bbb\aaa^\ddd\thick^\aaa\ddd^\ccc\thin\cdots\thin^{\bbb}\ccc^{\ddd}\thin\cdots\thin^\bbb\aaa^\ddd\thick^\aaa\ddd^\ccc\thin\cdots\thin^{\bbb}\ccc^{\ddd}\thin\cdots,\,(n=2,3,4),\]
which determines $T_1,T_2$ in Fig.\,\ref{abd-3}. Then $R(\aaa^n\ddd^n\cdots)=\ccc^{\frac{f}{2}-mn}$ and this $\ccc^{\frac{f}{2}-mn}$ determines $(\frac{f}{2}-mn)$ time zones ($(f-2mn)$ tiles). Then $\thin\bbb\thin\ccc_2\thin\cdots=\thin\bbb\thin^{\ddd}\ccc_2^{\bbb}\thin^{\ddd}\ccc^{\bbb}\thin\cdots=\bbb\ccc^m$. Then $R(\bbb\ccc_2\cdots)=R(\bbb\ccc_4\cdots)=\cdots=\ccc^{m-1}$ and each of $n$ such $\ccc^{m-1}$ determines $(m-1)$ time zones ($(2m-2)$ tiles). This tiling can also be obtained by applying the second flip modification $n$ times in Fig.\,\ref{flip1}.
\begin{figure}[htp]
	\centering
	\begin{tikzpicture}[>=latex,scale=0.43] 
		\draw (0,0)--(0,-2)--(8,-2)--(8,0)
		(6,-2)--(6,-5)--(4,-5)
		(6,-2)--(8,-5)--(10,-2)--(12,-2)--(12,0);
		
		\draw[line width=1.5] (4,0)--(4,-5)
		(8,-2)--(10,-2);
		\node at (2,0){\small $\aaa$};\node at (0.35,-1.6){\small $\bbb$};\node at (2,-1.6){\small $\ccc$};\node at (3.65,-1.6){\small$\ddd$};
		\node at (6,0){\small $\aaa$};\node at (4.35,-1.6){\small $\ddd$};\node at (6,-1.6){\small $\ccc$};\node at (7.65,-1.6){\small $\bbb$};
		\node at (10,0){\small $\ccc$};\node at (8.35,-1.6){\small $\ddd$};\node at (10,-1.6){\small $\aaa$};\node at (11.65,-1.6){\small $\bbb$};
		\node at (4.35,-2.4){\small $\ddd$};\node at (5.65,-2.4){\small $\ccc$};\node at (4.35,-4.6){\small $\aaa$};\node at (5.65,-4.6){\small $\bbb$};
		\node at (6.8,-2.45){\small $\bbb$};\node at (8,-2.4){\small $\aaa$};\node at (9.2,-2.45){\small $\ddd$};\node at (8,-4.4){\small $\ccc$};
		\node[draw,shape=circle, inner sep=0.5] at (2,-0.8) {\small $1$};
		\node[draw,shape=circle, inner sep=0.5] at (6,-0.8) {\small $2$};
		\node[draw,shape=circle, inner sep=0.5] at (10,-0.8) {\small $3$};
		\node[draw,shape=circle, inner sep=0.5] at (8,-3.3) {\small $4$};
		\node[draw,shape=circle, inner sep=0.5] at (5,-3.5) {\small $5$};
		
		\begin{scope}[xshift=2 cm]
			\draw 	
			(14,0)--(14,-2)--(26,-2)
			(17,-5)--(19,-5)--(20,-2)--(21,-5)
			(22,0)--(22,-2)--(23,-5);
			
			\draw[line width=1.5] 
			(18,0)--(18,-2)--(17,-5)
			(21,-5)--(23,-5)
			(26,0)--(26,-2);
			
			\node at (2+14,0){\small $\aaa$};\node at (0.35+14,-1.6){\small $\bbb$};\node at (2+14,-1.6){\small $\ccc$};\node at (3.65+14,-1.6){\small $\ddd$};
			\node at (6+14,0){\small $\aaa$};\node at (4.35+14,-1.6){\small $\ddd$};\node at (6+14,-1.6){\small $\ccc$};\node at (7.65+14,-1.6){\small $\bbb$};
			\node at (10+14,0){\small $\ddd$};\node at (8.35+14,-1.6){\small $\ccc$};\node at (10+14,-1.6){\small $\bbb$};\node at (11.65+14,-1.6){\small $\aaa$};
			\node at (18.2,-2.4){\small $\ddd$};\node at (19.5,-2.4){\small $\ccc$};\node at (17.5,-4.6){\small $\aaa$};\node at (18.8,-4.6){\small $\bbb$};
			\node at (20.6,-2.5){\small $\bbb$};\node at (21.7,-2.45){\small $\ccc$};\node at (21.2,-4.6){\small $\aaa$};\node at (22.5,-4.6){\small $\ddd$};
			
			\fill (0,-2) circle (0.05); \fill (16,-2) circle (0.05); \fill (24,-2) circle (0.05);

			\node[draw,shape=circle, inner sep=0.5] at (2+14,-0.8) {\small $1$};
			\node[draw,shape=circle, inner sep=0.5] at (6+14,-0.8) {\small $2$};
			\node[draw,shape=circle, inner sep=0.5] at (10+14,-0.8) {\small $3$};
			\node[draw,shape=circle, inner sep=0.5] at (21.5,-3.5) {\small $4$};
			\node[draw,shape=circle, inner sep=0.5] at (18.5,-3.5) {\small $5$};
		\end{scope}
	\end{tikzpicture}
	\caption{Two possibilities for $\protect\thin\aaa\protect\thick\aaa\protect\thin\cdots$.} \label{abd-2}
\end{figure}
\begin{figure}[htp]
	\centering	    	
	\begin{tikzpicture}[>=latex,scale=0.5] 
	\foreach \a in {0,1}
	{
		\begin{scope}[xshift=12*\a cm] 
		\draw (0,0)--(0,-8)
		(6,0)--(6,-8)
		(0,-6)--(2,-2)--(6,-2)--(4,-6)--(0,-6);
		\draw[line width=1.5] (2,0)--(2,-2)
		(4,-6)--(4,-8);
		\node at (1,0){\small $\aaa$};
		\node at (4,0){\small $\ddd$};
		\node at (0.3,-2){\small $\bbb$};
		\node at (1.7,-2){\small $\ddd$};
		\node at (2.3,-1.7){\small $\aaa$};
		\node at (4,-1.6){\small $\bbb$};
		\node at (5.7,-1.7){\small $\ccc$};
		\node at (0.25,-5.1){\small $\ccc$};
		\node at (5.75,-3.1){\small $\ccc$};
		\node at (0.3,-6.3){\small $\ccc$};
		\node at (2,-6.4){\small $\bbb$};
		\node at (3.7,-6.3){\small $\aaa$};
		\node at (4.3,-6){\small $\ddd$};
		\node at (5.7,-6){\small $\bbb$};
		\node at (2,-8){\small $\ddd$};
		\node at (5,-8){\small $\aaa$};
		
		\fill (4,-2) circle (0.04);
		\fill (2,-6) circle (0.04);

		\fill (1,-5.8) circle (0.05);
		\fill (0.8,-5.6) circle (0.05);
		\fill (0.6,-5.4) circle (0.05);
		
		\fill (1+4.4,-5.8+3.1) circle (0.05);
		\fill (0.8+4.4,-5.6+3.1) circle (0.05);
		\fill (0.6+4.4,-5.4+3.1) circle (0.05);
		
		\fill (6,-6) circle (0.04);
		\fill (0,-2) circle (0.04);
		\fill (8,-2) circle (0.04);
		\fill (10,-6) circle (0.04);
		
		\end{scope}
	}
	
	\foreach \b in {1,7}
	{
		\begin{scope}[xshift=2*\b cm] 
		\draw (6,0)--(6,-8)
		(8,0)--(8,-8);
		\draw[line width=1.5] (6,-6)--(8,-2);
		\node at (7,0){\small $\ccc$};
		\node at (6.3,-2){\small $\bbb$};
		\node at (7.7,-2){\small $\ddd$};
		\node at (6.3,-4.8){\small $\aaa$};
		\node at (7.7,-3.2){\small $\aaa$};
		\node at (6.3,-6){\small $\ddd$};
		\node at (7.7,-6){\small $\bbb$};
		\node at (7,-8){\small $\ccc$};
		
		\fill (8.6,-4) circle (0.05);
		\fill (8.9,-4) circle (0.05);
		\fill (9.2,-4) circle (0.05);
		
		\fill (4.6,-4) circle (0.05);
		\fill (4.9,-4) circle (0.05);
		\fill (5.2,-4) circle (0.05);
		\end{scope}
	}
	
	\fill (0,-2) circle (0.04);
	\fill (12,-6) circle (0.04);

	\node[draw,shape=circle, inner sep=0.5] at (1,-0.8) {\small $1$};
	\node[draw,shape=circle, inner sep=0.5] at (4,-0.8) {\small $2$};
	\node[draw,shape=circle, inner sep=0.5] at (13,-0.8) {\small $3$};
	\node[draw,shape=circle, inner sep=0.5] at (16,-0.8) {\small $4$};

	\end{tikzpicture}
	\caption{Many different tilings for  $\{(f-2n)\aaa\bbb\ddd,2\aaa^n\ccc^{\frac{f}{2}-mn}\ddd^n,2n\bbb\ccc^m\}$.} \label{abd-3}
\end{figure}

$\boldsymbol{2.\, \aaa+\ddd>1,\bbb<1}$. By $n_1=n_4$ and Lemma \ref{anglesum}, $\aaa\cdots=\aaa\ddd\cdots=\aaa\bbb\ddd$ or $\aaa\ccc^{\frac{f}{2}-m}\ddd$.  Therefore, \[\text{AVC}\sub\{\aaa\bbb\ddd,\aaa\ccc^{\frac{f}{2}-m}\ddd,\bbb^n\ccc^{\frac{f}{2}-\frac {nf}{2}+mn}\,(1\le n\le\tfrac{f}{f-2m})\}.\] 
By the AVC, we know that $\aaa\thin\ddd\cdots$ and $\ddd^2\cdots$ are not vertices. In Fig.\,\ref{abd-4}, we have the unique AAD $\bbb\ccc^m=\thin^{\ccc}\bbb_1^{\aaa}\thin^{\bbb}\ccc_2^{\ddd}\thin^{\bbb}\ccc_3^{\ddd}\thin\cdots$. Then $R(\bbb_1\cdots)=\ccc^m$ and this $\ccc^m$ determines $m$ time zones ($2m$ tiles). Then $R(\aaa_2\ddd_5\cdots)=\ccc^{\frac{f}{2}-m}$ and this $\ccc^{\frac{f}{2}-m}$ determines $(\frac{f}{2}-m)$ time zones (($f-2m$) tiles). This tiling can also be obtained by applying the first flip modification in Fig.\,\ref{flip1}.
\begin{figure}[htp]
	\centering
	\begin{tikzpicture}[>=latex,scale=0.5] 
		\foreach \a in {0,1}
		{
			\begin{scope}[xshift=3*\a cm] 
				\draw (4,0)--(4,-9)
				(7,0)--(7,-9);
				\draw[line width=1.5] (7,-2)--(4,-7);
				\node at (5.5,0){\small $\ccc$};
				\node at (4.3,-2){\small $\bbb$};
				\node at (6.7,-1.9){\small $\ddd$};
				\node at (6.7,-3.3){\small $\aaa$};
				\node at (4.3,-5.7){\small $\aaa$};
				\node at (6.7,-7){\small $\bbb$};
				\node at (4.3,-7.1){\small $\ddd$};
				\node at (5.5,-9){\small $\ccc$};
			\end{scope}
		}
		
		\draw (0,0)--(0,-9)
		(4,0)--(4,-9)
		(0,-2)--(3,-2)--(3,-4)--(4,-7)--(2,-5)
		(0,-2)--(2,-4)
		(0,-2)--(1,-5)--(1,-7)--(4,-7);
		\draw[line width=1.5] (3,-2)--(4,-2)
		(2,-4)--(3,-4)
		(1,-5)--(2,-5)
		(0,-7)--(1,-7);
		
		\node at (2,0){\small $\bbb$}; 
		\node at (0.3,-1.6){\small $\ccc$};
		\node at (3.7,-1.6){\small $\aaa$};
		\node at (3,-1.6){\small $\ddd$};
		\node at (1,-2.4){\small $\ccc$};
		\node at (2.7,-2.4){\small $\bbb$};
		\node at (3.25,-2.4){\small $\aaa$};
		\node at (3.75,-2.4){\small $\ddd$};
		\node at (2.1,-3.6){\small $\ddd$};
		\node at (2.75,-3.6){\small $\aaa$};
		\node at (3.25,-3.75){\small $\bbb$};
		\node at (0.25,-3.8){\small $\ccc$};
		\node at (3.75,-5.2){\small $\ccc$};
		\node at (0.7,-5){\small $\bbb$};
		\node at (0.25,-6.6){\small $\ddd$};
		\node at (0.75,-6.6){\small $\aaa$};
		\node at (1.3,-5.4){\small $\aaa$};
		\node at (2,-5.4){\small $\ddd$};
		\node at (1.3,-6.6){\small $\bbb$};
		\node at (3,-6.6){\small $\ccc$};
		\node at (0.25,-7.4){\small $\aaa$};
		\node at (1,-7.4){\small $\ddd$};
		\node at (3.7,-7.4){\small $\ccc$};
		\node at (2,-9){\small $\bbb$};

		\fill (0.56,-3.21) circle (0.05);
		\fill (0.68,-3.13) circle (0.05);
		\fill (0.80,-3.05) circle (0.05);
		
		\fill (0.56+2.55,-3.21-2.6) circle (0.05);
		\fill (0.68+2.55,-3.13-2.6) circle (0.05);
		\fill (0.80+2.55,-3.05-2.6) circle (0.05);
		
		\fill (11,-5) circle (0.05);
		\fill (11.3,-5) circle (0.05);
		\fill (11.6,-5) circle (0.05);
		
		\node[draw,shape=circle, inner sep=0.5] at (2,-1) {\small $1$};
		\node[draw,shape=circle, inner sep=0.5] at (5.5,-1) {\small $2$};
		\node[draw,shape=circle, inner sep=0.5] at (8.5,-1) {\small $3$};
		\node[draw,shape=circle, inner sep=0.5] at (8.5,-8) {\small $4$};
		\node[draw,shape=circle, inner sep=0.5] at (5.5,-8) {\small $5$};
		\node[draw,shape=circle, inner sep=0.5] at (2,-8) {\small $6$};
		
		\fill (0,-2) circle (0.15); \fill (4,-7) circle (0.15);
		
		\fill (10,-7) circle (0.04);
		
	\end{tikzpicture}
	\caption{$T((f-2)\,\aaa\bbb\ddd,2\aaa\ccc^{\frac{f}{2}-m}\ddd,2\bbb\ccc^m)$.} \label{abd-4}
\end{figure}
\subsubsection*{Case $\bbb^2\ccc^k$ appears ($\bbb\ccc^m,\ccc^{\frac{f}{2}}$ never appear)}	
Similarly, we deduce that $\text{AVC}\sub\{\aaa\bbb\ddd,\aaa\ccc^{\frac{f-2k}{4}}\ddd,\bbb^{n}\ccc^{\frac{(2-n)f+2kn}{4}}\,(n\in[1, \frac{2f}{f-2k}]\cap\mathbb{N})\}$. The tiling is shown in Fig.\,\ref{abd-5-}. Depending on the space between two flips, there are $\lfloor \frac{k+2}{2}\rfloor$ different tilings with the same set of vertices.

\begin{figure}[htp]
	\centering
	\begin{tikzpicture}[>=latex,scale=0.5] 
		\foreach \a in {0}
		{
			\begin{scope}[xshift=2 cm+11*\a cm] 
				\draw (4,0)--(4,-9)
				(7,0)--(7,-9);
				\draw[line width=1.5] (7,-2)--(4,-7);
				\node at (5.5,0){\small $\ccc$};
				\node at (4.3,-2){\small $\bbb$};
				\node at (6.7,-1.9){\small $\ddd$};
				\node at (6.7,-3.3){\small $\aaa$};
				\node at (4.3,-5.7){\small $\aaa$};
				\node at (6.7,-7){\small $\bbb$};
				\node at (4.3,-7.1){\small $\ddd$};
				\node at (5.5,-9){\small $\ccc$};
			\end{scope}
		}
		\foreach \b in {0,1}
		{
			\begin{scope}[xshift=11*\b cm] 
				\draw (0,0)--(0,-9)
				(4,0)--(4,-9)
				(0,-2)--(3,-2)--(3,-4)--(4,-7)--(2,-5)
				(0,-2)--(2,-4)
				(0,-2)--(1,-5)--(1,-7)--(4,-7);
				\draw[line width=1.5] (3,-2)--(4,-2)
				(2,-4)--(3,-4)
				(1,-5)--(2,-5)
				(0,-7)--(1,-7);
				
				\node at (2,0){\small $\bbb$}; 
				\node at (0.3,-1.6){\small $\ccc$};
				\node at (3.7,-1.6){\small $\aaa$};
				\node at (3,-1.6){\small $\ddd$};
				\node at (1,-2.4){\small $\ccc$};
				\node at (2.7,-2.4){\small $\bbb$};
				\node at (3.3,-2.4){\small $\aaa$};
				\node at (3.75,-2.4){\small $\ddd$};
				\node at (2.1,-3.6){\small $\ddd$};
				\node at (2.75,-3.6){\small $\aaa$};
				\node at (3.25,-3.7){\small $\bbb$};
				\node at (0.25,-3.8){\small $\ccc$};
				\node at (3.75,-5.2){\small $\ccc$};
				\node at (0.7,-5){\small $\bbb$};
				\node at (0.25,-6.6){\small $\ddd$};
				\node at (0.75,-6.6){\small $\aaa$};
				\node at (1.3,-5.4){\small $\aaa$};
				\node at (2,-5.5){\small $\ddd$};
				\node at (1.3,-6.6){\small $\bbb$};
				\node at (3,-6.6){\small $\ccc$};
				\node at (0.25,-7.4){\small $\aaa$};
				\node at (1,-7.4){\small $\ddd$};
				\node at (3.7,-7.4){\small $\ccc$};
				\node at (2,-9){\small $\bbb$};

			\end{scope}
		}
		\node[draw,shape=circle, inner sep=0.5] at (2,-1) {\small $1$};
		\node[draw,shape=circle, inner sep=0.5] at (5.5+2,-1) {\small $2$};
		\node[draw,shape=circle, inner sep=0.5] at (2,-8) {\small $3$};
		\node[draw,shape=circle, inner sep=0.5] at (5.5+2,-8) {\small $4$};
		
		\fill (0.56,-3.21) circle (0.05);
		\fill (0.68,-3.13) circle (0.05);
		\fill (0.80,-3.05) circle (0.05);
		
		\fill (0.56+2.55,-3.21-2.6) circle (0.05);
		\fill (0.68+2.55,-3.13-2.6) circle (0.05);
		\fill (0.80+2.55,-3.05-2.6) circle (0.05);
		
		\fill (8+1.6,-5) circle (0.05);
		\fill (8.3+1.6,-5) circle (0.05);
		\fill (8.6+1.6,-5) circle (0.05);
		
		\fill (8-3.2,-5) circle (0.05);
		\fill (8.3-3.2,-5) circle (0.05);
		\fill (8.6-3.2,-5) circle (0.05);

		\fill (9,-7) circle (0.04); \fill (6,-2) circle (0.04);

	\end{tikzpicture}
	\caption{Many different tilings for  $\{(f-4)\aaa\bbb\ddd,4\aaa\ccc^{\frac{f-2k}{4}}\ddd,2\bbb^2\ccc^k\}$.} \label{abd-5-}
\end{figure}

\subsubsection*{Case $\aaa\ccc^m\ddd$ appears ($\bbb^2\ccc^k,\bbb\ccc^m,\ccc^{\frac{f}{2}}$ never appear)}
By $\aaa\bbb\ddd$ and $\aaa\ccc^m\ddd$ ($m\ge2$ by Proposition $9$ in \cite{lw}), we have $\bbb=m\ccc$ and $\aaa+\ddd=(\frac{f}{2}-m)\ccc$.
By Irrational Angle Lemma, any vertex $\bn$ must satisfy $n_1=n_4$. If $\bbb\ge1$, then $\aaa+\ddd\le1$ by $\aaa\bbb\ddd$. By $n_1=n_4$, we deduce that $\bbb\cdots=\aaa\bbb\ddd$. This implies $\ccc\cdots=\ccc^{\frac{f}{2}}$, a contradiction. Therefore, $\bbb<1$.
By $\aaa\bbb\ddd$, $\aaa+\ddd>1$. By $n_1=n_4$, we deduce that $\aaa\cdots=\aaa\ddd\cdots=\aaa\bbb\ddd$ or $\aaa\ccc^m\ddd$ and $\bbb\cdots=\aaa\bbb\ddd$ or $\bbb^n\ccc^{\frac{f}{2}-mn}\,(n\in[1, \frac{f}{2m}]\cap\mathbb{N})$. Therefore, \[\text{AVC}\sub\{\aaa\bbb\ddd,\aaa\ccc^m\ddd,\bbb^n\ccc^{\frac{f}{2}-mn}\,(n\in[1, \tfrac{f}{2m}]\cap\mathbb{N})\}.\] 
By Lemma \ref{proposition-1}, we have $\aaa<1$ or $\ddd<1$. Then $2\bbb+\ccc>1$ or $\bbb+2\ccc>1$ by Lemma $\ref{geometry4}'$ and \ref{geometry4}. By $\bbb=m\ccc$ and $\bbb^n\ccc^{\frac{f}{2}-mn}$, we get $n\le4$ (we will show $n\le3$ in the end of this section using the moduli space).

By the AVC, we know that $\aaa\thin\ddd\cdots,\aaa^2\cdots$ and $\ddd^2\cdots$ are not vertices. In Fig.\,\ref{abd-5}, we have the unique AAD $\bbb^n\ccc^{\frac{f}{2}-mn}=\thin^{\ccc}\bbb_1^{\aaa}\thin\cdots\thin^{\bbb}\ccc_2^{\ddd}\thin\cdots\thin^{\ccc}\bbb^{\aaa}\thin\cdots\thin^{\bbb}\ccc^{\ddd}\thin\cdots$. Then $R(\bbb^n\cdots)=\ccc^{\frac{f}{2}-mn}$ and this $\ccc^{\frac{f}{2}-mn}$ determines $(\frac{f}{2}-mn)$ time zones ($(f-2mn)$ tiles). We deduce that $\ccc_1\cdots=\thin^{\bbb}\ccc_1^{\ddd}\thin^{\bbb}\ccc^{\ddd}\thin\cdots=\aaa\ccc^m\ddd$. Then $R(\aaa\ddd\cdots)=\ccc^m$ and each of $n$ such $\ccc^m$ determines $m$ time zones ($2m$ tiles). This tiling can also be obtained by applying the first flip modification $n$ times in Fig.\,\ref{flip1}.

\begin{figure}[htp]
	\centering
	\begin{tikzpicture}[>=latex,scale=0.5] 
	\foreach \a in {0,1}
	{
		\begin{scope}[xshift=2 cm+11*\a cm] 
		\draw (4,0)--(4,-9)
		(7,0)--(7,-9);
		\draw[line width=1.5] (7,-2)--(4,-7);
		\node at (5.5,0){\small $\ccc$};
		\node at (4.3,-2){\small $\bbb$};
		\node at (6.7,-1.9){\small $\ddd$};
		\node at (6.7,-3.3){\small $\aaa$};
		\node at (4.3,-5.7){\small $\aaa$};
		\node at (6.7,-7){\small $\bbb$};
		\node at (4.3,-7.1){\small $\ddd$};
		\node at (5.5,-9){\small $\ccc$};
		\end{scope}
	}
	\foreach \b in {0,1}
	{
		\begin{scope}[xshift=11*\b cm] 
		\draw (0,0)--(0,-9)
		(4,0)--(4,-9)
		(0,-2)--(3,-2)--(3,-4)--(4,-7)--(2,-5)
		(0,-2)--(2,-4)
		(0,-2)--(1,-5)--(1,-7)--(4,-7);
		\draw[line width=1.5] (3,-2)--(4,-2)
		(2,-4)--(3,-4)
		(1,-5)--(2,-5)
		(0,-7)--(1,-7);
		
		\node at (2,0){\small $\bbb$}; 
		\node at (0.3,-1.6){\small $\ccc$};
		\node at (3.7,-1.6){\small $\aaa$};
		\node at (3,-1.6){\small $\ddd$};
		\node at (1,-2.4){\small $\ccc$};
		\node at (2.7,-2.4){\small $\bbb$};
		\node at (3.3,-2.4){\small $\aaa$};
		\node at (3.75,-2.4){\small $\ddd$};
		\node at (2.1,-3.6){\small $\ddd$};
		\node at (2.75,-3.6){\small $\aaa$};
		\node at (3.25,-3.7){\small $\bbb$};
		\node at (0.25,-3.8){\small $\ccc$};
		\node at (3.75,-5.2){\small $\ccc$};
		\node at (0.7,-5){\small $\bbb$};
		\node at (0.25,-6.6){\small $\ddd$};
		\node at (0.75,-6.6){\small $\aaa$};
		\node at (1.3,-5.4){\small $\aaa$};
		\node at (2,-5.5){\small $\ddd$};
		\node at (1.3,-6.6){\small $\bbb$};
		\node at (3,-6.6){\small $\ccc$};
		\node at (0.25,-7.4){\small $\aaa$};
		\node at (1,-7.4){\small $\ddd$};
		\node at (3.7,-7.4){\small $\ccc$};
		\node at (2,-9){\small $\bbb$};

		\fill (0.56,-3.21) circle (0.05);
		\fill (0.68,-3.13) circle (0.05);
		\fill (0.80,-3.05) circle (0.05);
		
		\fill (0.56+2.55,-3.21-2.6) circle (0.05);
		\fill (0.68+2.55,-3.13-2.6) circle (0.05);
		\fill (0.80+2.55,-3.05-2.6) circle (0.05);
		
		\fill (8+1.6,-5) circle (0.05);
		\fill (8.3+1.6,-5) circle (0.05);
		\fill (8.6+1.6,-5) circle (0.05);
		
		\fill (8-3.2,-5) circle (0.05);
		\fill (8.3-3.2,-5) circle (0.05);
		\fill (8.6-3.2,-5) circle (0.05);

		\fill (9,-7) circle (0.04); \fill (6,-2) circle (0.04);
		\end{scope}
	}
	\node[draw,shape=circle, inner sep=0.5] at (2,-1) {\small $1$};
	\node[draw,shape=circle, inner sep=0.5] at (5.5+2,-1) {\small $2$};
	\node[draw,shape=circle, inner sep=0.5] at (2,-8) {\small $3$};
	\node[draw,shape=circle, inner sep=0.5] at (5.5+2,-8) {\small $4$};

	\end{tikzpicture}
	\caption{Many different tilings for  $\{(f-2n)\aaa\bbb\ddd,2n\,\aaa\ccc^m\ddd,2\bbb^n\ccc^{\frac{f}{2}-mn}\}$.} \label{abd-5}
\end{figure}
\end{proof}

\begin{remark}
	All quadrilaterals in the families of $2$-layer earth map $a^3b$-tilings have some irrational angles except the $12$ sporadic and $3$ infinite sequences of cases listed in \cite{lw}. In particular when   $(\aaa,\bbb,\ccc,\ddd)=(\frac{2}{f},\frac{4f-4}{3f},\frac{4}{f},\frac{2f-2}{3f})$ and $f=6k+4$ $(k\ge1)$, there are three extra tilings with all vertices being 
	$\{(f-6)\aaa\bbb\ddd,2\aaa\ddd^3,2\aaa^2\bbb\ccc^{\frac{f-4}{6}},4\bbb\ccc^{\frac{f+2}{6}}\}$.
\end{remark}

\subsection*{Moduli of $2$-layer earth map $a^3b$-tilings}
\begin{figure}[htp]
	\centering
	\begin{overpic}[scale=0.18]{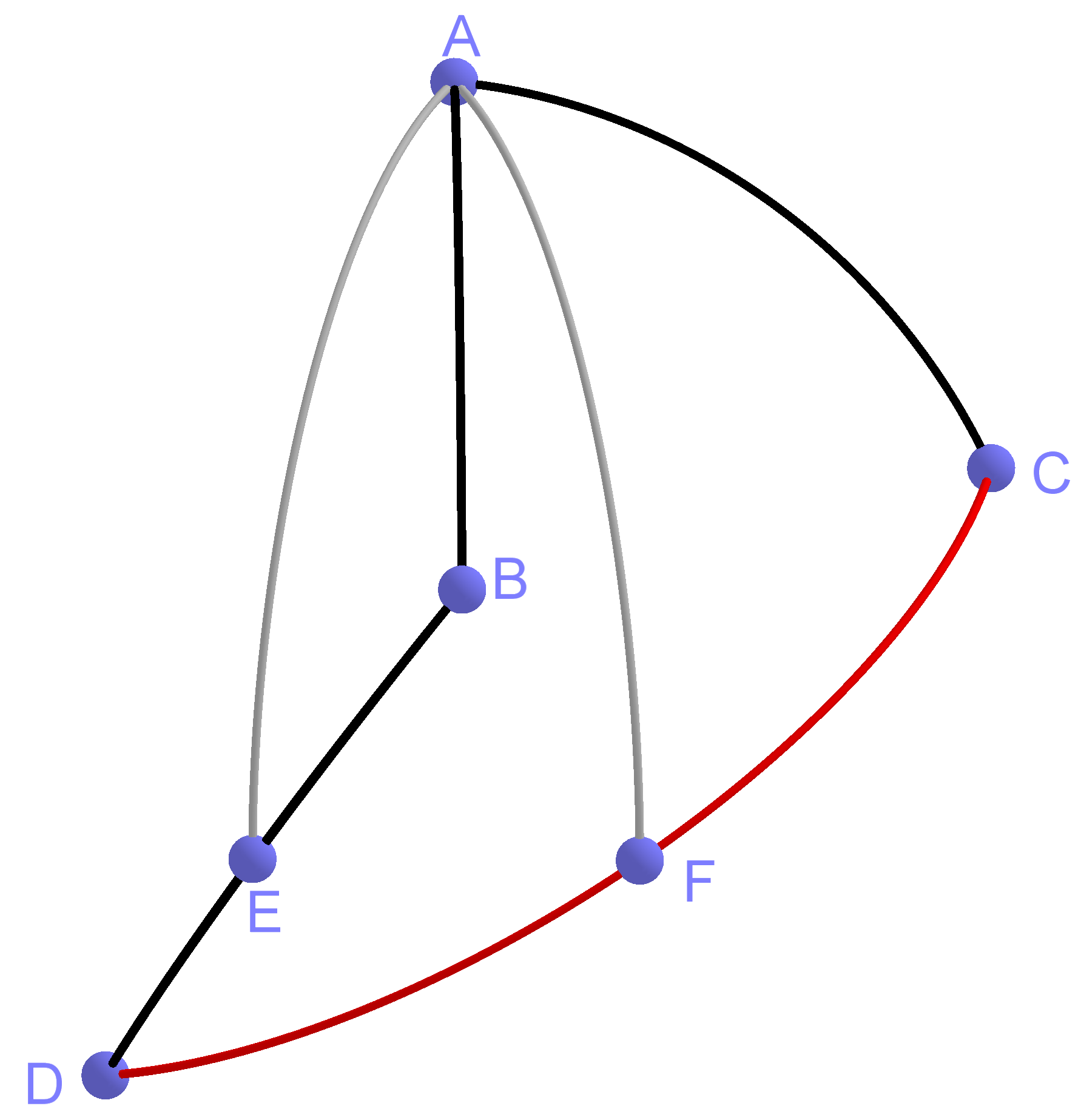}
		\put(24,92){\tiny $A$}
		\put(25,47){\tiny $B$}
		\put(75,58){\tiny $C$}
		\put(13,15){\tiny $E$}
		\put(49,15){\tiny $F$}
		\put(-8,0){\tiny $D$}
	\end{overpic}   \hspace{20pt}
    \begin{overpic}[scale=0.225]{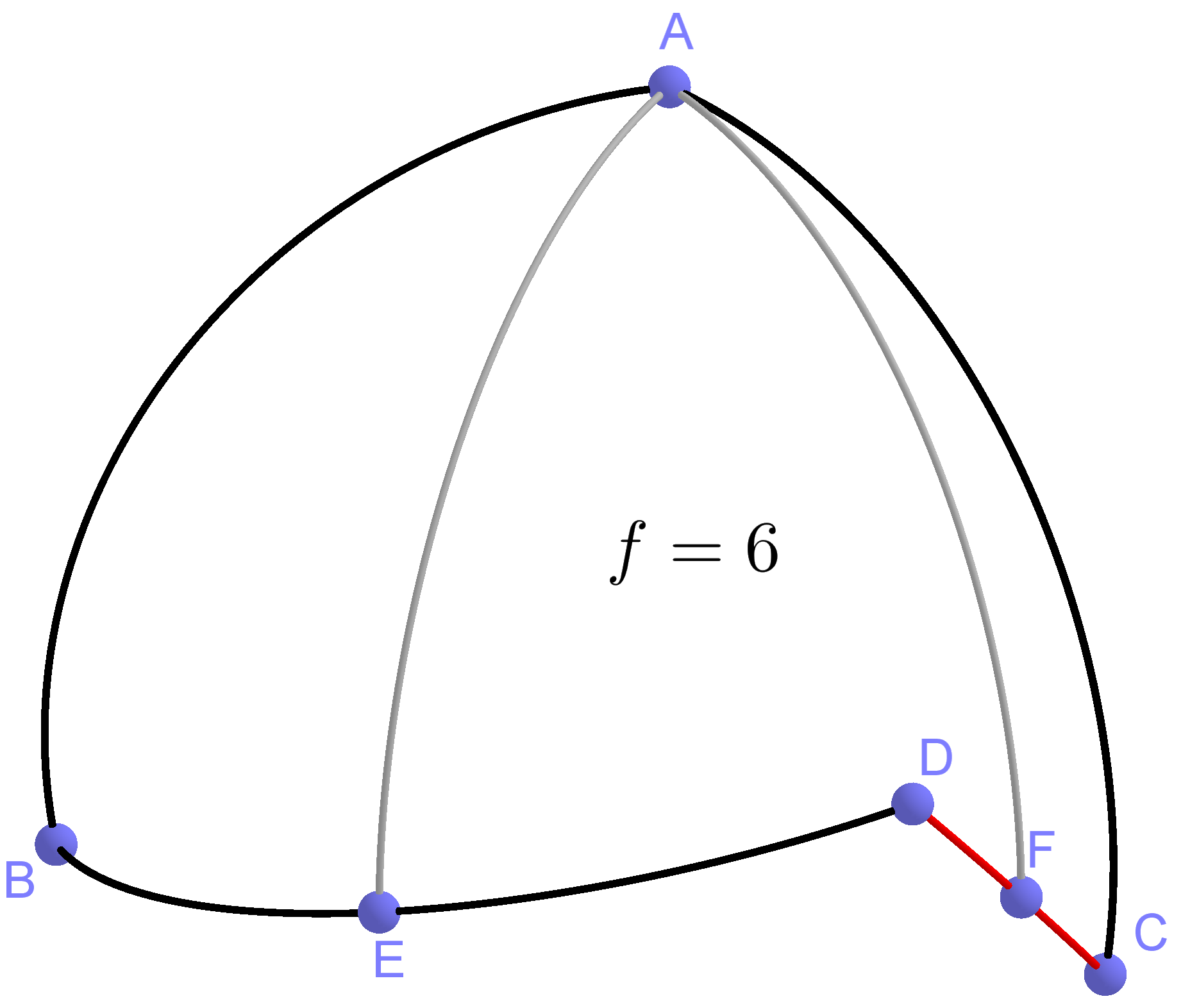}
    	\put(55,75){\tiny $A$}
    	\put(-4,8){\tiny $B$}
    	\put(100,0){\tiny $C$}
    	\put(28,1){\tiny $E$}
    	\put(86,2){\tiny $F$}
    	\put(76,11){\tiny $D$}
    \end{overpic}   \hspace{20pt}
    \begin{overpic}[scale=0.19]{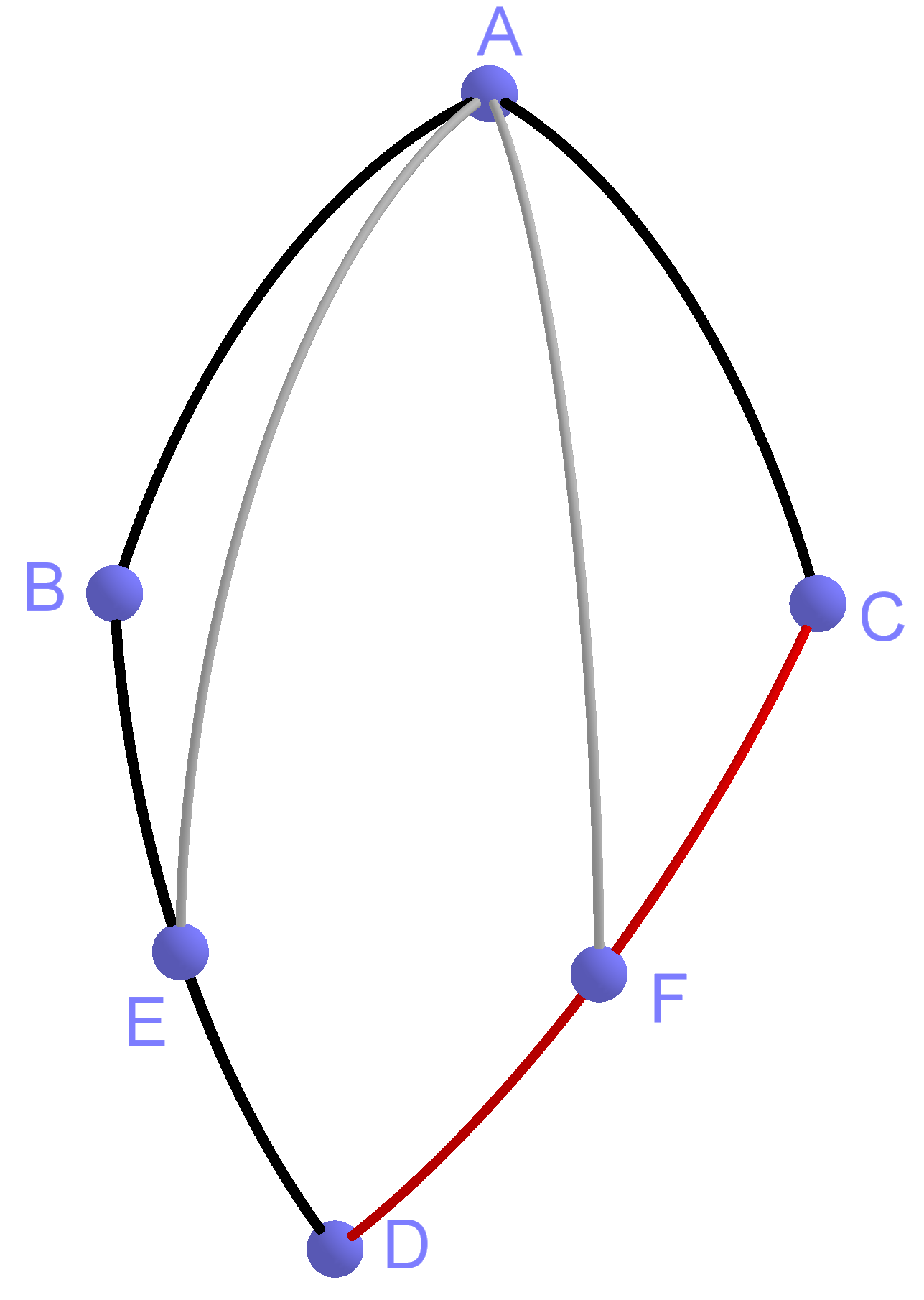}
    	\put(28,88){\tiny $A$}
    	\put(-7,52){\tiny $B$}
    	\put(63,52){\tiny $C$}
    	\put(-2,23){\tiny $E$}
    	\put(44,20){\tiny $F$}
    	\put(11,0){\tiny $D$}
    \end{overpic}   \hspace{20pt}
	\begin{overpic}[scale=0.2]{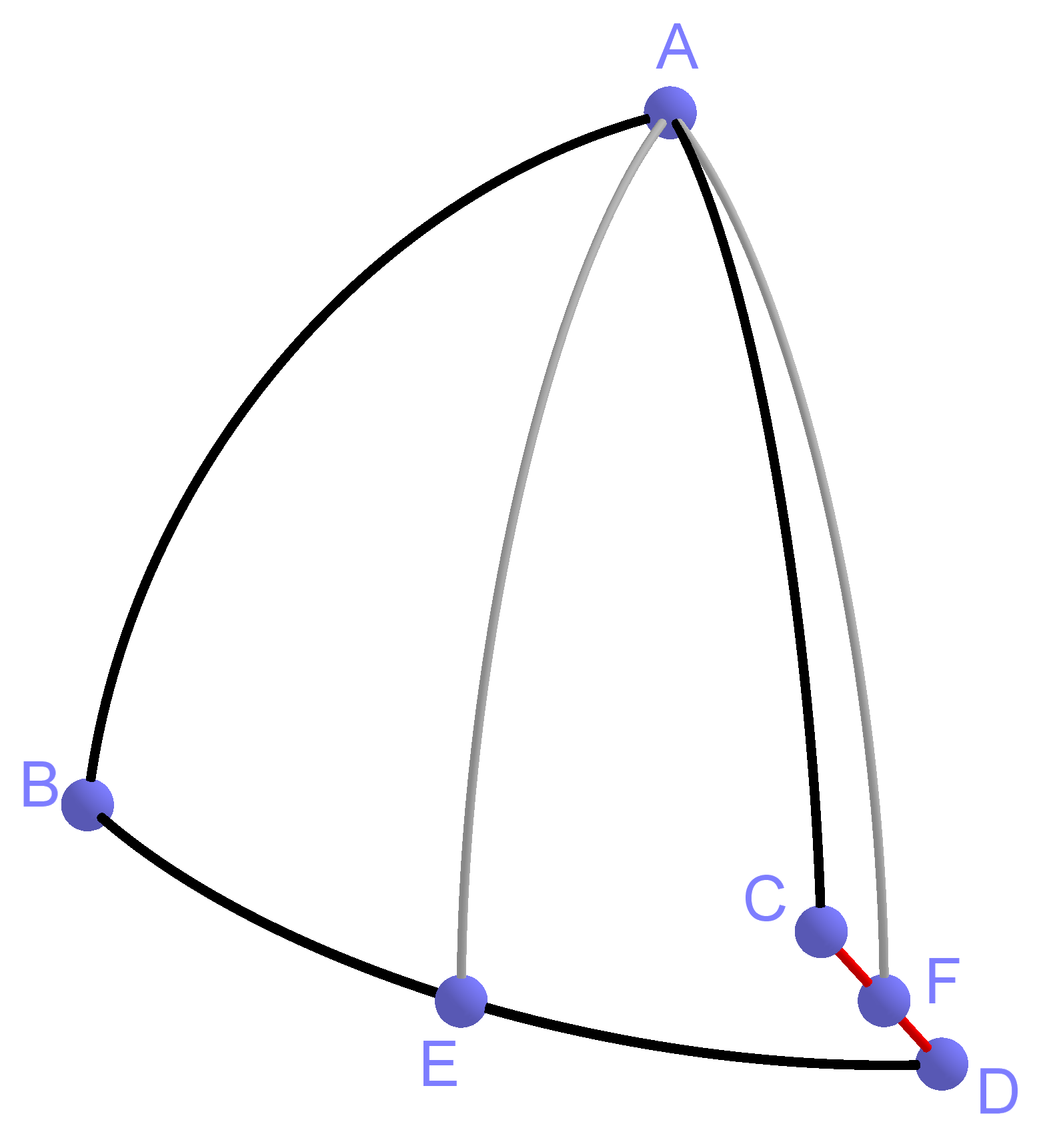}
		\put(55,85){\tiny $A$}
		\put(-8,25){\tiny $B$}
		\put(64,15){\tiny $C$}
		\put(42,10){\tiny $E$}
		\put(70,6){\tiny $F$}
		\put(88,0){\tiny $D$}
	\end{overpic}
	
	\caption{Construct $ABCD$ with fixed $A,E,F$ for $4$ typical positions of $D$.} \label{abd}
\end{figure}

The moduli of $2$-layer earth map tilings with $f=2n\ge6$ tiles was described in \cite{lpwx} Theorem $17$. Any tile can be constructed as shown in Fig.\,\ref{abd}. Fix a $\triangle AEF$ with $AE = AF = \frac{1}{2}$ and $EF = \frac{2}{f}$.
The quadrilateral is then determined by the location of $D$: extend $DE$ to $B$, such that $E$ is the middle point of $BD$; extend $DF$ to $C$, such that $F$ is the middle point of $CD$; then connect $A$ to $B$, $C$ to form the quadrilateral $\square ABDC$. Fig.\,\ref{abd} shows four typical positions of $D$. Then the moduli is all locations of $D$, such that the boundary of $\square ABDC$ has no self intersection, i.e. the interior of the union of two half lunes $\triangle AEF \cup \triangle A'PQ$  in Fig.\,\ref{modular}, where $A'$ is the antipode of $A$ and $P,Q$ are on the great arc $EF$ with $EP=FQ=\frac{1}{2}$.

The reduction of the quadrilateral from type $a^2 bc$ to type  $a^3 b$ ($a=b$ or equivalently $a=c$) is represented by the dotted curve $UVW$ inside the moduli, which is determined by the equation $BD=AB$. Let $t$ be the angle between two longitudes $AEA'$ and $ADA'$. Then $t\in (-\frac12,\frac12+\frac2f)$. Note that 
$AB=AC=BD=A'D=a$, $CD=b$, $AD=1-a$. Using the angle notations for type $a^3b$, denote  $\angle BDC=\aaa$, $\angle ABD=\bbb$, $\angle BAC=\ccc=\frac4f$, $\angle ACD=\ddd$. By the sine and cosine law for $\triangle ADE$ and $\triangle ABD$, we get 
\begin{align}
	\cos t=\frac{1}{2\sin\frac{a}{2}}   \label{a3b-2}
\end{align}
\begin{align}
	\cos \bbb=-\frac{\cos a}{1-\cos a}=-\cos(2t).  \label{a3b-1}
\end{align}
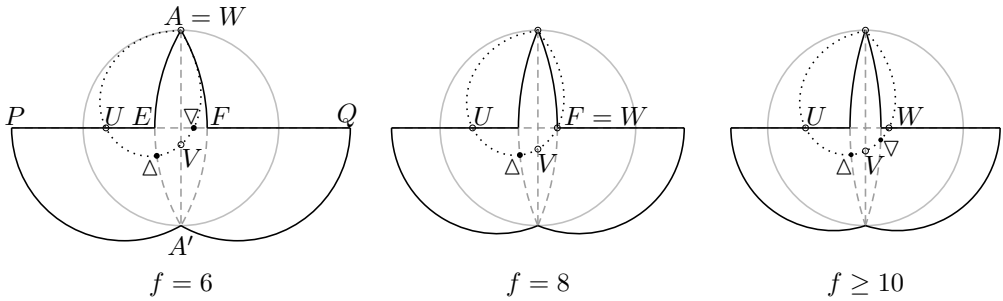
\begin{figure}[htp]
	\centering
	
	\begin{tikzpicture}[scale=0.65,line width=0.6pt,>=stealth]      
		
		\begin{scope}[xshift=0 cm] 
			\draw[gray!50] (0,0) circle (2);

			\draw[dotted]  plot[smooth,samples=558,domain=-60:60] ({2*(sqrt(4*cos(\x)*cos(\x)-1)/(2*cos(\x)*cos(\x)))*sin(\x-180/6)/(1+sqrt(4*cos(\x)*cos(\x)-1)/(2*cos(\x)*cos(\x))*cos(\x-180/6))},{2*	(1-2*cos(\x)*cos(\x))/(2*cos(\x)*cos(\x))/(1+(sqrt(4*cos(\x)*cos(\x)-1)/(2*cos(\x)*cos(\x)))*cos(\x-180/6))});
			

			\draw ({-2*sin(180/6)/(1+cos(180/6))},0) arc(180:150:{2*sin(180/6)/(1-cos(180/6)*cos(180/6))});
			\draw ({2*sin(180/6)/(1+cos(180/6))},0) arc(0:30:{2*sin(180/6)/(1-cos(180/6)*cos(180/6))});
			
			\draw[densely dashed, gray!75] ({-2*sin(180/6)/(1+cos(180/6))},0) arc(180:210:{2*sin(180/6)/(1-cos(180/6)*cos(180/6))});
			\draw[densely dashed, gray!75] ({2*sin(180/6)/(1+cos(180/6))},0) arc(0:-30:{2*sin(180/6)/(1-cos(180/6)*cos(180/6))});
			
			\draw[densely dashed, gray!75] (0,2)--(0,-2)
			({-2*cos(180/6)/(1-sin(180/6))},0)--({2*cos(180/6)/(1-sin(180/6))},0);

			\draw({-2*cos(180/6)/(1-sin(180/6))},0) arc(180:300:{2*cos(180/6)/(1-sin(180/6)*sin(180/6))});
			\draw ({2*cos(180/6)/(1-sin(180/6))},0) arc(0:-120:{2*cos(180/6)/(1-sin(180/6)*sin(180/6))});
			
			\draw({-2*cos(180/6)/(1-sin(180/6))},0)--({-2*sin(180/6)/(1+cos(180/6))},0)
			({2*sin(180/6)/(1+cos(180/6))},0)--({2*cos(180/6)/(1-sin(180/6))},0);

			\node at (0.5,2.3) {\small $A=W$};
			\node at (0,-2.4) {\small $A'$};
			\node at (-0.8,0.25) {\small $E$};
			\node at (0.8,0.25) {\small $F$};
			\node at (-3.4,0.25) {\small $P$};
			\node at (3.4,0.25) {\small $Q$};
			
			\node at (0,-3.2) {\small $f=6$};
			\node at (-0.63,-0.82) {\small $\vartriangle$};
			\node at (0.21,0.22) {\small $\triangledown$};
			\node at (-1.35,0.25) {\small $U$};
			\node at (0.23,-0.615) {\small $V$};

			\draw[line width=0.1] ({-2*cos(180/12)/(sin(180/12)+1)},0) circle (0.06);
			\fill ({-2*sqrt(3)/7},-4/7) circle (0.06);
			\draw[line width=0.1] (0,{-2/(2*sqrt(2)+3)}) circle (0.06);
			\fill ({2*cos(5*180/12)/(sin(5*180/12)+1)},0) circle (0.06);
			\draw[line width=0.1] (0,2) circle (0.06);
		\end{scope}
		\begin{scope}[xshift=7.3 cm] 
			\draw[gray!50] (0,0) circle (2);

			\draw[dotted]  plot[smooth,samples=558,domain=-60:60] ({2*(sqrt(4*cos(\x)*cos(\x)-1)/(2*cos(\x)*cos(\x)))*sin(\x-180/8)/(1+sqrt(4*cos(\x)*cos(\x)-1)/(2*cos(\x)*cos(\x))*cos(\x-180/8))},{2*	(1-2*cos(\x)*cos(\x))/(2*cos(\x)*cos(\x))/(1+(sqrt(4*cos(\x)*cos(\x)-1)/(2*cos(\x)*cos(\x)))*cos(\x-180/8))});
			

			\draw ({-2*sin(180/8)/(1+cos(180/8))},0) arc(180:157.3:{2*sin(180/8)/(1-cos(180/8)*cos(180/8))});
			\draw ({2*sin(180/8)/(1+cos(180/8))},0) arc(0:22.7:{2*sin(180/8)/(1-cos(180/8)*cos(180/8))});
			
			\draw[densely dashed, gray!75] ({-2*sin(180/8)/(1+cos(180/8))},0) arc(180:202.7:{2*sin(180/8)/(1-cos(180/8)*cos(180/8))});
			\draw[densely dashed, gray!75] ({2*sin(180/8)/(1+cos(180/8))},0) arc(0:-22.7:{2*sin(180/8)/(1-cos(180/8)*cos(180/8))});
			
			\draw[densely dashed, gray!75] (0,2)--(0,-2)
			({-2*cos(180/8)/(1-sin(180/8))},0)--({2*cos(180/8)/(1-sin(180/8))},0);
			
			\draw({-2*cos(180/8)/(1-sin(180/8))},0) arc(180:292.4:{2*cos(180/8)/(1-sin(180/8)*sin(180/8))});
			\draw ({2*cos(180/8)/(1-sin(180/8))},0) arc(0:-112.4:{2*cos(180/8)/(1-sin(180/8)*sin(180/8))});
			
			\draw({-2*cos(180/8)/(1-sin(180/8))},0)--({-2*sin(180/8)/(1+cos(180/8))},0)
			({2*sin(180/8)/(1+cos(180/8))},0)--({2*cos(180/8)/(1-sin(180/8))},0);

			\node at (-0.52,-0.8) {\small $\vartriangle$};
			\node at (-1.1,0.25) {\small $U$};
			\node at (1.4,0.25) {\small$F=W$};
			\node at (0.2,-0.71) {\small $V$};
			\node at (0,-3.2) {\small $f=8$};
			\draw[line width=0.1] ({-2*cos(180/8)/(sin(180/8)+1)},0) circle (0.06);
			\fill ({(-2*sqrt(3)*cos(3*180/8))/(2+sqrt(3)*cos(180/8))},{-2/(2+sqrt(3)*cos(180/8))}) circle (0.06);
			\draw[line width=0.1] (0,{(2-4*cos(180/8)*cos(180/8))/(2*cos(180/8)*cos(180/8)+sqrt(4*cos(180/8)*cos(180/8)-1))}) circle (0.06);
			\draw[line width=0.1] ({2*cos(3*180/8)/(sin(3*180/8)+1)},0) circle (0.06);
			\draw[line width=0.1] (0,2) circle (0.06);
		\end{scope}
		\begin{scope}[xshift=14 cm]

			\draw[gray!50] (0,0) circle (2);
			
			\draw[dotted]  plot[smooth,samples=558,domain=-60:60] ({2*(sqrt(4*cos(\x)*cos(\x)-1)/(2*cos(\x)*cos(\x)))*sin(\x-180/10)/(1+sqrt(4*cos(\x)*cos(\x)-1)/(2*cos(\x)*cos(\x))*cos(\x-180/10))},{2*	(1-2*cos(\x)*cos(\x))/(2*cos(\x)*cos(\x))/(1+(sqrt(4*cos(\x)*cos(\x)-1)/(2*cos(\x)*cos(\x)))*cos(\x-180/10))});	    	
			
			\draw ({-2*sin(180/10)/(1+cos(180/10))},0) arc(180:162:{2*sin(180/10)/(1-cos(180/10)*cos(180/10))});
			\draw ({2*sin(180/10)/(1+cos(180/10))},0) arc(0:18:{2*sin(180/10)/(1-cos(180/10)*cos(180/10))});
			
			\draw[densely dashed, gray!75] ({-2*sin(180/10)/(1+cos(180/10))},0) arc(180:198:{2*sin(180/10)/(1-cos(180/10)*cos(180/10))});
			\draw[densely dashed, gray!75]  ({2*sin(180/10)/(1+cos(180/10))},0) arc(0:-18:{2*sin(180/10)/(1-cos(180/10)*cos(180/10))});
			
			\draw[densely dashed, gray!75] (0,2)--(0,-2)
			({-2*cos(180/10)/(1-sin(180/10))},0)--({2*cos(180/10)/(1-sin(180/10))},0);
			
			\draw({-2*cos(180/10)/(1-sin(180/10))},0) arc(180:288:{2*cos(180/10)/(1-sin(180/10)*sin(180/10))});
			\draw ({2*cos(180/10)/(1-sin(180/10))},0) arc(0:-108:{2*cos(180/10)/(1-sin(180/10)*sin(180/10))});
			
			\draw({-2*cos(180/10)/(1-sin(180/10))},0)--({-2*sin(180/10)/(1+cos(180/10))},0)
			({2*sin(180/10)/(1+cos(180/10))},0)--({2*cos(180/10)/(1-sin(180/10))},0);
			
			\node at (-0.44,-0.785) {\small $\vartriangle$};
			\node at (0.53,-0.43) {\small $\triangledown$};
			\node at (-1,0.25) {\small $U$};
			\node at (0.92,0.25) {\small $W$};
			\node at (0.2,-0.75) {\small $V$};
			\node at (0,-3.2) {\small $f\ge10$};    	
			\draw[line width=0.1] ({-2*cos(3*180/20)/(sin(3*180/20)+1)},0) circle (0.06);
			\fill ({(-2*sqrt(3)*cos(2*180/5))/(2+sqrt(3)*sin(2*180/5))},{-2/(2+sqrt(3)*sin(2*180/5))}) circle (0.05);	  
			\draw[line width=0.1] (0,{(2-4*cos(180/10)*cos(180/10))/(2*cos(180/10)*cos(180/10)+sqrt(4*cos(180/10)*cos(180/10)-1))}) circle (0.06);  	
			\fill ({(2*cos(2*180/5)*sqrt(4*cos(180/5)*cos(180/5)-1))/(2*cos(180/5)*cos(180/5)+sin(2*180/5)*sqrt(4*cos(180/5)*cos(180/5)-1))},{(2-4*cos(180/5)*cos(180/5))/(2*cos(180/5)*cos(180/5)+sin(2*180/5)*sqrt(4*cos(180/5)*cos(180/5)-1))}) circle (0.05);
			\draw[line width=0.1] ({2*cos(7*180/20)/(sin(7*180/20)+1)},0) circle (0.06);
			\draw [line width=0.1] (0,2) circle (0.06);
		\end{scope}
	\end{tikzpicture}
	
	\caption{Moduli of $2$-layer earth map tilings.}\label{modular}
\end{figure}

The equation \eqref{a3b-2} determines the dotted curve $UVW$ representing the moduli of   $2$-layer earth map $a^3b$-tilings. The shape of this curve is unchanged for all even $f\ge6$. However, its position is getting closer to the middle as $f\to\infty$.
We remark that the point $V$ is the intersection of $UVW$ with the dashed line $AA'$, corresponding to a rhombus (type $a^4$) with $t=\frac1f$ and $\bbb=1-\frac2f$. 

\begin{figure}[htp]
	\centering	
	\begin{tikzpicture}[scale=0.4]   
%
		
		\node at (-19-3,0) {\small  $\cos a=-\frac{\cos \bbb}{1-\cos \bbb}$};
		\node at (-20.8-1.1,-2) {\small $\lim\limits_{f\to\infty}b=a\in[\frac13,\frac12)$};
		\node at (-21-0.6,-4) {\small  $\lim\limits_{f\to\infty}\aaa=\lim\limits_{f\to\infty}\ccc=0$};
		\node at (-21-1.6,-6) {\small  $\lim\limits_{f\to\infty}\ddd=2-\bbb$};
		
	\end{tikzpicture}\hspace{30pt}
    \begin{overpic}[scale=0.195]{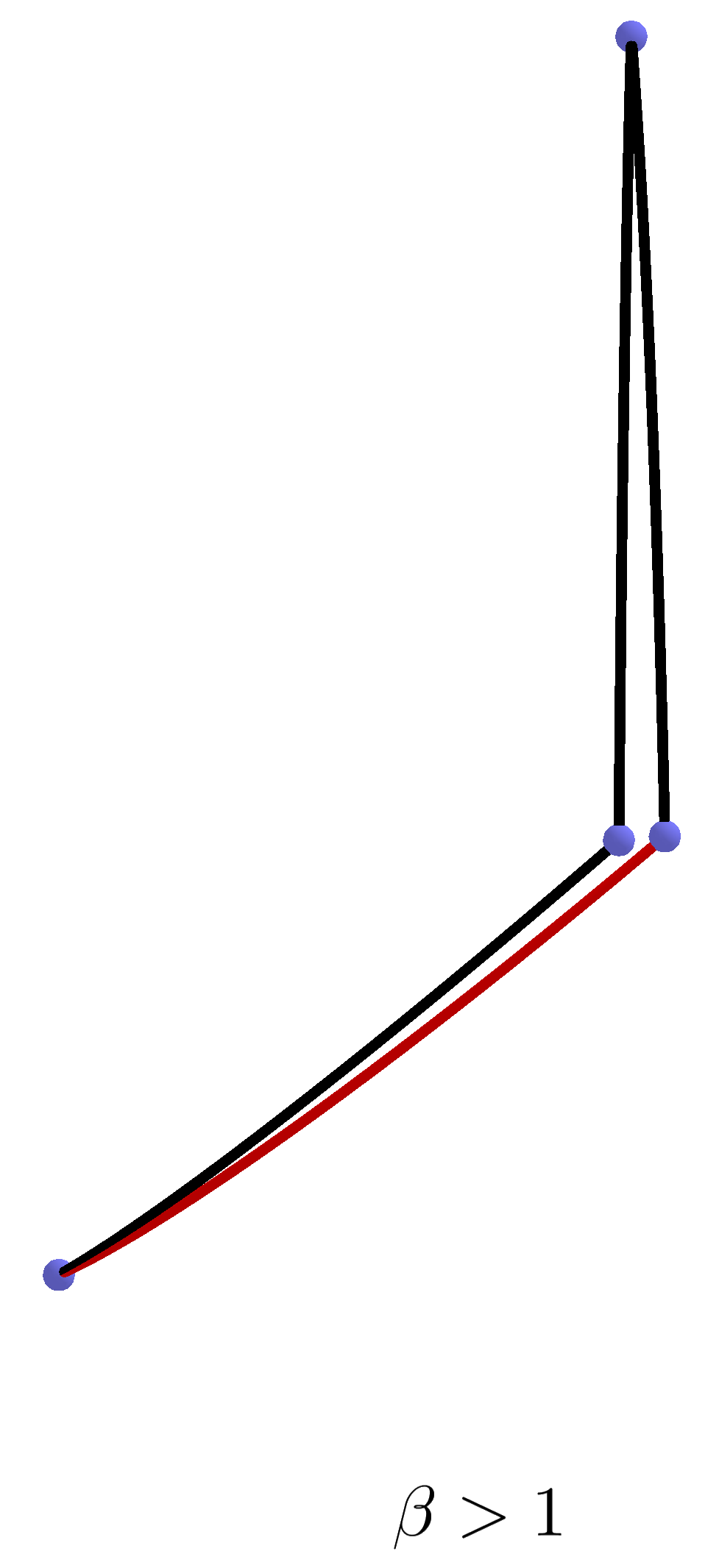}
    	\put(15,2){\small $\bbb>1$}
    \end{overpic}\hspace{20pt}			
    \begin{overpic}[scale=0.187]{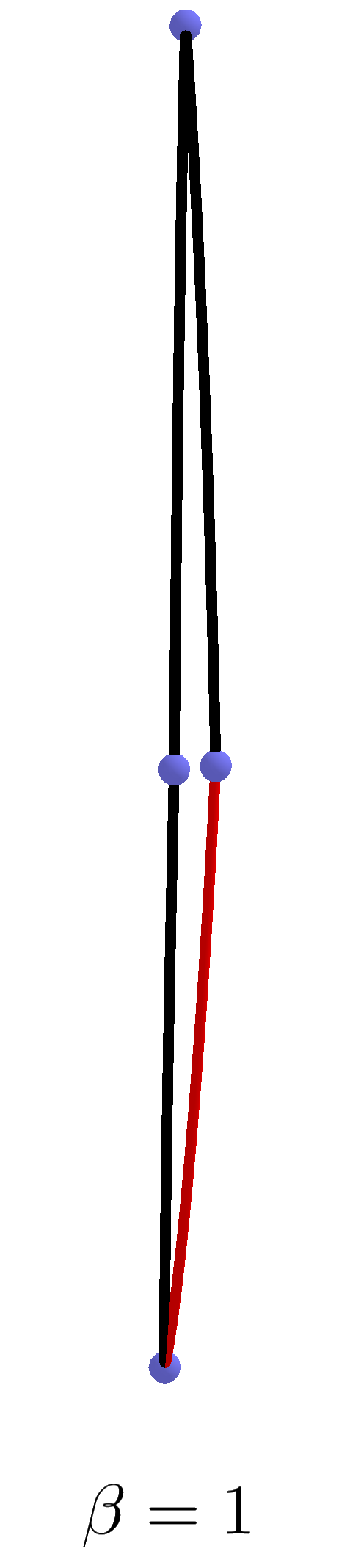} 
    	\put(5,2){\small $\bbb=1$}
    \end{overpic}\hspace{20pt}
    \begin{overpic}[scale=0.188]{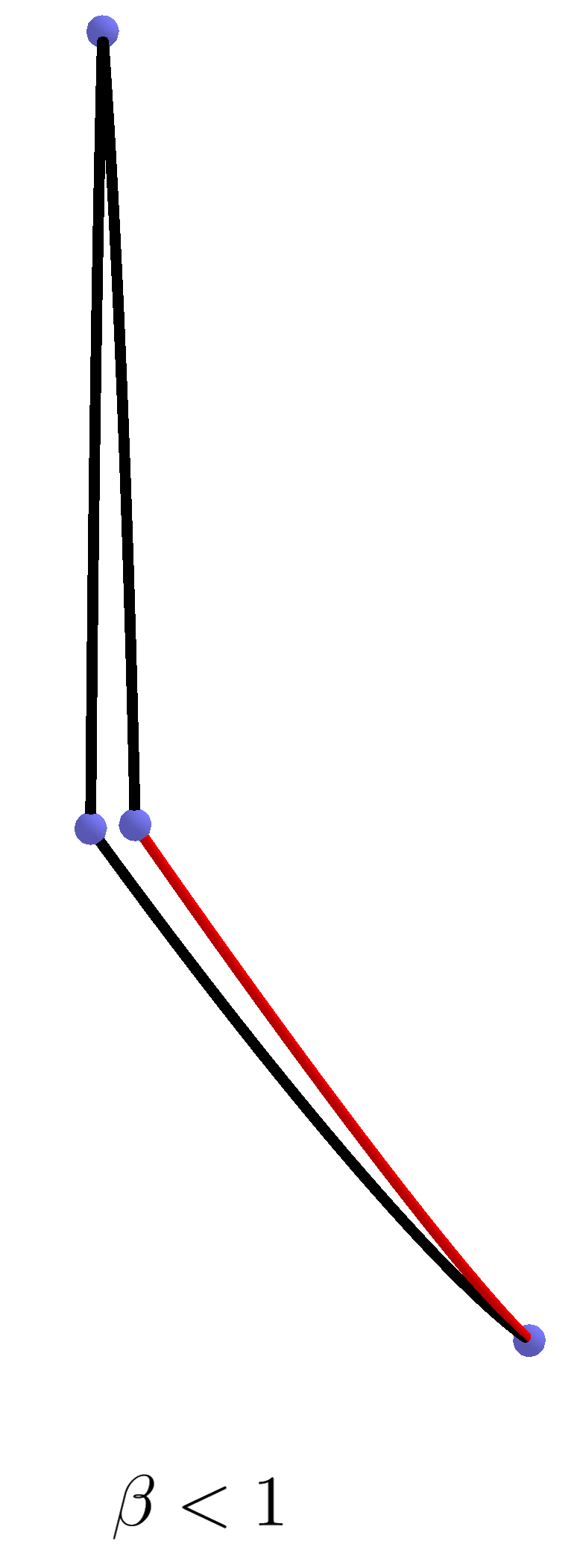}    
    	\put(0,2){\small $\bbb<1$}
    \end{overpic}\hspace{20pt}

	\caption{$a^3b$-quadrilaterals for $2$-layer earth map tilings when $f\to\infty$.}\label{3-limit}
\end{figure}

By \eqref{a3b-1}, we have $\bbb=1-2t$. Applying Lemma \ref{calculate-1}, we get exact formulas for all geometric data $a,b,\aaa,\ddd$ in Table \ref{tab-2}.
For $f=6$, we have $-\frac14<t<\frac13$, $\frac13<\bbb<\frac32$, $\frac13\le a<1$, $b_{\text{min}}\le b<\frac76$ and $0<\aaa<\frac53$. Here the minimum length $b_{\text{min}}=\arccos\frac{59+73\sqrt{73}}{768}\approx0.1514$ is reached when $\bbb=\arctan\frac{5\sqrt{3}+\sqrt{219}}{2}\approx0.4729$.
For $f=8$, we have $-\frac14<t<\frac14$, $\frac12<\bbb<\frac32$, $\frac13\le a<\frac12$ , $0<b<1$ and $0<\aaa<1$.
For $f\ge10$, we have $-\frac14<t<\frac14$, $\frac12<\bbb<\frac32$, $\frac13\le a<\frac12$ , $b_{\text{min}}\le b<\frac12+\frac4f$ and $0<\aaa<1$.  The exact formula for the minimum length $b_{\text{min}}$ is too long to be included here.
The $a^3b$-quadrilateral becomes thinner and thinner as $f\to\infty$ for some fixed $\bbb\in(\frac12,\frac32)$ in Fig.\,\ref{3-limit}.

\subsection*{Basic flips of special $2$-layer earth map $a^3b$-tilings}
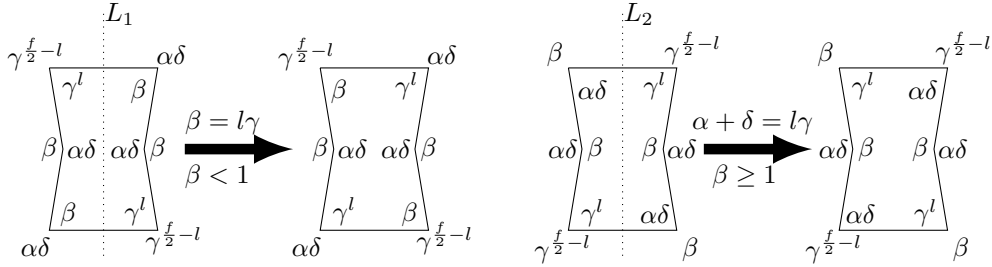
\begin{figure}[htp]
	\centering		
	\begin{tikzpicture}[>=latex,scale=0.3]
	%
	%
	%
	%
	%
	%
	%
	%
	%
	%
	%

	\begin{scope}[xshift=0 cm, yshift=0 cm, scale=1.2]
	\draw(-2,3)--(2,3)--(1.5,0)--(2,-3)--(-2,-3)--(-1.5,0)--(-2,3);
	
	\draw[line width=4pt, ->](3,0)--(7,0);
	
	\node at (4.4,1){\small $\bbb=l\ccc$};
	\node at (4.2,-1){\small $\bbb<1$};
	
	\draw[dotted] (0,5)--(0,-4);
	\node at (0.6,5){\small $L_1$};
	
	\node at (-1.1,2.35){\small $\ccc^{l}$};
	\node at (1.2,-2.3){\small $\ccc^{l}$};
	
	\node at (1.3,2.15){\small $\bbb$};
	\node at (-1.3,-2.4){\small $\bbb$};
	
	\node at (-0.8,0){\small $\aaa\ddd$};
	\node at (0.8,0){\small $\aaa\ddd$};
	
	\node at (-2,0){\small $\bbb$};
	\node at (2,0){\small $\bbb$};
	
	\node at (2.5,3.5){\small $\aaa\ddd$};
	\node at (-2.5,-3.6){\small $\aaa\ddd$};
	\node at (-2.5,3.5){\small $\ccc^{\frac{f}{2}-l}$};
	\node at (2.6,-3.3){\small $\ccc^{\frac{f}{2}-l}$};	    			
	\end{scope}	
	
	\begin{scope}[xshift=12 cm, yshift=0 cm, scale=1.2]
	\draw(-2,3)--(2,3)--(1.5,0)--(2,-3)--(-2,-3)--(-1.5,0)--(-2,3);
	
	\node at (-1.3,2.15){\small $\bbb$};
	\node at (1.3,-2.4){\small $\bbb$};
	
	\node at (1.2,2.35){\small $\ccc^{l}$};
	\node at (-1.1,-2.3){\small $\ccc^{l}$};
	
	\node at (-0.8,0){\small $\aaa\ddd$};
	\node at (0.8,0){\small $\aaa\ddd$};
	
	\node at (-2,0){\small $\bbb$};
	\node at (2,0){\small $\bbb$};
	
	\node at (2.5,3.5){\small $\aaa\ddd$};
	\node at (-2.5,-3.6){\small $\aaa\ddd$};
	\node at (-2.5,3.5){\small $\ccc^{\frac{f}{2}-l}$};
	\node at (2.6,-3.3){\small $\ccc^{\frac{f}{2}-l}$};
	\end{scope}	
	
	\begin{scope}[xshift=23 cm, yshift=0 cm, scale=1.2]
	\draw(-2,3)--(2,3)--(1.5,0)--(2,-3)--(-2,-3)--(-1.5,0)--(-2,3);
	
	\draw[dotted] (0,5)--(0,-4);
	\node at (0.6,5){\small $L_2$};
	
	\node at (1.2,2.35){\small $\ccc^{l}$};
	\node at (-1.1,-2.35){\small$\ccc^{l}$};
	
	\node at (-1.1,2.2){\small $\aaa\ddd$};
	\node at (1.2,-2.4){\small $\aaa\ddd$};
	
	\node at (-1,0){\small $\bbb$};
	\node at (1,0){\small $\bbb$};
	
	\node at (-2.2,0){\small $\aaa\ddd$};
	\node at (2.2,0){\small $\aaa\ddd$};
	
	\node at (-2.5,3.5){\small $\bbb$};
	\node at (2.5,-3.7){\small $\bbb$};
	\node at (2.6,3.75){\small$\ccc^{\frac{f}{2}-l}$};
	\node at (-2.2,-3.6){\small $\ccc^{\frac{f}{2}-l}$};
	\draw[line width=4pt, ->](3,0)--(7,0);
	
	\node at (4.8,1){\small $\aaa+\ddd=l\ccc$};
	\node at (4.5,-1){\small$\bbb\ge1$};
	\end{scope}	
	
	\begin{scope}[xshift=35 cm, yshift=0 cm, scale=1.2]
	\draw(-2,3)--(2,3)--(1.5,0)--(2,-3)--(-2,-3)--(-1.5,0)--(-2,3);
	
	\node at (1.1,2.2){\small $\aaa\ddd$};
	\node at (-1.2,-2.4){\small $\aaa\ddd$};
	
	\node at (-1.1,2.35){\small $\ccc^{l}$};
	\node at (1.2,-2.35){\small $\ccc^{l}$};
	
	\node at (-1,0){\small $\bbb$};
	\node at (1,0){\small $\bbb$};
	
	\node at (-2.2,0){\small $\aaa\ddd$};
	\node at (2.2,0){\small $\aaa\ddd$};
	
	\node at (-2.5,3.5){\small $\bbb$};
	\node at (2.5,-3.7){\small $\bbb$};
	\node at (2.6,3.75){\small $\ccc^{\frac{f}{2}-l}$};
	\node at (-2.2,-3.6){\small $\ccc^{\frac{f}{2}-l}$};
	\end{scope}

	%
	%
	%
	%
	%
	%

	
	\end{tikzpicture}
	
	\caption{Two basic flip modifications for certain $2$-layer earth map tilings.} \label{flip1}
\end{figure}

\begin{figure}[htp]
	\centering
	\begin{tikzpicture}[>=latex,scale=0.25] 
		\foreach \a in {0,1}
		{
			\begin{scope}[xshift=3*\a cm] 
				\draw (4,0)--(4,-9)
				(7,0)--(7,-9);
				\draw[line width=1.5] (7,-2)--(4,-7);
			\end{scope}
		}
		\foreach \b in {0}
		{
			\fill[gray!50]
			(0,0) -- (4,0) -- (4,-9)-- (0,-9)--(0,0);	
			
			\begin{scope}[xshift=6*\b cm] 
				\draw (0,0)--(0,-9)
				(4,0)--(4,-9)
				(0,-2)--(3,-2)--(3,-4)--(4,-7)--(2,-5)
				(0,-2)--(2,-4)
				(0,-2)--(1,-5)--(1,-7)--(4,-7);
				\draw[line width=1.5] (3,-2)--(4,-2)
				(2,-4)--(3,-4)
				(1,-5)--(2,-5)
				(0,-7)--(1,-7);

				\fill (0.56,-3.21) circle (0.05);
				\fill (0.68,-3.13) circle (0.05);
				\fill (0.80,-3.05) circle (0.05);
				
				\fill (0.56+2.55,-3.21-2.6) circle (0.05);
				\fill (0.68+2.55,-3.13-2.6) circle (0.05);
				\fill (0.80+2.55,-3.05-2.6) circle (0.05);

				\fill (0,-2) circle (0.15); \fill (4,-7) circle (0.15);
				
				\fill (10,-7) circle (0.04);
			\end{scope}
		}
		\fill (11,-5) circle (0.05);
		\fill (11.3,-5) circle (0.05);
		\fill (11.6,-5) circle (0.05);
		
		\foreach \c in {0,1}
		{

			\begin{scope}[xshift=7*\c cm+15 cm] 
				
				\fill[gray!50]
				(0,0) -- (4,0) -- (4,-9)-- (0,-9)--(0,0);
				\draw (0,0)--(0,-9)
				(4,0)--(4,-9)
				(0,-2)--(3,-2)--(3,-4)--(4,-7)--(2,-5)
				(0,-2)--(2,-4)
				(0,-2)--(1,-5)--(1,-7)--(4,-7);
				\draw[line width=1.5] (3,-2)--(4,-2)
				(2,-4)--(3,-4)
				(1,-5)--(2,-5)
				(0,-7)--(1,-7);

				\fill (0.56,-3.21) circle (0.05);
				\fill (0.68,-3.13) circle (0.05);
				\fill (0.80,-3.05) circle (0.05);
				
				\fill (0.56+2.55,-3.21-2.6) circle (0.05);
				\fill (0.68+2.55,-3.13-2.6) circle (0.05);
				\fill (0.80+2.55,-3.05-2.6) circle (0.05);

				\fill (0,-2) circle (0.15); \fill (4,-7) circle (0.15);

			\end{scope}
		}
		\fill (27-0.3,-5) circle (0.05);
		\fill (27.3-0.3,-5) circle (0.05);
		\fill (27.6-0.3,-5) circle (0.05);
		\foreach \l in {0}
		{
			\begin{scope}[xshift=3*\l cm+15 cm] 
				\draw (4,0)--(4,-9)
				(7,0)--(7,-9);
				\draw[line width=1.5] (7,-2)--(4,-7);
			\end{scope}
		}
		
		\foreach \d in {0,1}
		{
			
			\begin{scope}[xshift=7*\d cm+31 cm] 
				
				\fill[gray!50]
				(0,0) -- (4,0) -- (4,-9)-- (0,-9)--(0,0);
				\draw (0,0)--(0,-9)
				(4,0)--(4,-9)
				(0,-2)--(3,-2)--(3,-4)--(4,-7)--(2,-5)
				(0,-2)--(2,-4)
				(0,-2)--(1,-5)--(1,-7)--(4,-7);
				\draw[line width=1.5] (3,-2)--(4,-2)
				(2,-4)--(3,-4)
				(1,-5)--(2,-5)
				(0,-7)--(1,-7);

				\fill (0.56,-3.21) circle (0.05);
				\fill (0.68,-3.13) circle (0.05);
				\fill (0.80,-3.05) circle (0.05);
				
				\fill (0.56+2.55,-3.21-2.6) circle (0.05);
				\fill (0.68+2.55,-3.13-2.6) circle (0.05);
				\fill (0.80+2.55,-3.05-2.6) circle (0.05);

				\fill (0,-2) circle (0.15); \fill (4,-7) circle (0.15);

			\end{scope}
		}
		
		\begin{scope}[xshift=48 cm] 
			
			\fill[gray!50]
			(0,0) -- (4,0) -- (4,-9)-- (0,-9)--(0,0);
			\draw (0,0)--(0,-9)
			(4,0)--(4,-9)
			(0,-2)--(3,-2)--(3,-4)--(4,-7)--(2,-5)
			(0,-2)--(2,-4)
			(0,-2)--(1,-5)--(1,-7)--(4,-7);
			\draw[line width=1.5] (3,-2)--(4,-2)
			(2,-4)--(3,-4)
			(1,-5)--(2,-5)
			(0,-7)--(1,-7);

			\fill (0.56,-3.21) circle (0.05);
			\fill (0.68,-3.13) circle (0.05);
			\fill (0.80,-3.05) circle (0.05);
			
			\fill (0.56+2.55,-3.21-2.6) circle (0.05);
			\fill (0.68+2.55,-3.13-2.6) circle (0.05);
			\fill (0.80+2.55,-3.05-2.6) circle (0.05);

			\fill (0,-2) circle (0.15); \fill (4,-7) circle (0.15);

		\end{scope}
		
		\foreach \r in {0,1}
		{
			\begin{scope}[xshift=3*\r cm+38 cm] 
				\draw (4,0)--(4,-9)
				(7,0)--(7,-9);
				\draw[line width=1.5] (7,-2)--(4,-7);
			\end{scope}
		}

		\draw (31+4,0)--(31+4,-9)
		(31+7,0)--(31+7,-9);
		\draw[line width=1.5] (31+7,-2)--(31+4,-7);

		\fill (48+5-0.3,-5) circle (0.05);
		\fill (48+5.3-0.3,-5) circle (0.05);
		\fill (48+5.6-0.3,-5) circle (0.05);
		
	\end{tikzpicture}
	\caption{Applying $1,2$ or $3$ basic flips for special $2$-layer earth map tilings} \label{a1}
\end{figure}
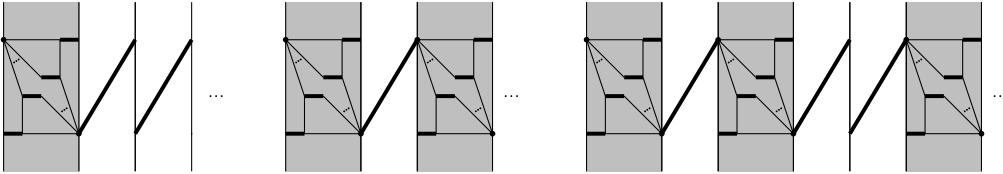

As explained in \cite{lw} Section $1$, the quadrilaterals admitting any flip modifications are special points on the curve $UVW$ when $\bbb$ is an integer multiple of $\ccc$.  There are two types of basic flip modifications in Fig.\,\ref{flip1}, depending on $\bbb<1$ or $\bbb\ge1$. In Table \ref{tab-3},  $\aaa+\ddd=m\ccc$ is equivalent to $\bbb=(\frac{f}{2}-m)\ccc$ for some integer $m\ge2$ and $f\ge8$, and each basic flip modification flips $l=\min\{m,\frac{f}{2}-m\}$ time zones ($2l$ tiles). Then this  basic flip may be applied several times as long as there are still $l$ continuous time zones left, as shown in Fig.\,\ref{a1}. By $\ccc=\frac4f$ and $\bbb\in(\frac12,\frac32)$,  we have $m\in(\frac f8,\frac{3f}{8})$, and such $m,f$ also  determine the quadrilateral uniquely.

If $m=\frac f4$, then $\bbb=1$, $D=\vartriangle$ in Fig.\,\ref{modular} and the quadrilateral reduces to a triangle. If $m=\frac{f+2}{4}$, then $\bbb=1-\frac2f$, $D=V$ in Fig.\,\ref{modular} and the quadrilateral reduces to a rhombus. If $m=\frac{f+4}{4}$, then $\ddd=1$, $D=\triangledown$ in Fig.\,\ref{modular} and the quadrilateral reduces to a triangle. The quadrilateral is convex when $D$ is on the curve $UVW$ between $\vartriangle$ and $\triangledown$.  We conclude that convex quadrilaterals in $2$-layer earth map $a^3b$-tilings never admit any flips since there are no other integers between $\frac f4$, $\frac{f+4}{4}$.

The tilings deduced in Fig.\,\ref{abd-3} and \ref{abd-5} can be viewed as applying the basic flips in Fig.\,\ref{flip1} several times. 
By the AVC $\{(f-2n)\aaa\bbb\ddd,2\aaa^n\ccc^{\frac{f}{2}-mn}\ddd^n,2n\bbb\ccc^m\}$ in Fig.\,\ref{abd-3}, we get $\frac{f}{2}-mn\ge0$. Then $m>\frac f8$ implies that $n\le3$.
We may apply at most $n=3$ basic flips along $L_2$ in Fig.\,\ref{flip1} simultaneously when $\frac f8<m\le\frac{f}{6}$; and apply at most $n=2$ such basic flips simultaneously when $\frac f6<m\le\frac f4$. These are the first three cases in Table \ref{tab-3}.
By the AVC $\{(f-2n)\aaa\bbb\ddd,2n\,\aaa\ccc^m\ddd,2\bbb^n\ccc^{\frac{f}{2}-mn}\}$ in Fig.\,\ref{abd-5}, we get $\frac{f}{2}-(\frac{f}{2}-m)n\ge0$ (to be consistent with $\aaa+\ddd=m\ccc$ in this section, we need to switch $m$ and $\frac{f}{2}-m$).  Then $m<\frac {3f}{8}$ implies that $n\le3$.
We may apply at most $n=3$ basic flips along $L_1$ in Fig.\,\ref{flip1} simultaneously when $\frac f3\le m<\frac{3f}{8}$; and apply at most $n=2$ such basic flips simultaneously when $\frac {f+4}{4}\le m<\frac f3$. These are the last two cases in Table \ref{tab-3}.

Counting $a^3b$-quadrilaterals on the moduli curve $UVW$ admitting flips, there are totally $\mathcal{Q}_{1}(f)=2\lfloor \frac{f-4}{8}\rfloor+1$ of them. By the number of rational ones $\mathcal{Q}_{2}(f)$ deduced from \cite{lw}, we get $\mathcal{Q}_{3}(f)=\mathcal{Q}_{1}(f)-\mathcal{Q}_{2}(f)$ for the number of general ones (i.e. with some irrational angle) shown in Table \ref{Tab-1.3}.

In the last column of Table \ref{tab-3}, the total number of different tilings by $1,2$ or $3$ simultaneous flips are counted using unordered integer partitions, where the function $\langle x \rangle$ gives the closest integer to $x$.

\section{Tilings with two degree $3$ vertex types}
\label{sec-two3}

\begin{proposition}\label{33}
	In an $a^3b$-tiling with some irrational angle, if there exist two degree $3$ vertex types different from $\aaa\bbb\ddd$ and $\aaa\ccc\ddd$, then the tiling is $T(6\aaa^2\bbb,2\ccc^3,6\bbb\ccc\ddd^2)$ with $12$ tiles. 
\end{proposition}

\begin{proof}
	
	If there are two degree 3 vertex types with only $a$-edges, we always have $\bbb=\ccc$, contradicting Proposition $9$ in \cite{lw}. If there are two degree $3$ vertex types with a $b$-edge, by Lemma \ref{geometry2}, we just need to consider $\{\aaa^2\bbb,\ccc\ddd^2\}$. All other cases will have one vertex with a $b$-edge and the other with only $a$-edges. Up to the symmetry of interchanging $\aaa\leftrightarrow\ddd$ and $\bbb\leftrightarrow\ccc$, we will choose one from   $\{\aaa^2\bbb,\bbb\ddd^2\}$ and one from $\{\bbb^3,\bbb^2\ccc,\bbb\ccc^2,\ccc^3\}$, giving $8$ combinations of two degree $3$ vertices. However $\{\aaa^2\bbb,\bbb\ccc^2\}$ implies $\aaa=\ccc$ and $\ddd=\frac{4}{f}<1$, contradicting $\ddd=1$ in Lemma \ref{geometry3}. The case $\{\bbb^3,\bbb\ddd^2\}$ implies $\bbb=\ddd=\frac{2}{3}$. By Lemma \ref{geometry3}, $\aaa=1>\ddd$. By Lemma \ref{anglesum}, $\ccc=\frac{4}{f}-\frac{1}{3}\le\frac13<\bbb$, contradicting Lemma \ref{geometry1}. The other seven cases are summarized in Table \ref{Tab-2}. 
	
	The AAD $\aaa^2\bbb=\thin^{\aaa}\bbb^{\ccc}\thin^{\bbb}\aaa^{\ddd}\thick^{\ddd}\aaa^{\bbb}\thin$ implies $\bbb\ccc\cdots$ and $\ddd^2\cdots$ are vertices.
	The AAD $\bbb\ddd^2=\thin^{\aaa}\bbb^{\ccc}\thin^{\ccc}\ddd^{\aaa}\thick^{\aaa}\ddd^{\ccc}\thin$ implies $\aaa^2\cdots,\aaa\ccc\cdots$ and $\ccc^2\cdots$ are vertices.

	\begin{table*}[htp]                        
		\centering     
		      
		~\\ 
		\resizebox{\textwidth}{18mm}{\begin{tabular}{c|c|c|c}	 			
				Cases & Angles & Irrational angle lemma & Conclusion \\
				\hline  
				$\aaa^2\bbb,\bbb^3$&$(\frac23,\frac23,\ccc,\frac23-\ccc+\frac4f)$, $0<\ccc<\frac23+\frac4f$&$n_3=n_4$&No $\bbb\ccc\cdots$, a contradiction \\
				\hline 
				$\aaa^2\bbb,\bbb^2\ccc$&$(\aaa,2-2\aaa,-2+4\aaa,2-3\aaa+\frac4f)$, $\frac12<\aaa<\frac23+\frac{4}{3f}$&$n_1+4n_3=2n_2+3n_4$&No $\ddd^2\cdots$ , a contradiction\\			
				\hline 
				$\aaa^2\bbb,\ccc^3$&$(\aaa,2-2\aaa,\frac23,-\frac23+\aaa+\frac4f)$, $\frac23-\frac4f<\aaa<1$&$n_1+n_4=2n_2$&$\text{AVC}\sub\{\aaa^2\bbb,\ccc^3,\bbb\ccc\ddd^2,\bbb^2\ddd^4\}$\\
				\hline 
				\multirow{1}{*}{$\aaa^2\bbb,\ccc\ddd^2$}&$(\aaa,2-2\aaa,-2+2\aaa+\frac8f,2-\aaa-\frac4f)$, $1-\frac4f<\aaa<1$&$n_1+2n_3=2n_2+n_4$&$\text{AVC}\sub\{\aaa^2\bbb,\ccc\ddd^2,\bbb^{\frac f4}\ccc^{\frac f4}\}$\\
				\hline
				$\bbb^2\ccc,\bbb\ddd^2$&$(-1+\frac{3\bbb}{2}+\frac4f,\bbb,2-2\bbb,1-\frac\bbb2)$, $\frac23-\frac{8}{3f}<\bbb<1$&$3n_1+2n_2=4n_3+n_4$ &No $\aaa^2\cdots$, a contradiction\\
				\hline 
				$\bbb\ccc^2,\bbb\ddd^2$&$(\frac4f,\bbb,1-\frac\bbb2,1-\frac\bbb2)$, $0<\bbb<2$& $2n_2=n_3+n_4$&  No $\aaa\ccc\cdots$, a contradiction \\ 	
				\hline 
				$\bbb\ddd^2,\ccc^3$&$(\aaa,\frac23-2\aaa+\frac8f,\frac23,\frac23+\aaa-\frac4f)$, $0<\aaa<\frac13+\frac4f$& $n_1+n_4=2n_2$&$\text{AVC}\sub\{\bbb\ddd^2,\ccc^3,\aaa^2\bbb\ccc\}$\\
				\hline
		\end{tabular}}
	\caption{Two degree 3 vertex types without $\aaa\bbb\ddd$, $\aaa\ccc\ddd$.}\label{Tab-2}  
	\end{table*}

	\subsubsection*{Case \{$\aaa^2\bbb,\bbb^3$\}}
	
	The AAD of $\aaa^2\bbb$ gives $\bbb\ccc\cdots$, denoted by $\bn$ with $n_2,n_3\ge1$. By $n_4=n_3\ge1$, Lemma \ref{anglesum} and Parity Lemma, we get $n_1=0$ and $n_3=n_4\ge2$, contradicting $\bbb+2\ccc+2\ddd=2+\frac8f>2$.
	
	\subsubsection*{Case \{$\aaa^2\bbb,\bbb^2\ccc$\}}
	The AAD of $\aaa^2\bbb$ gives $\ddd^2\cdots$, denoted by $\bn$, $n_4\ge2$.
	
	If $n_1=0$, by $4n_3=2n_2+3n_4\ge6$, then $n_3\ge2$. By $\bbb+2\ccc+2\ddd=2+\frac8f>2$, we get $n_2=0$. Therefore, $4n_3=3n_4$. Then $n_3\ge3$ and $n_4\ge4$, contradicting $3\ccc+4\ddd=2+\frac{16}{f}>2$.
	
	If $n_1=1$, by Parity Lemma, we get $n_4\ge3$. By $1+4n_3=2n_2+3n_4$, we have $n_3\ge2$, contradicting $\aaa+2\ccc+3\ddd=2+\frac{12}{f}>2$.
	
	If $n_1\ge2$, by $n_1+4n_3=2n_2+3n_4\ge6$, we get either $n_3\ge1$ or $n_1\ge6$, contradicting $2\aaa+\ccc+2\ddd=2+\frac{8}{f}>2$ and $6\aaa+2\ddd=4+\frac{8}{f}>2$. 
	

	\subsubsection*{Case \{$\aaa^2\bbb,\ccc^3$\}}
	We will find all possible vertices $\bn$. 
	
	If $n_2=0$, then $n_1=n_4=0$. By $\ccc^3$, we get $\bn=(0\,0\,3\,0)$.
	
	If $n_2=1$, by $n_1+n_4=2$, then $\bn=(2\,1\,n_3\,0),(1\,1\,n_3\,1)$ or $(0\,1\,n_3\,2)$.  By $\aaa^2\bbb$, $\bbb+2\ccc+2\ddd>2$ and Lemma \ref{anglesum}, we get $\bn=(2\,1\,0\,0)$ or $(0\,1\,1\,2)$. 
	
	If $n_2=2$, by $\aaa^2\bbb$ and $n_1+n_4=4$, we deduce that $\bn=(1\,2\,n_3\,3)$ or $(0\,2\,n_3\,4)$. By $\aaa+2\bbb+3\ddd>2$ and $2\bbb+\ccc+4\ddd>2$, we get $\bn=(0\,2\,0\,4)$. 
	
	If $n_2\ge3$, by $\aaa^2\bbb$ and $n_1+n_4=2n_2$, we get either $n_1=1,n_4\ge5$ or $n_1=0,n_4\ge6$, contradicting $\aaa+3\bbb+5\ddd>2$ and $3\bbb+6\ddd>2$. 
	
	Therefore, $\text{AVC}\sub\{\aaa^2\bbb,\ccc^3,\bbb\ccc\ddd^2,\bbb^2\ddd^4\}$. 	
	If $\bbb^2\ddd^4$ is a vertex, then we get the AAD $\bbb^{\aaa}\thin^{\aaa}\bbb\cdots,\bbb^{\aaa}\thin^{\ccc}\bbb\cdots$ or $\bbb^{\aaa}\thin^{\ccc}\ddd\cdots$ at $\bbb^2\ddd^4$. This gives a vertex $\aaa\thin\aaa\cdots$ or $\aaa\ccc\cdots$, which is not in the AVC. Therefore, $\ddd\cdots=\bbb\ccc\ddd^2$. By Balance Lemma, $\bbb\ccc\ddd^2$ is a vertex. By Lemma \ref{anglesum}, $f=12$. There is only one solution satisfying Balance Lemma: $\text{AVC}=\{6\aaa^2\bbb,2\ccc^3,6\bbb\ccc\ddd^2\}$.

	In Fig.\,\ref{a2b,c3}, the unique AAD $\thin^\ddd\ccc_1^\bbb\thin^\ddd\ccc_{1'}^\bbb\thin^\ddd\ccc_{1''}^\bbb\thin$ determines $T_1,T_{1'},T_{1''}$. Then $\aaa_1\cdots=\aaa_1\aaa_2\bbb_3$ determines $T_2$. Then $\aaa_{1'}\cdots=\aaa_{1'}\aaa_{2'}\bbb_{3'}$ determines $T_{2'}$. Then $\bbb_1\ddd_{1'}\ddd_{2'}\cdots=\bbb_1\ccc_3\ddd_{1'}\ddd_{2'}$ determines $T_3$. Then $\bbb_{1'}\ddd_{1''}\cdots=\bbb_{1'}\ccc_{3'}\ddd_{1''}\cdots$ determines $T_{3'}$. Then $\aaa_3\bbb_2\cdots=\aaa_3\aaa_4\bbb_2$ determines $T_4$. Then $\ccc_{2'}\ddd_3\ddd_4\cdots=\bbb_{4'}\ccc_{2'}\ddd_3\ddd_4$ and $\aaa_{3'}\bbb_{2'}\cdots=\aaa_{3'}\aaa_{4'}\bbb_{2'}$ determine $T_{4'}$. Similarly, we can determine $T_{2''},T_{3''},T_{4''}$. We get a unique earth map tiling $T(6\aaa^2\bbb,2\ccc^3,6\bbb\ccc\ddd^2)$.
	\begin{figure}[htp]
		\centering
		\begin{tikzpicture}[>=latex,scale=0.6]
			
			\foreach \b in {0,1,2}
			{
				\begin{scope}[xshift=2*\b cm]	
					\draw (0,-2)--(0,-6)
					(-2,0)--(-2,-4)
					(0,0)--(0,-2)--(-1,-2)--(-1,-4)--(-2,-4)--(-2,-6);
					
					\draw[line width=1.5]	
					(-2,-2)--(-1,-2)
					(-1,-4)--(0,-4);

					\node at (-1,-0.25) {\small $\ccc$};
					\node at (-1,-1.75) {\small $\aaa$};
					\node at (-0.25,-1.7) {\small $\bbb$};
					\node at (-1.75,-1.65) {\small $\ddd$};
					\node at (-0.75,-2.4) {\small $\bbb$};
					\node at (-1.25,-2.3) {\small $\aaa$};
					\node at (-0.25,-2.35) {\small $\ccc$};
					\node at (-1.75,-2.35) {\small $\ddd$};
					
					\node at (-0.75,-3.75) {\small $\aaa$};
					\node at (-1.25,-3.75) {\small $\bbb$};
					\node at (-0.25,-3.65) {\small $\ddd$};
					\node at (-1.75,-3.7) {\small $\ccc$};
					
					\node at (-0.25,-4.35) {\small $\ddd$};
					\node at (-1.75,-4.35) {\small $\bbb$};
					\node at (-1,-4.35) {\small $\aaa$};
					\node at (-1,-5.8) {\small $\ccc$};

				\end{scope}
			}

			\node[draw,shape=circle, inner sep=0.5] at (-1,-1) {\small $1$};
			\node[draw,shape=circle, inner sep=0.5] at (1,-1) {\small $1'$};
			\node[draw,shape=circle, inner sep=0.5] at (-0.5,-3) {\small $3$};
			\node[draw,shape=circle, inner sep=0.5] at (0.5,-3) {\small $2'$};
			\node[draw,shape=circle, inner sep=0.5] at (-1.5,-3) {\small $2$};
			\node[draw,shape=circle, inner sep=0.5] at (1.5,-3) {\small $3'$};
			\node[draw,shape=circle, inner sep=0.5] at (-1,-5) {\small $4$};
			\node[draw,shape=circle, inner sep=0.5] at (1,-5) {\small $4'$};
			
			\node[draw,shape=circle, inner sep=0.5] at (-1+4,-1) {\small $1''$};
			\node[draw,shape=circle, inner sep=0.5] at (-0.5+4,-3) {\small $3''$};
			\node[draw,shape=circle, inner sep=0.5] at (-1.5+4,-3) {\small $2''$};
			\node[draw,shape=circle, inner sep=0.5] at (-1+4,-5) {\small $4''$};

		\end{tikzpicture}\hspace{30pt}
		\caption{$f=12$, $T(6\aaa^2\bbb,  2\ccc^3, 6\bbb\ccc\ddd^2)$.} 
		\label{a2b,c3}
	\end{figure}

	By Lemma \ref{calculate-1}, we get
	\[
		\aaa=1-\arcsin{\tfrac{\sqrt6}{4}},\,\,\,
		a=\arccos{\tfrac{2\sqrt5-3}{3}},\,\,\,
		b=\arccos{\left(3\sqrt5-6\right)}
	\]
	$\aaa\approx 0.7902,a\approx 0.3367,b\approx 0.2495$.
			
	\subsubsection*{Case \{$\aaa^2\bbb,\ccc\ddd^2$\}}
	
	We will find all possible vertices $\bn$.
	
	If $n_1=0$, by $2n_3=2n_2+n_4$, $\ccc\ddd^2$ and Parity Lemma, we deduce that $n_3>0$, $n_4=0$ or $2$. Then $\bn=(0\,0\,1\,2)$ or $(0\,\frac{f}{4}\,\frac{f}{4}\,0)$. 
	
	If $n_3=0$, by $n_1=2n_2+n_4$ and $\aaa^2\bbb$, we deduce that $n_1\ge2$, $n_2=0$ or $1$. When $n_2=0$, we get $n_1=n_4\ge2$, contradicting  $\aaa+\ddd=2-\frac4f\ge\frac43$. When $n_2=1$, we get $\bn=(2\,1\,0\,0)$.

	If $n_1n_3\neq0$, then Lemma \ref{anglesum} implies that either $n_2=0,n_4\ge3$ contradicting $\ccc\ddd^2$, or $n_4=0,n_2\ge2$ contradicting $\aaa^2\bbb$ by $n_1=2n_2-2n_3\ge2$.
	
	Therefore, $\text{AVC}\subset\{\aaa^2\bbb,\ccc\ddd^2,\bbb^\frac{f}{4}\ccc^\frac{f}{4}
	\}$. There is only one solution satisfying Balance Lemma: $\{\frac{f}{2}\,\aaa^2\bbb,\frac{f}{2}\,\ccc\ddd^2,2\bbb^\frac{f}{4}\ccc^\frac{f}{4}\}$.
	By Lemma \ref{geometry4} and Lemma \ref{geometry4}$'$, we get $\frac{24}{f}=3\bbb+3\ccc>2$. Therefore, $f<12$. Then $f=8,\aaa+\ddd=\frac32$ and $\bbb+\ccc=1$. By Lemma \ref{relation}, $\aaa+\ddd<1+\bbb$ and $\aaa+\ddd<1+\ccc$, contradicting $2\aaa+2\ddd=3=2+\bbb+\ccc$. 
	
	\subsubsection*{Case \{$\bbb\ddd^2,\bbb^2\ccc$\}}
	
	The AAD of $\bbb\ddd^2$ gives $\aaa^2\cdots$, denoted by $\bn$, $n_1\ge2$.
	
	If $n_4=0$, by $4n_3=3n_1+2n_2\ge6$, then $n_3\ge2$. By $2\aaa+\bbb+2\ccc>2$, we get $n_2=0$. Therefore, $3n_1=4n_3$. Then $n_1\ge4$ and $n_3\ge3$, contradicting $4\aaa+3\ccc>2$. 
	
	If $n_4=1$, by Parity Lemma, we get $n_1\ge3$. By $3n_1+2n_2=4n_3+1$, we get $n_3\ge2$, contradicting $3\aaa+2\ccc+\ddd>2$.
	
	If $n_4\ge2$, by $3n_1+2n_2=4n_3+n_4$, we get either $n_3\ge1$ or $n_4\ge6$, contradicting $2\aaa+\ccc+2\ddd>2$ and $2\aaa+6\ddd>2$.
	
	
	\subsubsection*{Case \{$\bbb\ddd^2,\bbb\ccc^2$\}}
	The AAD of $\bbb\ddd^2$ gives $\aaa\ccc\cdots$, denoted by $\bn$ with  $n_1,n_3\ge1$. By $2n_2=n_3+n_4$, we get $n_2\ge1$. By Lemma \ref{anglesum} and Parity Lemma, we get $n_4=0$, $n_3=2n_2\ge2$ and $n_1\ge2$, contradicting $2\aaa+\bbb+2\ccc>2$.

	\subsubsection*{Case \{$\bbb\ddd^2,\ccc^3$\}}
	 By Lemma \ref{geometry4}, we have $2\aaa+\bbb>1$.
	 We will find all possible vertices $\bn$.
	 
	 If $n_2=0$, then $n_1=n_4=0$. By $\ccc^3$, we get $\bn=(0\,0\,3\,0)$.
	 
	 If $n_2=1$, by $n_1+n_4=2$, then $\bn=(2\,1\,n_3\,0),(1\,1\,n_3\,1)$ or $(0\,1\,n_3\,2)$.  By $\bbb\ddd^2$, $2\aaa+\bbb+2\ccc>2$ and Lemma \ref{anglesum}, we get $\bn=(2\,1\,1\,0)$ or $(0\,1\,0\,2)$.

	 If $n_2\ge2$, by $\bbb\ddd^2$ and $n_1+n_4=2n_2$, we get either $n_4=1,n_1\ge3$ or $n_4=0,n_1\ge4$, contradicting $3\aaa+2\bbb+\ddd>2$ and $4\aaa+2\bbb>2$. 
	 
	 Therefore, $\text{AVC}\sub\{\bbb\ddd^2,\ccc^3,\aaa^2\bbb\ccc\}$. By Balance Lemma, $\aaa^2\bbb\ccc$ is a vertex. Then we get $f=12$ and $0<\aaa<\tfrac23$. If it satisfies \eqref{4-7}, we have $\sin\aaa \sin(\aaa+\frac\pi3)=0$.  If it satisfies \eqref{4-8}, we have $\sin^2(\aaa+\frac\pi3)=0$. Both contradict $0<\aaa<\tfrac23$. 
\end{proof}
\begin{remark}
	By the classification of rational case in \cite{lw}, Proposition \ref{33} still holds for $a^3b$-quadrilaterals with arbitrary angles. All other tilings with $\ge2$ different degree $3$ vertex types must involve $\aaa\bbb\ddd$ or $\aaa\ccc\ddd$, and are related to $2$-layer earth map tilings by Proposition \ref{111}. We pick all such tilings out of Table \ref{tab-2} and \ref{tab-3} and obtain the complete list in Table \ref{Tab-33}.
			
	\begin{table*}[htp]                        
		\centering     
	  
		~\\ 
		\begin{tabular}{c|c|c}	 
			
			Vertex type & f & tilings \\
			\hline 
			$\aaa\bbb\ddd,\bbb^3$& $6k(k\ge2)$& $T((6k-6)\aaa\bbb\ddd,2\bbb^3,6\aaa\ccc^k\ddd)$\\
			\hline 
			\multirow{11}{*}{$\aaa\bbb\ddd,\bbb\ccc^2$}&\multirow{2}{*}{$8$}&$T(6\aaa\bbb\ddd,2\bbb\ccc^2,2\aaa\ccc^2\ddd)$\\
			&& $T(4\aaa\bbb\ddd,4\bbb\ccc^2,2\aaa^2\ddd^2)$\\
			\cline{2-3}
			&\multirow{3}{*}{$10$}& $T(8\aaa\bbb\ddd,2\bbb\ccc^2,2\aaa\ccc^3\ddd)$\\ &&$T(6\aaa\bbb\ddd,4\bbb\ccc^2,2\aaa^2\ccc\ddd^2)$\\ &&$\{4\aaa\bbb\ddd,4\bbb\ccc^2,2\aaa^2\bbb\ccc,2\aaa\ddd^3\}$: $3$ tilings\\
			\cline{2-3}
			&\multirow{3}{*}{$12$}& $T(10\aaa\bbb\ddd,2\bbb\ccc^2,2\aaa\ccc^4\ddd)$\\ &&$\{8\aaa\bbb\ddd,4\bbb\ccc^2,2\aaa^2\ccc^2\ddd^2\}$: $2$ tilings\\
			&&$T(6\aaa\bbb\ddd,6\bbb\ccc^2,2\aaa^3\ddd^3)$\\
			\cline{2-3}
			&\multirow{3}{*}{$14$}&$T(12\aaa\bbb\ddd,2\bbb\ccc^2,2\aaa\ccc^5\ddd)$\\ &&$\{10\aaa\bbb\ddd,4\bbb\ccc^2,2\aaa^2\ccc^3\ddd^2\}$: $2$ tilings\\ 
			&&$T(8\aaa\bbb\ddd,6\bbb\ccc^2,2\aaa^3\ccc\ddd^3)$\\
			\hline 
			$\aaa\bbb\ddd,\ccc^3$&$6$&$T(6\aaa\bbb\ddd,2\ccc^3)$\\
			\hline 
			$\aaa^2\bbb,\ccc^3$&$12$&$T(6\aaa^2\bbb,2\ccc^3,6\bbb\ccc\ddd^2)$\\
			\hline
		\end{tabular}
		\caption{All $a^3b$-tilings with two degree $3$ vertex types.}\label{Tab-33}      
	\end{table*}
	
	By \cite{lw}, there are only six rational quadrilaterals in Table \ref{Tab-33}: 		
    $(6,3,4,3)/6$, $(1,8,4,3)/6$, $(12,4,6,2)/9$, $(3,10,6,5)/9$ with $f=6$;      
    $(1,6,2,3)/5$ with $f=10$; 	
    $(1,2,1,3)/3$ with $f=12$.  
\end{remark}

\section{Tilings with a unique degree $3$ vertex type}
\label{sec-one3}

By \eqref{vcountf}, all vertices have degree $3$ when $f=6$. The standard spherical cube is the only tiling with a unique degree $3$ vertex type when $f=6$, and it is of type $a^4$ and all angles are $\frac23$. For tilings of type $a^3b$ with a unique degree $3$ vertex type, we will always assume $f\ge8$ henceforth.

By Proposition \ref{111}, we only need to consider $\aaa^2\bbb,\bbb^3,\bbb\ccc^2$ or $\bbb\ddd^2$ to be the unique degree $3$ vertex type, up to the symmetry of exchanging $\aaa\leftrightarrow\ddd$ and $\bbb\leftrightarrow\ccc$. Each of these four cases turns out to produce a unique tiling with irrational angles. 

\begin{lemma} \label{tildev21}
	In an $a^3b$-tiling, if all degree $3$ vertices have no $b$-edge, then there exists some degree $4$ vertex with two $b$-edges: one of $\{\aaa^4,\aaa^3\ddd,\aaa^2\ddd^2,\aaa\ddd^3,\ddd^4\}$ must appear, and  $\aaa^2\cdots$, $\ddd^2\cdots$ are both vertices.
\end{lemma}

\begin{proof}
	
	Let $\tilde{v}_k$ be the number of degree $k$ vertices and $\tilde{f}$ be the number of hexagons obtained by removing all $b$-edges in the tiling. Similar counting via Euler's formula for this combinatorial hexagonal tiling gives: 
	
	\begin{equation}
		f=2\tilde{f}=4+\sum_{k=3}^{\infty}(k-2)\tilde{v}_k
		=4+\tilde{v}_3+2\tilde{v}_4+3\tilde{v}_5+\cdots, \label{vcountf'} 
	\end{equation}
	\begin{equation}
		\tilde{v}_{2}=6+\sum_{k\ge 4}^{\infty}(k-3)\tilde{v}_k
		=6+\tilde{v}_4+2\tilde{v}_5+3\tilde{v}_6+\cdots, \label{vcountv'} 
	\end{equation}
	which implies $\tilde{v}_2\ge6$. Such lowest degree $2$ vertices only come from two cases in Fig.\,\ref{tildev2}: an old degree 3 vertex with one $b$-edge or an old degree 4 vertex with two $b$-edges, which is one of $\{\aaa^4,\aaa^3\ddd,\aaa^2\ddd^2,\aaa\ddd^3,\ddd^4\}$. Then by Balance Lemma, both $\aaa^2\cdots$ and $\ddd^2\cdots$ must also appear.	
	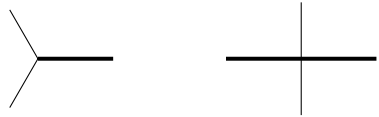
\begin{figure}[htp]
		\centering
		\begin{tikzpicture}[>=latex,scale=0.5] 
			
			\draw (0,0)--(-0.75,-1.3)
			(0,0)--(-0.75,1.3)
			(7,1.5)--(7,-1.5);
			
			\draw[line width=1.5] (0,0)--(2,0)
			(5,0)--(9,0);
			
		\end{tikzpicture}
		\caption{Two cases for $\tilde{v}_2$.}
		\label{tildev2}
	\end{figure}		
\end{proof}

\begin{lemma} \label{a2b4}
	In an $a^3b$-tiling, if all degree $3$ vertices are $\aaa^2\bbb$, then there exists a degree $4$ vertex without $\aaa$: one of  $\{\bbb^4,\bbb^3\ccc,\bbb^2\ccc^2,\bbb^2\ddd^2,\bbb\ccc^3,\bbb\ccc\ddd^2,\ccc^4,\ccc^2\ddd^2,\ddd^4\}$ must appear.
\end{lemma}
\begin{proof}
	The counting formula \eqref{vcountf} and \eqref{vcountv} implies
	\[
	2v_3+v_4\\=16+v_4+2v_5+4v_6+\cdots\ge f+10. 
	\]
	However, there are $2v_3$ $\aaa$'s in all degree $3$ vertices $\aaa^2\bbb$, which must be $\le f$. Then $v_4\ge10$. If all degree $4$ vertices have $\aaa$, we get $f= \# \aaa \ge 2v_3+v_4$, a contradiction. Then one of  $\{\bbb^4,\bbb^3\ccc,\bbb^2\ccc^2,\bbb^2\ddd^2,\bbb\ccc^3,\bbb\ccc\ddd^2,\ccc^4,\ccc^2\ddd^2,\ddd^4\}$ must appear by Parity Lemma. 
\end{proof}

\begin{proposition}
	In an $a^3b$-tiling with some irrational angle, if all degree $3$ vertices are $\bbb\ccc^2$, then the tiling is $T(8\bbb\ccc^2,6\aaa^2\ddd^2,4\aaa\bbb^2\ddd)$ with $16$ tiles. 
\end{proposition}

\begin{proof}
	Lemma \ref{tildev21} implies that one of $\{\aaa^4,\aaa^3\ddd,\aaa^2\ddd^2,\aaa\ddd^3,\ddd^4\}$ must appear. Applying the Irrational Angle Lemma to each case, we get Table \ref{Tab-bc2}.

	\begin{table*}[htp]                        
		\centering     
	   
		~\\ 
		\resizebox{\textwidth}{14mm}{\begin{tabular}{c|c|c|c}	 
				
				Cases & Angles &Irrational angle Lemma& Conclusion \\
				\hline 
				$\bbb\ccc^2,\aaa^4$&$(\frac12,\bbb,1-\frac\bbb2,\frac12-\frac\bbb2+\frac4f)$, $0<\bbb<1+\frac8f$
				&$2n_2=n_3+n_4$ &No $\ddd^2\cdots$, a contradiction\\
				\hline 
				$\bbb\ccc^2,\aaa^3\ddd$&$(\frac12+\frac\bbb4-\frac2f,\bbb,1-\frac\bbb2,\frac12-\frac{3\bbb}{4}+\frac6f)$, $0<\bbb<\frac23+\frac8f$
				& $n_1+4n_2=2n_3+3n_4$&No $\ddd^2\cdots$, a contradiction\\
				\hline 
				$\bbb\ccc^2,\aaa^2\ddd^2$& $(\aaa,\frac8f,1-\frac4f,1-\aaa)$, $0<\aaa<1$
				&$n_1=n_4$ &$\text{AVC}\sub\{\bbb\ccc^2,\aaa^2\ddd^2,\aaa\bbb^2\ddd\}$\\
				\hline 
				$\bbb\ccc^2,\aaa\ddd^3$& $(\frac12-\frac{3\bbb}{4}+\frac6f,\bbb,1-\frac\bbb2,\frac12+\frac\bbb4-\frac2f)$, $0<\bbb<\frac23+\frac8f$
				&$3n_1+2n_3=4n_2+n_4$ & No $\aaa^2\cdots$, a contradiction\\
				\hline 
				$\bbb\ccc^2,\ddd^4$&$(\frac12-\frac\bbb2+\frac4f,\bbb,1-\frac\bbb2,\frac12)$, $0<\bbb<1+\frac8f$
				&$n_1+n_3=2n_2$ &No $\aaa^2\cdots$, a contradiction\\
				\hline 
		\end{tabular}}
		\caption{All degree $3$ vertices being $\bbb\ccc^2$.}\label{Tab-bc2}     
	\end{table*}

	\subsubsection*{Case \{$\bbb\ccc^2,\aaa^4$\}}
	
	Denote the vertex $\ddd^2\cdots$ by $\bn$, $n_4\ge2$. By $2n_2=n_3+n_4$, we have $n_2\ge1$.  
	
	If $n_2=1$, by $2=n_3+n_4$, then $n_3=0$ and $n_4=2$. By Parity Lemma and its degree being $>3$, we get $n_1\ge2$, contradicting $2\aaa+\bbb+2\ddd>2$.
	
	If $n_2\ge2$, by $2n_2=n_3+n_4$ and $\bbb\ccc^2$, then we get either $n_3=0,n_4\ge4$ or $n_3=1,n_4\ge3$, contradicting $2\bbb+4\ddd>2$ and $2\bbb+\ccc+3\ddd>2$.
	
			
	\subsubsection*{Case \{$\bbb\ccc^2,\aaa^3\ddd$\}}

	Denote the vertex $\ddd^2\cdots$ by $\bn$, $n_4\ge2$.
	
	If $n_1=0$, by $4n_2=2n_3+3n_4\ge6$, then $n_2\ge2$. By $2\bbb+\ccc+2\ddd>2$, we get $n_3=0$. Therefore, $4n_2=3n_4$. Then $n_2\ge3$ and $n_4\ge4$, contradicting $3\bbb+4\ddd>2$.
	
	If $n_1=1$, by Parity Lemma, then $n_4\ge3$. By $1+4n_2=2n_3+3n_4$, we get $n_2\ge2$. contradicting $\aaa+2\bbb+3\ddd>2$.
	
	If $n_1\ge2$, by $n_1+4n_2=2n_3+3n_4\ge6$, we get either $n_2\ge1$ or $n_1\ge6$,  contradicting $2\aaa+\bbb+2\ddd>2$ and $6\aaa+2\ddd>2$.
	

	\subsubsection*{Case \{$\bbb\ccc^2,\aaa^2\ddd^2$\}}

	By $f\ge8$, we get $\ccc\ge\frac12$. We will find all possible vertices $\bn$.
	
	If $n_1=0$, then $n_4=n_1=0$. By $\bbb\ccc^2$ and $\ccc\ge\frac12$, we get $\bn=(0\,\frac f4\,0\,0)$, $(0\,\frac{f+4}{8}\,1\,0)$, $(0\,1\,2\,0)$ or $(0\,0\,4\,0)$.
	
	If $n_1\ge1$, by $\aaa^2\ddd^2$, $n_1=n_4$ and Lemma \ref{anglesum}, we get $\bn=(2\,0\,0\,2)$, $(1\,\frac f8\,0\,1)$ or $(1\,0\,2\,1)$.

    Therefore, $\text{AVC}\subset\{\bbb\ccc^2,\aaa^2\ddd^2,\aaa\ccc^2\ddd,\ccc^4,\aaa\bbb^{\frac f8}\ddd,\bbb^{\frac f4},\bbb^{\frac {f+4}{8}}\ccc\}$. If $\aaa\ccc^2\ddd$ or $\ccc^4$ is a vertex, then $\bbb=1,\ccc=\frac12$ and $f=8$. We deduce that $\text{AVC}\sub\{\bbb\ccc^2,\aaa^2\ddd^2,\aaa\ccc^2\ddd,\ccc^4\}$. Then $\#\bbb<\#\ccc$, we get a contradiction. Then \[\text{AVC}\subset\{\bbb\ccc^2,\aaa^2\ddd^2,\aaa\bbb^{\frac f8}\ddd,\bbb^{\frac f4},\bbb^{\frac{f+4}{8}}\ccc\}.\]
	By the AVC, there is no $\aaa\ccc\cdots$ vertex. If $\bbb^{\frac f4}$ is a vertex, we have the unique AAD $\bbb^{\frac f4}=\thin\bbb^{\ccc}\thin^{\ccc}\bbb\thin\cdots$. This gives a vertex $\thin\ccc^{\bbb}\thin^{\bbb}\ccc\thin\cdots$. Then $\thin\ccc^{\bbb}\thin^{\bbb}\ccc\thin\cdots=\bbb\ccc^2=\thin^{\aaa}\bbb^{\ccc}\thin^{\ddd}\ccc^{\bbb}\thin^{\bbb}\ccc^{\ddd}\thin$. This gives a vertex $\ccc\ddd\cdots$, which is not in the AVC, a contradiction. Similarly, $\bbb^{\frac{f+4}{8}}\ccc$ is not a vertex and $\aaa\bbb^{\frac f 8}\ddd=\aaa\bbb^2\ddd=\thin^{\ccc}\bbb^{\aaa}\thin^{\aaa}\bbb^{\ccc}\thin^{\bbb}\aaa^{\ddd}\thick^{\aaa}\ddd^{\ccc}\thin$. Therefore, $\text{AVC}\sub\{\bbb\ccc^2,\aaa^2\ddd^2,\aaa\bbb^2\ddd\}$. By Balance Lemma,  $\aaa\bbb^2\ddd$ is a vertex, and we get $f=16$ by Lemma \ref{anglesum}. There is only one solution satisfying Balance Lemma: $\{ 8\bbb\ccc^2,6\aaa^2\ddd^2,4\aaa\bbb^2\ddd\}$.	
	
	In Fig.\,\ref{ab2,c2d2}, $\aaa_1\ddd_2\bbb_3\bbb_4$ determines $T_1,T_2,T_3,T_4$. Then $\aaa_3\aaa_4\cdots=\aaa^2\ddd^2$ determines $T_5,T_6$. Then $\bbb_1\ccc_3\cdots=\bbb_1\ccc_3\ccc_7$ and $\ccc_1\cdots=\bbb_7\ccc_1\cdots$ determine $T_7$; $\aaa_5\ddd_3\ddd_7\cdots=\aaa^2\ddd^2$ determines $T_{8}$. Similarly, $T_9,T_{10},\dots,T_{16}$ are determined.

	\begin{figure}[htp]
		\centering
		\begin{tikzpicture}[>=latex]

			\foreach \b in {0,1,2,3,4,5}
			{
				\begin{scope}[rotate=60*\b ] 
					\draw (0,1)--(0,2)--(0,3)--(-2.598,1.5)--(-1.732,1)--(-0.866,0.5)--(0,1);
					
					\draw[line width=1.5] 
					(0,2)--(-1.732,1);

				\end{scope}
			}
			\draw[line width=1.5] 
			(0,-1)--(0,1)
			(-2.598,1.5)--(-3.464,2)
			(2.598,-1.5)--(3.464,-2);

			\node at (-0.25,0.6) {\small $\aaa$};
			\node at (-0.25,-0.6) {\small $\ddd$};
			\node at (0.25,0.6) {\small $\ddd$};
			\node at (0.25,-0.6) {\small $\aaa$};
			
			\node at (-0.65,0.3) {\small $\bbb$};
			\node at (0.65,0.3) {\small $\ccc$};
			\node at (-0.65,-0.3) {\small $\ccc$};
			\node at (0.65,-0.3) {\small $\bbb$};
			
			\node at (-0.9,0.7) {\small $\ccc$};
			\node at (0.9,0.7) {\small $\ccc$};
			\node at (-0.9,-0.8) {\small $\ccc$};
			\node at (0.9,-0.8) {\small $\ccc$};
			
			\node at (-0.25,-1.2) {\small $\bbb$};
			\node at (0.25,-1.2) {\small $\bbb$};
			\node at (-0.25,1.2) {\small $\bbb$};
			\node at (0.25,1.2) {\small $\bbb$};
			
			\node at (-1.1,0.3) {\small $\ccc$};
			\node at (1.1,0.3) {\small $\bbb$};
			\node at (-1.1,-0.3) {\small $\bbb$};
			\node at (1.1,-0.3) {\small $\ccc$};
			
			\node at (-1.3,1) {\small $\ddd$};
			\node at (1.3,1) {\small $\ddd$};
			\node at (-1.3,-1) {\small $\ddd$};
			\node at (1.3,-1) {\small $\ddd$};
			
			\node at (-1.5,-0.7) {\small $\aaa$};
			\node at (-1.5,0.6) {\small $\ddd$};
			\node at (1.5,-0.7) {\small $\ddd$};
			\node at (1.5,0.7) {\small $\aaa$};
			
			\node at (-0.25,-1.6) {\small $\aaa$};
			\node at (-0.25,1.6) {\small $\aaa$};
			\node at (0.25,-1.6) {\small $\aaa$};
			\node at (0.25,1.6) {\small $\aaa$};
			
			\node at (-0.25,-2.1) {\small $\ddd$};
			\node at (-0.25,2.1) {\small $\ddd$};
			\node at (0.25,-2.1) {\small $\ddd$};
			\node at (0.25,2.1) {\small $\ddd$};
			
			\node at (-0.25,-2.7) {\small $\ccc$};
			\node at (-0.25,2.6) {\small $\ccc$};
			\node at (0.25,2.6) {\small $\ccc$};
			\node at (0.25,-2.7) {\small $\ccc$};
			
			\node at (-2.1,-1.5) {\small $\bbb$};
			\node at (2.1,1.5) {\small $\bbb$};
			\node at (2.1,-1.5) {\small $\bbb$};
			\node at (-2.1,1.5) {\small $\bbb$};
			
			\node at (-2.45,-1.2) {\small $\ccc$};
			\node at (2.45,1.2) {\small $\ccc$};
			\node at (2.45,-1.2) {\small $\bbb$};
			\node at (-2.45,1.1) {\small $\bbb$};
			
			\node at (0,-3.3) {\small $\bbb$};
			\node at (0,3.3) {\small $\bbb$};
			
			\node at (-2.6,1.8) {\small $\aaa$};
			\node at (2.6,-1.8) {\small $\aaa$};
			\node at (-2.8,1.3) {\small $\ddd$};
			\node at (2.8,-1.3) {\small $\ddd$};
			
			\node at (-2.7,-1.8) {\small $\ccc$};
			\node at (2.7,1.8) {\small $\ccc$};
			
			\node at (-1.7,1.2) {\small $\aaa$};
			\node at (-1.7,-1.2) {\small $\aaa$};
			\node at (1.7,1.2) {\small $\aaa$};
			\node at (1.7,-1.2) {\small $\aaa$};
			
			\node at (-1.9,0.9) {\small $\aaa$};
			\node at (1.9,-0.9) {\small $\aaa$};
			\node at (-1.9,-0.9) {\small $\ddd$};
			\node at (1.9,0.9) {\small $\ddd$};

			\node[draw,shape=circle, inner sep=0.5] at (-0.7,1.2) {\small $3$};
			\node[draw,shape=circle, inner sep=0.5] at (0.7,1.2) {\small $4$};
			\node[draw,shape=circle, inner sep=0.5] at (-1.2,1.8) {\small $5$};
			\node[draw,shape=circle, inner sep=0.5] at (1.2,1.8) {\small $6$};
			\node[draw,shape=circle, inner sep=0.5] at (2.2,0) {\small $14$};
			\node[draw,shape=circle, inner sep=0.5] at (1.4,0) {\small $13$};
			\node[draw,shape=circle, inner sep=0.5] at (-1.4,0) {\small $7$};
			\node[draw,shape=circle, inner sep=0.5] at (-2.2,0) {\small $8$};
			\node[draw,shape=circle, inner sep=0.5] at (-0.4,0) {\small $1$};
			\node[draw,shape=circle, inner sep=0.5] at (0.4,0) {\small $2$};
			\node[draw,shape=circle, inner sep=0.5] at (-0.7,-1.2) {\small $9$};
			\node[draw,shape=circle, inner sep=0.5] at (-1.2,-1.8) {\small $10$};
			\node[draw,shape=circle, inner sep=0.5] at (0.7,-1.2) {\small $11$};
			\node[draw,shape=circle, inner sep=0.5] at (1.2,-1.8) {\small $12$};
			\node[draw,shape=circle, inner sep=0.5] at (1.8,3.1) {\small $15$};
			\node[draw,shape=circle, inner sep=0.5] at (-1.8,-3.1) {\small $16$};

		\end{tikzpicture}
		\caption{$f=16$, $T(8\bbb\ccc^2,6\aaa^2\ddd^2,4\aaa\bbb^2\ddd)$.} 
		\label{ab2,c2d2}
	\end{figure}
	By Lemma \ref{calculate-1},  we get 
	\begin{equation*}
		\aaa=1-\arcsin{\tfrac{\sqrt{10+4\sqrt 2}}{\sqrt{17}}}\approx 0.5906,\,\,\,\,
		a=\tfrac14,  \,\,\,\,
		b=\arccos{\tfrac{2\sqrt 2-1}{4}}\approx 0.3488.
	\end{equation*}

	\subsubsection*{Case \{$\bbb\ccc^2,\aaa\ddd^3$\}}

	Denote the vertex $\aaa^2\cdots$ by $\bn$, $n_1\ge2$.
	
	If $n_4=0$, by $4n_2=3n_1+2n_3\ge6$, then $n_2\ge2$. By $2\aaa+2\bbb+\ccc>2$, we get $n_3=0$. Therefore, $3n_1=4n_2$. Then $n_1\ge4$ and $n_2\ge3$, contradicting $4\aaa+3\bbb>2$.  
	
	If $n_4=1$, by Parity Lemma, we get $n_1\ge3$. By $3n_1+2n_3=4n_2+1$, we get $n_2\ge2$, contradicting $3\aaa+2\bbb+\ddd>2$.
	
	If $n_4\ge2$, by $4n_2+n_4=3n_1+2n_3\ge6$, we get either $n_2\ge1$ or $n_4\ge6$,  contradicting $2\aaa+\bbb+2\ddd>2$ and $2\aaa+6\ddd>2$.
	
	
	\subsubsection*{Case \{$\bbb\ccc^2,\ddd^4$\}}
	
	Denote the vertex $\aaa^2\cdots$ by $\bn$, $n_1\ge2$. By $n_1+n_3=2n_2$, we have $n_2\ge1$. 
	
	If $n_2=1$, by $n_1+n_3=2$ and $n_1\ge2$, then $\bn=(2\,1\,0\,n_4)$. By Parity Lemma, we get $n_4\ge2$, contradicting $2\aaa+\bbb+2\ddd>2$.
	
	If $n_2\ge2$, by $n_1+n_3=2n_2$, we get $n_1\ge4$, or $n_1\ge3$ and $n_3\ge1$, or $n_1\ge2$ and $n_3\ge2$, contradicting $4\aaa+2\bbb>2$, $3\aaa+2\bbb+\ccc>2$ and $2\aaa+2\bbb+2\ccc>2$.	
\end{proof}

\begin{remark}
	
	Combined with the classification of rational $a^3b$-tilings in \cite{lw}, the only other tiling with $\bbb\ccc^2$ as the unique degree $3$ vertex type is	
	$T(18\bbb\ccc^2$, $6\aaa^3\ddd,6\aaa^2\bbb^2,6\aaa\bbb\ddd^3,2\ddd^6)$ with four angles being $(5,4,7,3)/9$ and $f=36$.
\end{remark}

\begin{proposition}\label{special9}
	In an $a^3b$-tiling with some irrational angle, if all degree $3$ vertices are  $\bbb^3$, then the tiling is a special quadrilateral subdivision of the octahedron  $T(8\bbb^3,12\aaa^2\ccc^2,6\ddd^4)$ with $24$ tiles. 
\end{proposition}

\begin{proof}
	 Lemma \ref{tildev21} implies that one of $\{\aaa^4,\aaa^3\ddd,\aaa^2\ddd^2,\aaa\ddd^3,\ddd^4\}$ must appear.  Applying the Irrational Angle Lemma to each case, we get Table \ref{Tab-b3}. 
	
	\begin{table*}[htp]                        
		\centering     
		
		~\\ 
		\resizebox{\textwidth}{14mm}{\begin{tabular}{c|c|c|c}	 
				
				Cases & Angles &Irrational angle Lemma & Conclusion \\
				\hline 
				$\bbb^3,\aaa^4$&$(\frac12,\frac23,\ccc,\frac56-\ccc+\frac4f)$, $0<\ccc<\frac56+\frac4f$
				&$n_3=n_4$& No $\aaa\ccc\cdots$, a contradiction\\
				\hline 
				$\bbb^3,\aaa^3\ddd$&$(\frac13+\frac\ccc2-\frac2f,\frac23,\ccc,1-\frac{3\ccc}{2}+\frac6f)$, $0<\ccc<\frac23+\frac4f$
				&$n_1+2n_3=3n_4$ &No $\ddd^2\cdots$, a contradiction\\
				\hline 
				$\bbb^3,\aaa^2\ddd^2$&$(\aaa,\frac23,\frac13+\frac4f,1-\aaa)$, $0<\aaa<1$
				&$n_1=n_4$& No $\aaa\ccc\cdots$, a contradiction\\
				\hline 
				$\bbb^3,\aaa\ddd^3$&$(1-\frac{3\ccc}{2}+\frac6f,\frac23,\ccc,\frac13+\frac\ccc2-\frac2f)$, $0<\ccc<\frac23+\frac4f$
				& $3n_1=2n_3+n_4$&No $\aaa^2\cdots$, a contradiction\\
				\hline 
				$\bbb^3,\ddd^4$&$(\frac56-\ccc+\frac4f,\frac23,\ccc,\frac12)$, $0<\ccc<\frac56+\frac4f$
				&$n_1=n_3$& $\text{AVC}\sub\{\bbb^3,\aaa^2\ccc^2,\ddd^4\}$.\\
				\hline 
		\end{tabular}}
	\caption{All degree $3$ vertices being $\bbb^3$.}\label{Tab-b3}        
	\end{table*}

	\subsubsection*{Case \{$\bbb^3,\aaa^4$\}} 
	
    By Lemma \ref{aadlemma}, the AAD of $\bbb^3$ gives $\aaa\ccc\cdots$, denoted by $\bn$ with $n_1,n_3\ge1$.   
	
	If $n_3=1$, then $n_4=n_3=1$. By Lemma \ref{anglesum} and Parity Lemma, we get $n_2=0$ and $n_1\ge3$, contradicting $3\aaa+\ccc+\ddd>2$.
	
	If $n_3\ge2$, then $n_4=n_3\ge2$, contradicting $\aaa+2\ccc+2\ddd>2$.			
	
	
	\subsubsection*{Case \{$\bbb^3,\aaa^3\ddd$\}}	
	Denote the vertex $\ddd^2\cdots$ by $\bn$, $n_4\ge2$. 
	
	If $n_1=0$, by $2n_3=3n_4\ge6$, then $n_3\ge3$, contradicting $3\ccc+2\ddd>2$. 
	
	If $n_1=1$, by $1+2n_3=3n_4$ and Parity Lemma, we get $n_4\ge3$ and $n_3\ge4$, contradicting $\aaa+4\ccc+3\ddd>2$.
	
	If $n_1\ge2$, by $\aaa^3\ddd$ and $n_1+2n_3=3n_4\ge6$, then $n_1=2$ and $n_3\ge2$, contradicting  $2\aaa+2\ccc+2\ddd>2$.	
	
	
	\subsubsection*{Case \{$\bbb^3,\aaa^2\ddd^2$\}}		
	The AAD of $\bbb^3$ gives $\aaa\ccc\cdots$, denoted by $\bn$ with $n_1,n_3\ge1$. By $n_4=n_1\ge1$ and $R(\aaa\ccc\ddd\cdots)<\aaa+\ddd,\bbb,2\ccc$, we deduce that $\bn=(1\,0\,2\,1)$.  Then $\ccc=\frac12$, and its AAD gives $\bbb\thin\ccc\cdots$ or $\bbb\thin\ddd\cdots$. By $\aaa+\ddd=1$, $\bbb=\frac23$, $\ccc=\frac12$, $n_1=n_4$ and Lemma \ref{anglesum}, there is no $\bbb\ccc\cdots$ and $\bbb\ddd\cdots$.	
	
	\subsubsection*{Case \{$\bbb^3,\aaa\ddd^3$\}}	
	Denote the vertex $\aaa^2\cdots$ by $\bn$, $n_1\ge2$. 
	
	If $n_4=0$, by $2n_3=3n_1\ge6$, then $n_3\ge3$, contradicting $2\aaa+3\ccc>2$. 
	
	If $n_4=1$, by $2n_3+1=3n_1\ge6$ and Parity Lemma, we get $n_1\ge3$ and $n_3\ge4$, contradicting $3\aaa+4\ccc+\ddd>2$.
	
	If $n_4\ge2$, by $\aaa\ddd^3$ and $2n_3+n_4=3n_1\ge6$, then $n_4=2$ and $n_3\ge2$, contradicting  $2\aaa+2\ccc+2\ddd>2$.	
	
				
	\subsubsection*{Case \{$\bbb^3,\ddd^4$\}}	
	We will find all possible vertices $\bn$.	
	If $n_1=0$, then $n_3=n_1=0$. By $\bbb=\frac23$ and $\ddd=\frac12$, we get $\bn=(0\,3\,0\,0)$ or $(0\,0\,0\,4)$.
		
	If $n_1=1$, then $n_3=n_1=1$. By Parity Lemma and Lemma \ref{anglesum}, we get $n_2=0$ and $n_4\ge3$, contradicting $\aaa+\ccc+3\ddd>2$.
	
	If $n_1\ge2$, then $n_3=n_1\ge2$ and $R(\aaa^2\ccc^2\cdots)<\aaa+\ccc,\bbb,\ddd$, we deduce that $\bn=(2\,0\,2\,0)$.
	
	Therefore, $\text{AVC}\sub\{\bbb^3,\aaa^2\ccc^2,\ddd^4\}$.
	By Balance Lemma, $\aaa^2\ccc^2$ is a vertex. Then by Lemma \ref{anglesum}, we get $f=24$. There is only one solution satisfying Balance Lemma:  $\{8\bbb^3,12\aaa^2\ccc^2,6\ddd^4\}$.  
	
	In the first picture of Fig. \ref{case a3 d4}, the unique AAD $\ddd^4=\thin^\ccc\ddd^\aaa\thick^\aaa\ddd^\ccc\thin^\ccc\ddd^\aaa\thick^\aaa\ddd^\ccc\thin$ determines $T_1,T_2,T_3,T_4$. Then $\aaa_1\aaa_4\cdots=\aaa_1\aaa_4\ccc_5\ccc_6$ and $\bbb_4\cdots=\bbb_4\bbb_5\bbb_{12}$ determine $T_5$; $\bbb_1\cdots=\bbb_1\bbb_6\bbb_7$ determines $T_6$; $\ccc_1\ccc_2\cdots=\aaa_7\aaa_8\ccc_1\ccc_2$ determines $T_7,T_8$; $\bbb_2\bbb_8\cdots=\bbb_2\bbb_8\bbb_9$ and $\aaa_2\aaa_3\cdots=\aaa_2\aaa_3\ccc_9\ccc_{10}$ determine $T_9$; $\bbb_3\cdots=\bbb_3\bbb_{10}\bbb_{11}$ determines $T_{10}$; $\ccc_3\ccc_4\cdots=\aaa_{11}\aaa_{12}\ccc_3\ccc_4$ determines $T_{11},T_{12}$. Similarly, we can determine $T_{13},T_{14},\dots,T_{24}$. 
	
	\begin{figure}[htp]
		\centering
		\begin{minipage}[c]{0.5\linewidth}
			\begin{tikzpicture}[>=latex,scale=2.1]
				\foreach \a in {0,1,2,3}
				\draw[rotate=90*\a]
				(0.4,0.4)--(0,0.55)--(0,0.9)
				(0,0.55)--(-0.4,0.4)
				(-1.3,1.3)--(0,1.3)--(1.3,1.3)
				(0,0.9)--(0,1.3);
				\foreach \a in {0,2}
				\draw[line width=1.5, rotate=90*\a]   
				(0,0)--(0.4,0.4)
				(-0.4,0.4)--(-0.9,0.9)
				(0,0.9)--(0.9,0.9)
				(0.9,0)--(0.9,0.9)
				(-0.9,0.9)--(-1.3,1.3)
				(1.3,1.3)--(1.6,1.6);
				\foreach \a in {0,2}
				\draw[rotate=90*\a]   
				(0,0)--(-0.4,0.4)
				(0.4,0.4)--(0.9,0.9)
				(-0.9,0)--(-0.9,0.9)--(0,0.9)
				(0.9,0.9)--(1.3,1.3)
				(-1.3,1.3)--(-1.6,1.6);
				
				\node[draw,shape=circle, inner sep=0.5] at (0,0.29) {\small $1$};%
				\node[draw,shape=circle, inner sep=0.5] at (0.3,0.7) {\small $6$};%
				\node[draw,shape=circle, inner sep=0.5] at (-0.3,0.7) {\small $7$};%
				\node[draw,shape=circle, inner sep=0.5] at (0.3,0) {\small $4$};%
				\node[draw,shape=circle, inner sep=0.5] at (0,-0.28) {\small $3$};%
				\node[draw,shape=circle, inner sep=0.5] at (-0.3,0) {\small $2$};%
				\node[draw,shape=circle, inner sep=0.5] at (-0.7,0.3) {\small $8$};%
				\node[draw,shape=circle, inner sep=0.5] at (-1.1,0.5) {\footnotesize $17$};%
				\node[draw,shape=circle, inner sep=0.5] at (-0.5,1.1) {\footnotesize $16$};%
				\node[draw,shape=circle, inner sep=0.5] at (0.7,0.3) {\small $5$};%
				\node[draw,shape=circle, inner sep=0.5] at (1.1,0.5) {\footnotesize $14$};%
				\node[draw,shape=circle, inner sep=0.5] at (0.5,1.1) {\footnotesize $15$};%
				
				\node[draw,shape=circle, inner sep=0.5] at (0,1.6) {\footnotesize $23$};
				\node[draw,shape=circle, inner sep=0.5] at (1.6,0) {\footnotesize $22$};
				\node[draw,shape=circle, inner sep=0.5] at (0.7,-0.3) {\footnotesize $12$};%
				\node[draw,shape=circle, inner sep=0.5] at (1.1,-0.5) {\footnotesize $13$};%
				\node[draw,shape=circle, inner sep=0.5] at (0.3,-0.7) {\footnotesize $11$};%
				\node[draw,shape=circle, inner sep=0.5] at (0.5,-1.1) {\footnotesize $20$};%
				\node[draw,shape=circle, inner sep=0.5] at (0,-1.6) {\footnotesize $21$};%
				
				\node[draw,shape=circle, inner sep=0.5] at (-1.6,0) {\footnotesize $24$};
				\node[draw,shape=circle, inner sep=0.5] at (-0.7,-0.3) {\small $9$};%
				\node[draw,shape=circle, inner sep=0.5] at (-1.1,-0.5) {\footnotesize $18$};%
				\node[draw,shape=circle, inner sep=0.5] at (-0.3,-0.7) {\footnotesize $10$};%
				\node[draw,shape=circle, inner sep=0.5] at (-0.5,-1.1) {\footnotesize $19$};%
				
				\node at (0,0.12){\small $\ddd$}; \node at (0.12,0){\small $\ddd$};
				\node at (0,-0.12){\small $\ddd$}; \node at (-0.12,0){\small $\ddd$};   
				\node at (0,0.44){\small $\bbb$}; \node at (0.46,0){\small $\bbb$};
				\node at (0,-0.46){\small $\bbb$}; \node at (-0.46,0){\small $\bbb$};
				\node at (0.25,0.36){\small $\aaa$}; \node at (0.37,0.24){\small $\aaa$};
				\node at (-0.25,-0.37){\small $\aaa$}; \node at (-0.37,-0.24){\small $\aaa$};   
				\node at (-0.22,0.36){\small $\ccc$}; \node at (-0.38,0.24){\small $\ccc$};
				\node at (0.22,-0.36){\small $\ccc$}; \node at (0.38,-0.24){\small $\ccc$};
				\node at (0.08,0.62){\small $\bbb$}; \node at (-0.08,0.6){\small $\bbb$};
				\node at (0.08,-0.62){\small $\bbb$}; \node at (-0.08,-0.62){\small $\bbb$};
				\node at (0.6,0.08){\small $\bbb$}; \node at (0.6,-0.1){\small $\bbb$};
				\node at (-0.6,0.08){\small $\bbb$}; \node at (-0.6,-0.1){\small $\bbb$};   
				\node at (0.08,0.8){\small $\aaa$}; \node at (-0.08,0.8){\small $\ccc$};
				\node at (0.08,-0.8){\small $\ccc$}; \node at (-0.08,-0.8){\small $\aaa$};
				\node at (0.8,0.08){\small $\aaa$}; \node at (0.8,-0.1){\small $\ccc$};
				\node at (-0.8,0.1){\small $\ccc$}; \node at (-0.8,-0.08){\small $\aaa$};
				\node at (0.38,0.53){\small $\ccc$}; \node at (0.51,0.37){\small $\ccc$};
				\node at (-0.4,-0.52){\small $\ccc$}; \node at (-0.51,-0.37){\small $\ccc$};
				\node at (0.38,-0.53){\small $\aaa$}; \node at (0.51,-0.37){\small $\aaa$};
				\node at (-0.38,0.51){\small $\aaa$}; \node at (-0.51,0.37){\small $\aaa$};
				\node at (0.7,0.8){\small $\ddd$}; \node at (0.8,0.7){\small $\ddd$};
				\node at (-0.7,-0.8){\small $\ddd$}; \node at (-0.8,-0.7){\small $\ddd$};
				\node at (-0.69,0.81){\small $\ddd$}; \node at (-0.825,0.69){\small $\ddd$};
				\node at (0.69,-0.825){\small $\ddd$}; \node at (0.825,-0.7){\small $\ddd$};   
				\node at (0.9,1){\small $\ddd$}; \node at (1,0.9){\small $\ddd$};
				\node at (-0.9,-1){\small $\ddd$}; \node at (-1,-0.9){\small $\ddd$};
				\node at (-0.85,1){\small $\ddd$}; \node at (-0.98,0.85){\small $\ddd$};
				\node at (0.85,-1){\small $\ddd$}; \node at (1,-0.85){\small $\ddd$};
				\node at (0.08,1){\small $\aaa$}; \node at (-0.08,1){\small $\ccc$};
				\node at (0.08,-1){\small $\ccc$}; \node at (-0.08,-1){\small $\aaa$};
				\node at (1,0.08){\small $\aaa$}; \node at (1,-0.1){\small $\ccc$};
				\node at (-1,0.1){\small $\ccc$}; \node at (-1,-0.08){\small $\aaa$};
				\node at (0.08,1.2){\small $\bbb$}; \node at (-0.08,1.2){\small $\bbb$};
				\node at (0.08,-1.2){\small $\bbb$}; \node at (-0.08,-1.2){\small $\bbb$};
				\node at (1.2,0.08){\small $\bbb$}; \node at (1.2,-0.1){\small $\bbb$};
				\node at (-1.2,0.08){\small $\bbb$}; \node at (-1.2,-0.1){\small $\bbb$};   
				\node at (0,1.38){\small $\bbb$}; \node at (0,-1.4){\small $\bbb$}; \node at (1.38,0){\small $\bbb$};  \node at (-1.4,0){\small $\bbb$};
				\node at (0,1.85){\small $\ddd$}; \node at (0,-1.85){\small $\ddd$}; \node at (1.85,0){\small $\ddd$};  \node at (-1.85,0){\small $\ddd$};
				\node at (1.1,1.2){\small $\ccc$}; \node at (1.22,1.1){\small $\ccc$};
				\node at (-1.1,-1.2){\small $\ccc$}; \node at (-1.22,-1.1){\small $\ccc$};
				\node at (1.08,-1.23){\small $\aaa$}; \node at (1.22,-1.1){\small $\aaa$};
				\node at (-1.1,1.2){\small $\aaa$}; \node at (-1.23,1.08){\small $\aaa$};
				\node at (1.25,1.4){\small $\aaa$}; \node at (1.4,1.25){\small $\aaa$};
				\node at (-1.25,-1.4){\small $\aaa$}; \node at (-1.4,-1.25){\small $\aaa$};
				\node at (-1.25,1.4){\small $\ccc$}; \node at (-1.4,1.25){\small $\ccc$};
				\node at (1.25,-1.4){\small $\ccc$}; \node at (1.4,-1.25){\small $\ccc$};
			\end{tikzpicture}
		\end{minipage}%
		\begin{minipage}[c]{0.5\linewidth}\hspace*{4em}
			\includegraphics[scale=0.22]{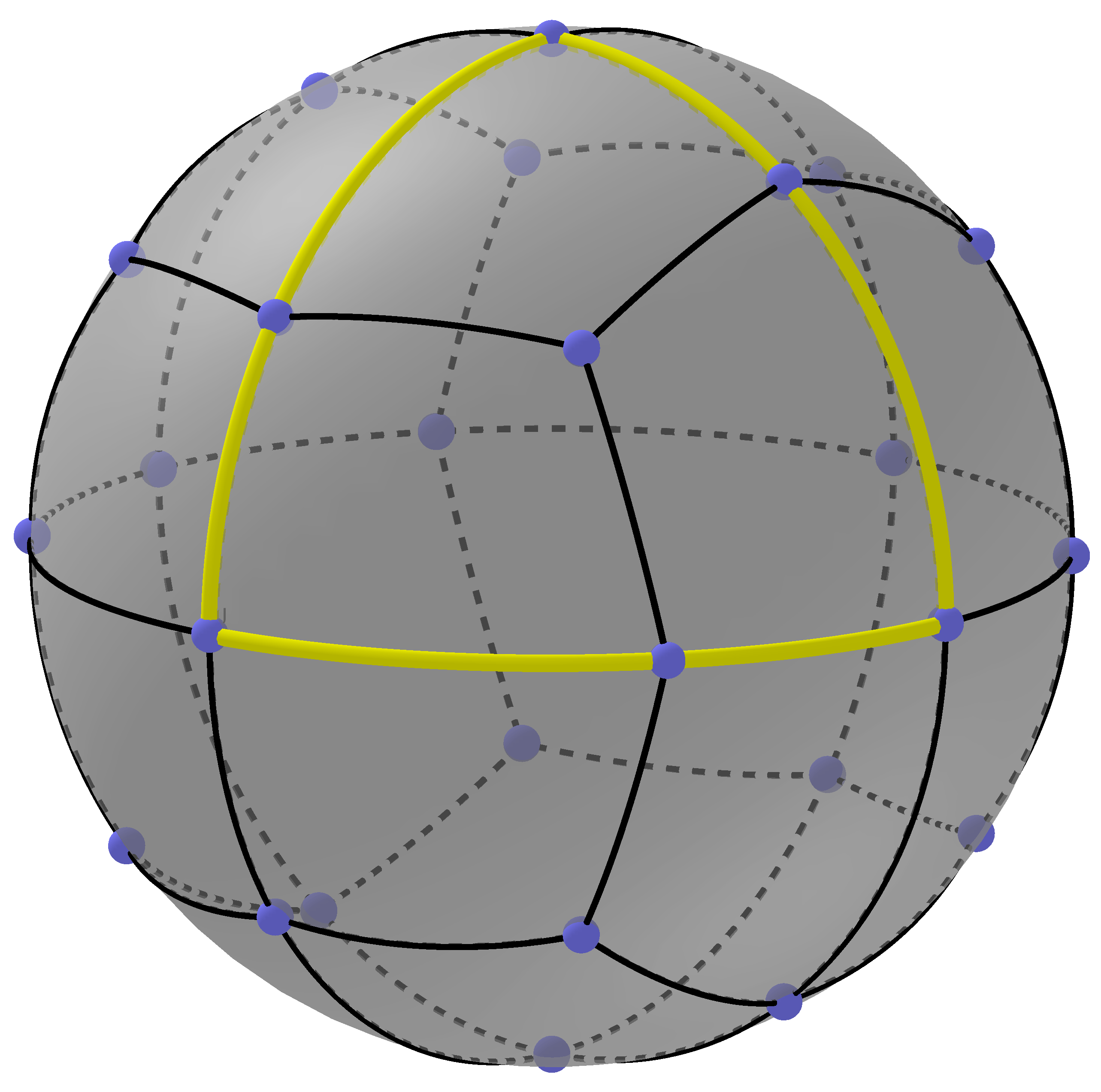}
		\end{minipage}
		\caption{A special quadrilateral subdivision of the octahedron.} \label{case a3 d4}
	\end{figure}
	By Lemma \ref{calculate-1}, we get
	\[
		\aaa=\arcsin{\tfrac{\sqrt{4+\sqrt3}}{\sqrt6}}\approx 0.4322, 
		a=\arcsin{\tfrac{\sqrt2}{\sqrt{4+\sqrt3}}}\approx 0.2011,  
		b=\tfrac12-a\approx 0.2988
	\]
    It is a special quadrilateral subdivision of the octahedron, whose $3$D picture is shown in the second of Fig. \ref{case a3 d4}. 
\end{proof}

\begin{remark}
	Combined with the classification of rational $a^3b$-tilings in \cite{lw}, Proposition \ref{special9} still holds for $a^3b$-quadrilaterals with arbitrary angles.
\end{remark}

\begin{proposition}\label{special5}
	In an $a^3b$-tiling with some irrational angle, if all degree $3$ vertices are $\aaa^2\bbb$, then the tiling is  $T(8\aaa^2\bbb,8\bbb\ccc\ddd^2,2\ccc^4)$ with $16$ tiles.
\end{proposition}

\begin{proof}
	By Lemma \ref{a2b4}, one of $\{\bbb^4,\bbb^3\ccc,\bbb^2\ccc^2,\bbb^2\ddd^2,\bbb\ccc^3,\bbb\ccc\ddd^2,\ccc^4,\ccc^2\ddd^2,\ddd^4\}$ must appear. By Irrational Angle Lemma for each case, we get Table \ref{Tab-a2b}.

	\begin{table*}[htp]                        
		\centering     
	    
		~\\ 
		\resizebox{\textwidth}{21mm}{\begin{tabular}{c|c|c|c}	 
				
				Cases &Angles &Irrational angle Lemma & Conclusion \\
				\hline 
				$\aaa^2\bbb,\bbb^4$&$(\frac34,\frac12,\ccc,\frac34-\ccc+\frac4f)$, $0<\ccc<\frac34+\frac4f$&$n_3=n_4$
				&No $\bbb\ccc\cdots$, a contradiction\\
				\hline 
				$\aaa^2\bbb,\bbb^3\ccc$&$(1-\frac\bbb2,\bbb,2-3\bbb,-1+\frac{5\bbb}{2}+\frac4f)$, $\frac25-\frac{8}{5f}<\bbb<\frac23$
				& $n_1+6n_3=2n_2+5n_4$& No $\ddd^2\cdots$, a contradiction \\
				\hline  
				$\aaa^2\bbb,\bbb^2\ccc^2$&$(1-\frac\bbb2,\bbb,1-\bbb,\frac\bbb2+\frac4f)$, $0<\bbb<1$
				&$n_1+2n_3=2n_2+n_4$&No $\ddd^2\cdots$, a contradiction\\
				\hline
				$\aaa^2\bbb,\bbb^2\ddd^2$&$(1-\frac\bbb2,\bbb,\frac\bbb2+\frac4f,1-\bbb)$, $0<\bbb<1$&
				$n_1+2n_4=2n_2+n_3$&No $\bbb\ccc\cdots$, a contradiction\\
				\hline 
				$\aaa^2\bbb,\bbb\ccc^3$&$(1-\frac\bbb2,\bbb,\frac23-\frac\bbb3,\frac13-\frac\bbb6+\frac4f)$, $0<\bbb<2$
				&$3n_1+2n_3+n_4=6n_2$&No $\ddd^2\cdots$, a contradiction\\
				\hline
				$\aaa^2\bbb,\bbb\ccc\ddd^2$&$(1-\frac\bbb2,\bbb,\frac8f,1-\frac\bbb2-\frac4f)$, $0<\bbb<2-\frac8f$&
				$n_1+n_4=2n_2$&$\text{AVC}=\{\aaa^2\bbb,\bbb\ccc\ddd^2,\ccc^{\frac f4}\}$\\
				\hline  
				$\aaa^2\bbb,\ccc^4$&$(1-\frac\bbb2,\bbb,\frac12,\frac12-\frac\bbb2+\frac4f)$, $0<\bbb<1+\frac8f$
				&$n_1+n_4=2n_2$&$\text{AVC}\sub\{\aaa^2\bbb,\bbb\ccc\ddd^2,\ccc^4\}$\\
				\hline				
				$\aaa^2\bbb,\ccc^2\ddd^2$&$(1-\frac4f,\frac8f,\ccc,1-\ccc)$, $0<\ccc<1$&$n_3=n_4$&No $\bbb\ccc\cdots$, a contradiction \\				
				\hline 
				$\aaa^2\bbb,\ddd^4$&$(1-\frac\bbb2,\bbb,\frac12-\frac\bbb2+\frac4f,\frac12)$, $0<\bbb<1+\frac8f$
				&$n_1+n_3=2n_2$&No $\bbb\ccc\cdots$, a contradiction\\
				\hline 
		\end{tabular}}
		\caption{All degree $3$ vertices being $\aaa^2\bbb$.}\label{Tab-a2b}    
	\end{table*}

	\subsubsection*{Case \{$\aaa^2\bbb,\bbb^4$\}}	
	The AAD of $\aaa^2\bbb$ gives $\bbb\ccc\cdots$, denoted by $\bn$ with $n_2,n_3\ge1$. By $n_4=n_3\ge1$, Lemma \ref{anglesum} and Parity Lemma, we get $n_1=0$ and $n_3=n_4\ge2$, contradicting $\bbb+2\ccc+2\ddd>2$.
	
	\subsubsection*{Case \{$\aaa^2\bbb,\bbb^3\ccc$\}}
	The AAD of $\aaa^2\bbb$ gives $\ddd^2\cdots$, denoted by $\bn$, $n_4\ge2$. 
	
	If $n_1=0$, by $6n_3=2n_2+5n_4\ge10$, then $n_3\ge2$. By $\bbb+2\ccc+2\ddd>2$, we get $n_2=0$. Therefore, $6n_3=5n_4$. Then $n_3\ge5$ and $n_4\ge6$, contradicting $5\ccc+6\ddd>2$.
	
	If $n_1=1$, by Parity Lemma, we get $n_4\ge3$. By $1+6n_3=2n_2+5n_4$, we have $n_3\ge3$, contradicting $\aaa+3\ccc+3\ddd=4-2\bbb+\frac{12}{f}>2$.	
	
	If $n_1\ge2$, by $n_1+6n_3=2n_2+5n_4\ge10$, we get either $n_3\ge1$ or $n_1\ge10$, contradicting $2\aaa+\ccc+2\ddd>2$ and $10\aaa+2\ddd>2$.
	

	\subsubsection*{Case \{$\aaa^2\bbb,\bbb^2\ccc^2$\}}	
	The AAD of $\aaa^2\bbb$ gives $\ddd^2\cdots$, denoted by $\bn$, $n_4\ge2$. 
	
	If $n_1=0$, by $2n_3=2n_2+n_4\ge2$ and $\bbb^2\ccc^2$, we deduce that $n_3\ge2$ and $n_2=0$ or $1$. When $n_2=0$, we get $n_4\ge4$, contradicting $2\ccc+4\ddd>2$. When $n_2=1$, we get $n_3\ge2$, contradicting  $\bbb+2\ccc+2\ddd>2$. 
	
	If $n_1\ge1$, by $2\aaa+2\ddd>2$, then $n_1=1$. By Parity Lemma, we get $n_4\ge3$. By $1+2n_3=2n_2+n_4$, we have $n_3\ge1$, contradicting $\aaa+\ccc+3\ddd>2$. 
	
	
	\subsubsection*{Case \{$\aaa^2\bbb,\bbb^2\ddd^2$\}}		
	The AAD of $\aaa^2\bbb$ gives $\bbb\ccc\cdots$, denoted by $\bn$ with $n_2,n_3\ge1$. Then we get $n_1=0$ or $n_4=0$ by Lemma \ref{anglesum}.
	
	If $n_1=0$, by $2n_4=2n_2+n_3\ge3$, then $n_4\ge2$ and $n_3\ge2$, contradicting $\bbb+2\ccc+2\ddd>2$. 
	
	If $n_4=0$, by $n_1=2n_2+n_3\ge3$, then $n_1\ge3$, contradicting $3\aaa+\bbb+\ccc>2$.
	
	
	\subsubsection*{Case \{$\aaa^2\bbb,\bbb\ccc^3$\}}
	The AAD of $\aaa^2\bbb$ gives $\ddd^2\cdots$, denoted by $\bn$, $n_4\ge2$.
	By $3n_1+2n_3+n_4=6n_2$, we get $n_2\ge1$. 
	
	If $n_1=0$, by $2n_3+n_4=6n_2\ge6$, we get $n_3\ge2$, or $n_3=1$ and $n_4\ge4$, or $n_3=0$ and $n_4\ge6$, contradicting $\bbb+2\ccc+2\ddd>2$, $\bbb+\ccc+4\ddd>2$ and $\bbb+6\ddd>2$. 
	
	If $n_1\ge1$, by $\aaa^2\bbb$, then $n_1=1$. By $n_4\ge2$ and Parity Lemma, we get $n_4\ge3$, contradicting $\aaa+\bbb+3\ddd>2$.
	
	
	\subsubsection*{Case \{$\aaa^2\bbb,\bbb\ccc\ddd^2$\}}	
	We will find all possible vertices $\bn$. 	
	By $n_1+n_4=2n_2$, if $n_2=0$, then $n_1=n_4=0$ and $\bn=(0\,0\,\frac f4\,0)$ if $f\ge16$. 
	
	If $n_2=1$, by $n_1+n_4=2$, then $\bn=(2\,1\,n_3\,0)$, $(1\,1\,n_3\,1)$ or $(0\,1\,n_3\,2)$. By $\aaa^2\bbb$, $\bbb+2\ccc+2\ddd>2$ and Lemma \ref{anglesum}, we get $\bn=(2\,1\,0\,0)$ or $(0\,1\,1\,2)$.
	
	If $n_2\ge2$, by $\aaa^2\bbb$ and $n_1+n_4=2n_2$, we get either  $n_1=1,n_4\ge3$ or $n_1=0,n_4\ge4$. By $\aaa+2\bbb+3\ddd=4-\frac{12}{f}\ge\frac52$ and $2\bbb+4\ddd=4-\frac{16}{f}\ge2$, we deduce that $\bn=(0\,2\,0\,4)$ if $f=8$.
	
	Therefore, $\text{AVC}\sub\{\aaa^2\bbb,\bbb^2\ddd^4,\ccc^{\frac f4}\}$. By Balance Lemma, $\ccc^{\frac f4}$ is a vertex. Then $f\ge16$ and $\text{AVC}\sub\{\aaa^2\bbb,\bbb\ccc\ddd^2,\ccc^{\frac f4}\}$.
	There is only one solution satisfying Balance Lemma: $\{\frac{f}{2}\,\aaa^2\bbb,\frac{f}{2}\,\bbb\ccc\ddd^2,2\ccc^{\frac f4}\}$.
	
	We have $(\aaa,\bbb,\ccc,\ddd)=(\aaa,2-2\aaa,\frac2k,\aaa-\frac1k)$ and $f=4k(k\ge4)$.
	By $\bbb>0$, we have $\aaa<1$.  By $\aaa>\ddd$ and Lemma \ref{geometry1}, we get $\bbb<\ccc$, which implies that $\aaa>1-\frac1k\ge\frac34$ and the quadrilateral is convex. By Lemma \ref{relation}, we get $\aaa+\ddd<1+\bbb$. Therefore, $\aaa<\frac{3k+1}{4k}$. Then by $\frac{3k+1}{4k}>1-\frac1k$, we get $k<5$. Therefore, $k=4$ and $f=16$.	
	By Lemma \ref{calculate-1}, we get 
	\begin{align*}
		\aaa&=1-\arcsin{\tfrac{\sqrt{3\sqrt{10}-3\sqrt5-3\sqrt2+15}}{6}}\approx 0.7898,a=\arccos{\tfrac{\sqrt{10}+\sqrt{5}-\sqrt2-3}{2}}\approx 0.3362\\
		b&=\arccos{\tfrac{(27\sqrt5-43)\sqrt2+23\sqrt5-27}{(196\sqrt5-420)\sqrt2-267\sqrt5+623}}\approx 0.1052.
	\end{align*}

	In Fig.\,\ref{case1jie3}, the unique AAD $\ccc^4=\thin^\ddd\ccc^\bbb\thin^\ddd\ccc^\bbb\thin^\ddd\ccc^\bbb\thin\cdots$ determines $T_1,T_{1'},T_{1''}$. Then $\aaa_1\cdots=\aaa_1\aaa_2\bbb_3$ determines $T_2$. Then $\aaa_{1'}\cdots=\aaa_{1'}\aaa_{2'}\bbb_{3'}$ determines $T_{2'}$. Then $\bbb_1\ddd_{1'}\ddd_{2'}\cdots=\bbb_1\ccc_3\ddd_{1'}\ddd_{2'}$ determines $T_3$. Then $\bbb_{1'}\ddd_{1''}\cdots=\bbb_{1'}\ccc_{3'}\ddd_{1''}\cdots$ determines $T_{3'}$. Then $\aaa_3\bbb_2\cdots=\aaa_3\aaa_4\bbb_2$ determines $T_4$. Then $\ccc_{2'}\ddd_3\ddd_4\cdots=\bbb_{4'}\ccc_{2'}\ddd_3\ddd_4$. By $\bbb_4$, we get $\ccc_4\cdots=\ccc_4\ccc_{4'}\cdots$. This determines $T_{4'}$. Similarly, after repeating the process $2$ times, we get a unique earth map tiling $T(8\aaa^2\bbb,8\bbb\ccc\ddd^2,2\ccc^4)$.

	\begin{figure}[htp]
		\centering
		\begin{tikzpicture}[>=latex,scale=0.62]
			
			\foreach \b in {0,1,2,3}
			{
				\begin{scope}[xshift=2*\b cm]	
					\draw (0,-2)--(0,-6)
					(-2,0)--(-2,-4)
					(0,0)--(0,-2)--(-1,-2)--(-1,-4)--(-2,-4)--(-2,-6);
					
					\draw[line width=1.5]	
					(-2,-2)--(-1,-2)
					(-1,-4)--(0,-4);

					\node at (-1,-0.25) {\small $\ccc$};
					\node at (-1,-1.75) {\small $\aaa$};
					\node at (-0.25,-1.75) {\small $\bbb$};
					\node at (-1.75,-1.65) {\small $\ddd$};
					\node at (-0.75,-2.45) {\small $\bbb$};
					\node at (-1.25,-2.3) {\small $\aaa$};
					\node at (-0.25,-2.35) {\small $\ccc$};
					\node at (-1.75,-2.35) {\small $\ddd$};
					
					\node at (-0.75,-3.75) {\small $\aaa$};
					\node at (-1.25,-3.75) {\small $\bbb$};
					\node at (-0.25,-3.65) {\small $\ddd$};
					\node at (-1.75,-3.7) {\small $\ccc$};
					
					\node at (-0.25,-4.35) {\small $\ddd$};
					\node at (-1.75,-4.35) {\small $\bbb$};
					\node at (-1,-4.35) {\small $\aaa$};
					\node at (-1,-5.8) {\small $\ccc$};

				\end{scope}
			}

			\node[draw,shape=circle, inner sep=0.5] at (-1,-1) {\small $1$};
			\node[draw,shape=circle, inner sep=0.5] at (1,-1) {\small $1'$};
			\node[draw,shape=circle, inner sep=0.5] at (-0.5,-3) {\small $3$};
			\node[draw,shape=circle, inner sep=0.5] at (0.5,-3) {\small $2'$};
			\node[draw,shape=circle, inner sep=0.5] at (-1.5,-3) {\small $2$};
			\node[draw,shape=circle, inner sep=0.5] at (1.5,-3) {\small $3'$};
			\node[draw,shape=circle, inner sep=0.5] at (-1,-5) {\small $4$};
			\node[draw,shape=circle, inner sep=0.5] at (1,-5) {\small $4'$};
			
			\node[draw,shape=circle, inner sep=0.5] at (-1+4,-1) {\small $1''$};
			\node[draw,shape=circle, inner sep=0.5] at (-0.5+4,-3) {\small $3''$};
			\node[draw,shape=circle, inner sep=0.5] at (-1.5+4,-3) {\small $2''$};
			\node[draw,shape=circle, inner sep=0.5] at (-1+4,-5) {\small $4''$};

		\end{tikzpicture}\hspace{30pt}
		\caption{$f=16$, $T(8\aaa^2\bbb,8\bbb\ccc\ddd^2, 2\ccc^4)$.} 
		\label{case1jie3}
	\end{figure}	
	
	\subsubsection*{Case \{$\aaa^2\bbb,\ccc^4$\}}
	We will find all possible vertices $\bn$. By $n_1+n_4=2n_2$,	
	if $n_2=0$, then $n_1=n_4=0$. By $\ccc^4$, we get $\bn=(0\,0\,4\,0)$. 
	
	If $n_2=1$, by $n_1+n_4=2$, then $\bn=(2\,1\,n_3\,0)$, $(1\,1\,n_3\,1)$ or $(0\,1\,n_3\,2)$. By $\aaa^2\bbb$, $\bbb+2\ccc+2\ddd>2$ and Lemma \ref{anglesum}, we get $\bn=(2\,1\,0\,0)$ or $(0\,1\,1\,2)$.
	
	If $n_2\ge2$, by $\aaa^2\bbb$ and $n_1+n_4=2n_2$, we get either  $n_1=1,n_4\ge3$ or $n_1=0,n_4\ge4$, contradicting  $\aaa+2\bbb+3\ddd>2$ and $2\bbb+4\ddd>2$.

    Therefore, $\text{AVC}\sub\{\aaa^2\bbb,\bbb\ccc\ddd^2,\ccc^4\}$, which has been discussed in Case $\{\aaa^2\bbb,\bbb\ccc\ddd^2\}$ and gives a unique tiling in Fig.\,\ref{case1jie3}.	
    
    \subsubsection*{Case \{$\aaa^2\bbb,\ccc^2\ddd^2$\}}	   
    The AAD of $\aaa^2\bbb$ gives $\bbb\ccc\cdots$, denoted by $\bn$ with $n_2,n_3\ge1$. By $n_4=n_3\ge1$, Lemma \ref{anglesum} and Parity Lemma, we get $n_1=0$ and $n_3=n_4\ge2$, contradicting $\bbb+2\ccc+2\ddd>2$.
    
    \subsubsection*{Case \{$\aaa^2\bbb,\ddd^4$\}}	    
    The AAD of $\aaa^2\bbb$ gives $\bbb\ccc\cdots$, denoted by $\bn$ with $n_2,n_3\ge1$.
    
    If  $n_2=1$, by $n_1+n_3=2$ and Parity Lemma, we get either $n_1=0,n_3=2,n_4\ge2$ or $n_1=n_3=1,n_4\ge1$, contradicting $\bbb+2\ccc+2\ddd>2$ or Lemma \ref{anglesum}.
    
    If $n_2\ge2$, by $n_1+n_3=2n_2$, $\aaa^2\bbb$ and Parity Lemma, we get either $n_1=0,n_3\ge4$ or $n_1=1,n_3\ge3,n_4\ge1$, contradicting $2\bbb+4\ccc>2$ or Lemma \ref{anglesum}. 
\end{proof}

\begin{remark}
	Combined with the classification of rational $a^3b$-tilings in \cite{lw}, the only other tiling with $\aaa^2\bbb$ as the unique degree $3$ vertex type is
	$T(14\aaa^2\bbb$, $8\aaa\ddd^3,10\bbb\ccc^3,6\bbb^2\ccc\ddd^2)$ with four angles being $(15,6,10,7)/18$ and $f=36$.
\end{remark}

\begin{proposition}\label{special4}
	In an $a^3b$-tiling with some irrational angle, if all degree $3$ vertices are  $\bbb\ddd^2$, then the tiling is $T(8\bbb\ddd^2,8\aaa^2\ccc^2,2\bbb^{4})$ with $16$ tiles. 
\end{proposition}

\begin{proof}
	By Lemma $\ref{a2b4}'$, one of $\{\aaa^4,\aaa^2\bbb^2,\aaa^2\bbb\ccc,\aaa^2\ccc^2,\bbb^4,\bbb^3\ccc,\bbb^2\ccc^2,\bbb\ccc^3,\ccc^4\}$ must appear.  By Irrational Angle Lemma for each case, we get Table \ref{Tab-bd2}. By $\bbb\ddd^2$ and Lemma \ref{geometry4}, we always have $2\aaa+\bbb>1$.

	\begin{table*}[htp]                        
		\centering     
	   
		~\\ 
		\resizebox{\textwidth}{23mm}{\begin{tabular}{c|c|c|c}	 
				
				Cases & Angles &Irrational angle Lemma & Conclusion \\
				\hline 
				$\bbb\ddd^2,\aaa^4$&$(\frac12,\bbb,\frac12-\frac\bbb2+\frac4f,1-\frac\bbb2)$, $0<\bbb<1+\frac8f$ 
				&$2n_2=n_3+n_4$& No $\ccc^2\cdots$, a contradiction\\
				\hline
				$\bbb\ddd^2,\aaa^2\bbb^2$&$(1-\bbb,\bbb,\frac\bbb2+\frac4f,1-\frac\bbb2)$, $0<\bbb<1$
				&$2n_1+n_4=2n_2+n_3$&No $\ccc^2\cdots$, a contradiction\\
				\hline 
				$\bbb\ddd^2,\aaa^2\bbb\ccc$&$(\aaa,2-2\aaa-\frac8f,\frac8f,\aaa+\frac4f)$, $0<\aaa<1-\frac4f$
				&$n_1+n_4=2n_2$&$\text{AVC}\sub\{\bbb\ddd^2,\aaa^2\bbb\ccc,\ccc^{\frac f4}\}$\\
				\hline 
				$\bbb\ddd^2,\aaa^2\ccc^2$&$(\aaa,\frac8f,1-\aaa,1-\frac4f)$, $0<\aaa<1$
				&$n_1=n_3$&
				$\text{AVC}\sub\{\bbb\ddd^2,\aaa^2\ccc^2,\bbb^{\frac f4}\}$\\
				\hline 
				$\bbb\ddd^2,\bbb^4$&$(\aaa,\frac12,\frac34-\aaa+\frac4f,\frac34)$, $0<\aaa<\frac34+\frac4f$
				&$n_1=n_3$&$\text{AVC}\sub\{\bbb\ddd^2,\aaa^2\ccc^2,\bbb^4\}$\\
				\hline 
				$\bbb\ddd^2,\bbb^3\ccc$&$(-1+\frac{5\bbb}{2}+\frac4f,\bbb,2-3\bbb,1-\frac\bbb2)$, $\frac25-\frac{8}{5f}<\bbb<\frac23$
				&$5n_1+2n_2=6n_3+n_4$&No $\ccc^2\cdots$, a contradiction\\
				\hline  
				$\bbb\ddd^2,\bbb^2\ccc^2$&$(\frac\bbb2+\frac4f,\bbb,1-\bbb,1-\frac\bbb2)$, $0<\bbb<1$
				& $n_1+2n_2=2n_3+n_4$&No $\aaa^2\cdots$, a contradiction \\
				\hline  
				$\bbb\ddd^2,\bbb\ccc^3$&$(\frac13-\frac\bbb6+\frac4f,\bbb,\frac23-\frac\bbb3,1-\frac\bbb2)$, $0<\bbb<2$
				& $n_1+2n_3+3n_4=6n_2$&No $\aaa^2\cdots$, a contradiction\\
				\hline 
				$\bbb\ddd^2,\ccc^4$&$(\frac12-\frac\bbb2+\frac4f,\bbb,\frac12,1-\frac\bbb2)$, $0<\bbb<1+\frac8f$
				&$n_1+n_4=2n_2$&$\text{AVC}\sub\{\bbb\ddd^2,\aaa^2\bbb\ccc,\ccc^4\}$\\
				\hline 
		\end{tabular}}
		\caption{All degree $3$ vertices being $\bbb\ddd^2$.}\label{Tab-bd2}     
	\end{table*}

	\subsubsection*{Case \{$\bbb\ddd^2,\aaa^4$\}}	
	The AAD of $\bbb\ddd^2$ gives $\ccc^2\cdots$, denoted by $\bn$, $n_3\ge2$.	
	By $2n_2=n_3+n_4$, we get $n_2\ge1$. 
	
	If $n_2=1$, by $2=n_3+n_4\ge2$, then $n_3=2$ and $n_4=0$, By Parity Lemma, we get $n_1\ge2$, contradicting $2\aaa+\bbb+2\ccc>2$.
	
	If $n_2\ge2$, by $2n_2=n_3+n_4$ and $\bbb\ddd^2$, we get either  $n_4=0,n_3\ge4$ or $n_4=1,n_3\ge3$, contradicting $2\bbb+4\ccc>2$ or $2\bbb+3\ccc+\ddd>2$.	
	
	\subsubsection*{Case \{$\bbb\ddd^2,\aaa^2\bbb^2$\}}
	The AAD of $\bbb\ddd^2$ gives $\ccc^2\cdots$, denoted by $\bn$, $n_3\ge2$.	
	
	If $n_1=0$, by $n_4=2n_2+n_3\ge2$, then $n_4\ge2$, contradicting $2\ccc+2\ddd>2$. 
	
	If $n_1=1$, by Parity Lemma and Lemma \ref{anglesum}, we get $n_4\ge1$ and $n_2=0$. By $2+n_4=n_3$, we have $n_3\ge3$, contradicting $\aaa+3\ccc+\ddd>2$.
	
	If $n_1\ge2$, by $2n_1+n_4=2n_2+n_3$, we get either $n_2\ge1$ or $n_3\ge4$, contradicting $2\aaa+\bbb+2\ccc>2$ and $2\aaa+4\ccc>2$.
		
	
	\subsubsection*{Case \{$\bbb\ddd^2,\aaa^2\bbb\ccc$\}}	
	We will find all possible vertices $\bn$.  By $n_1+n_4=2n_2$, if $n_2=0$, then $n_1=n_4=0$ and $\bn=(0\,0\,\frac f4\,0)$ for $f=4k(k\ge4)$.
	
	If $n_2=1$, by $n_1+n_4=2$, we get $\bn=(2\,1\,n_3\,0),(1\,1\,n_3\,1)$ or $(0\,1\,n_3\,2)$. 
    By $\bbb\ddd^2$, $\aaa^2\bbb\ccc$ and Lemma \ref{anglesum}, we get $\bn=(2\,1\,1\,0)$ or $(0\,1\,0\,2)$. 
	
	If $n_2\ge2$, by $\bbb\ddd^2$ and $n_1+n_4=2n_2$, we get either $n_4=1,n_1\ge3$ or $n_4=0,n_1\ge4$, contradicting $3\aaa+2\bbb+\ddd>2$ and $4\aaa+2\bbb>2$. 
	
	Therefore, $\text{AVC}\sub\{\bbb\ddd^2,\aaa^2\bbb\ccc,\ccc^{\frac{f}{4}}(f\ge16)\}$. There is only one solution satisfying Balance Lemma: $\{\frac f2\,\bbb\ddd^2,\frac f2\,\aaa^2\bbb\ccc,2\ccc^{\frac{f}{4}}\}$.
	
	We have $(\aaa,\bbb,\ccc,\ddd)=(\aaa,2-2\aaa-\frac2k,\frac2k,\aaa+\frac1k)$ and $f=4k(k\ge4)$.
	By Lemma $\ref{geometry1}'$, we get $\aaa<1-\frac2k$. Then $\ddd=\frac1k+\aaa<1$. By Lemma \ref{geometry4} and Lemma $\ref{geometry4}'$, we get $\frac12-\frac2k<\aaa<\frac12+\frac1k$. Therefore, $\frac12-\frac2k<\aaa<\text{min}\{\frac12+\frac1k,1-\frac2k\}$.		
	
	By $\sin(\aaa+\frac\ccc2)\sin\frac\bbb2\neq0$ and $-\sin\frac\ccc2 \sin(\ddd+\frac{\bbb}{2})=0$, the equation \eqref{4-8} never holds, and it must satisfy \eqref{4-7}: $\sin(\aaa-\frac\ccc2)\sin\frac\bbb2=\sin\frac\ccc2 \sin(\ddd-\frac{\bbb}{2})$. This implies the following quadratic equation for $\cos^2\aaa$.
	\[A\cos^{4}\aaa+B\cos^{2}\aaa+C=0,\]
	\begin{align*}			
	&A=8 \cos^{3} \tfrac{\pi}{k}-4 \cos^{2} \tfrac{\pi}{k}-8 \cos \tfrac{\pi}{k}+5
	,\quad C=\cos^{2}\tfrac{\pi}{k} \left(2 \cos^{2} \tfrac{\pi}{k}+\cos \tfrac{\pi}{k}-2\right)^{2},\\
	&B=-8 \cos^5 \tfrac{\pi}{k}+4 \cos^{3}\tfrac{\pi}{k}+\cos^{2}\tfrac{\pi}{k}+4\cos\tfrac{\pi}{k}-4
	.
    \end{align*}	
	
	Then we get four roots for $\cos\aaa$:
	
	\[x_{1,2}=\pm\cos \! \tfrac{\pi}{k},\quad x_{3,4}=\pm\frac{2 \cos^2 \frac{\pi}{k}+\cos \frac{\pi}{k}-2}{\sqrt{8 \cos ^3 \frac{\pi}{k}-4 \cos^2\frac{\pi}{k}-8 \cos\frac{\pi}{k}+5}}.\]
	
	If $k\ge6$, then $\aaa>\frac12-\frac2k\ge\frac1k=\frac\ccc2$. By \eqref{4-7}, we get $\ddd>\frac\bbb2$, which implies  $\aaa>\frac12-\frac1k$. Then we have $\frac12-\frac1k<\aaa<\frac12+\frac1k$, or equivalently $|\cos\aaa|<\cos (\tfrac12-\tfrac{1}{k})\pi$. However, $|x_{1,2}|=\cos\frac\pi k>\cos (\tfrac12-\tfrac{1}{k})\pi$ when $k\ge6$; and $|x_{3,4}|>\cos (\tfrac12-\tfrac{1}{k})\pi$ by the graph of $y(k)=|x_{3,4}|-\cos (\tfrac12-\tfrac{1}{k})\pi$ in Fig.\,\ref{420}. Then there is no proper solution for $\aaa$ and such quadrilateral does not exist for $k\ge6$.

	\begin{figure}[htp]
		\centering
		\includegraphics[scale=0.4]{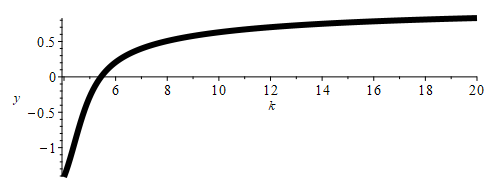}						
		\caption{The graph for the positive and increasing function $y(k)$, $k\ge6$.} 
		\label{420}	
	\end{figure}
	
	 If $k=5$, then $\aaa=\frac15,\frac25,\frac35$ or $\frac45$, contradicting \eqref{4-7} or $\frac{1}{10}<\aaa<\frac35$. 
	 
	 If $k=4$, then we get two double roots $\aaa=\frac14$ or $\frac34$. By  $\aaa<1-\frac2k=\frac12$, we get $\aaa=\frac14$, $\bbb=1$, $\ccc=\ddd=\frac12$. This  is Case $(1,4,2,2)/4$ of Table $1$ in \cite{lw}, which admits two tilings with the same vertices $\{8\bbb\ddd^2,8\aaa^2\bbb\ccc,2\ccc^4\}$. 
	
					
	\subsubsection*{Case \{$\bbb\ddd^2,\aaa^2\ccc^2$\}}	
	By $f\ge8$, we have $\ddd\ge\frac{1}{2}$. We will find all possible vertices $\bn$.
	
	If $n_1=0$, then $n_3=n_1=0$. By $\bbb\ddd^2,\ddd\ge\frac12$ and Parity Lemma, we get $\bn=(0\,\frac f4\,0\,0),(0\,1\,0\,2)$ or $(0\,0\,0\,4)$.

	If $n_1=1$, then $n_3=n_1=1$. By Parity Lemma and Lemma \ref{anglesum}, we get $n_2=0$ and $n_4\ge3$, contradicting $\aaa+\ccc+3\ddd>2$.
	
	If $n_1\ge2$, then $n_3=n_1\ge2$. By $\aaa^2\ccc^2$, we get $\bn=(2\,0\,2\,0)$.
	
	Therefore, $\text{AVC}\subset\{\bbb\ddd^2,\aaa^2\ccc^2,\ddd^4,\bbb^{\frac f4}(f\ge16)\}$. By Balance Lemma, $\bbb^{\frac f4}$ must appear and $f\ge16$. Then $\ddd^4$ is not a vertex. There is only one solution satisfying Balance Lemma: $\{\frac{f}{2}\,\bbb\ddd^2,\frac{f}{2}\,\aaa^2\ccc^2,2\bbb^{\frac f4}\}$. Let $f=4k(k\ge4)$, then $(\aaa,\bbb,\ccc,\ddd)=(\aaa,\frac2k,1-\aaa,1-\frac1k)$ and
	the quadrilateral is convex. By $k\ge4$, we have $\bbb=\frac2k\le\frac12<\frac34\le1-\frac1k=\ddd$. By Lemma $\ref{geometry3}'$, we get $\aaa>\ccc$, which implies $\aaa>\frac12$. By Lemma \ref{relation}, we have $\aaa+1-\frac1k=\aaa+\ddd<1+\bbb=1+\frac2k$. Then $\aaa<\frac{3}{k}$. By $\aaa>\frac12$, we get $k=4,5$.	
	If $k=5$, by \eqref{4-7}, we deduce that $\cos\aaa=\frac{1-\sqrt{5}}{4}$. Then
	$\aaa=\frac35=\ddd$, contradicting Proposition $9$ in \cite{lw}. Therefore, $k=4\,(f=16)$.	
	By Lemma \ref{calculate-1}, we get
	\begin{equation*}
		\aaa=1-\arcsin{\tfrac{\sqrt{2\sqrt2+1}}{2}},\,\,\,\,
		a=\arccos{\tfrac{1}{\sqrt{2\sqrt2+1}}},\,\,\,\,
		b=\arccos{\tfrac{7}{(2\sqrt2+1)^{\frac32}}},
	\end{equation*}
	$\aaa\approx 0.5664,a\approx 0.3292,b\approx 0.1158$.	
	
	In Fig. \ref{ad2 b2c2}, the unique AAD $\bbb^4=\thin^{\ccc}\bbb_1^{\aaa}\thin^{\ccc}\bbb_2^{\aaa}\thin\cdots$ determines $T_1,T_1'$. Then $\ddd_1\cdots=\bbb_2\ddd_1\ddd_3$ and $\thick\aaa\thin\ccc_1\thin\cdots=\thick\aaa\thin\ccc_1\thin\ccc_2\thin\cdots$ determine $T_2,T_3$. Similarly, we can determines $T_{2'},T_{3'}$. Then $\bbb_3\ddd_{2'}\cdots=\bbb_3\ddd_{2'}\ddd_4$ determines $T_4$; Similarly, we can determines $T_{4'}$. After repeating the process $2$ times, we get a unique earth map tiling $T(8\bbb\ddd^2,8\aaa^2\ccc^2,2\bbb^4)$.
	\begin{figure}[htp]
		\centering
		\begin{tikzpicture}[>=latex,scale=0.5]
			
			\foreach \b in {0,1,2,3}
			{
				\begin{scope}[xshift=4*\b cm]	
					\draw (0,4)--(0,0)
					(0,2)--(2,2)--(2,0)--(4,0)--(4,4)
					(2,0)--(2,-2)
					(6,0)--(6,-2);
					
					\draw[line width=1.5]	
					(0,0)--(2,0)
					(2,2)--(4,2)
					(4,0)--(6,0);

					\node at (2,4.2) {\small $\bbb$};\node at (2,2.45) {\small $\ddd$};
					\node at (0.4,2.45) {\small $\ccc$}; \node at (3.6,2.45) {\small $\aaa$};
					\node at (1.6,1.55) {\small $\bbb$}; \node at (2.4,1.55) {\small $\ddd$};
					\node at (0.4,1.55) {\small $\ccc$}; \node at (0.4,0.45) {\small $\ddd$};
					\node at (1.6,0.45) {\small $\aaa$}; \node at (2.4,0.45) {\small $\ccc$};
					\node at (3.6,0.45) {\small $\bbb$};\node at (3.6,1.55) {\small $\aaa$};
					\node at (2.4,-0.45) {\small $\ccc$};  \node at (4,-0.45) {\small $\ddd$}; 
					\node at (5.6,-0.45) {\small $\aaa$};  \node at (4,-2.2) {\small $\bbb$}; 
					
				\end{scope}
			}

			\node[draw,shape=circle, inner sep=0.5] at (2,3.3) {\small $1$};
			\node[draw,shape=circle, inner sep=0.5] at (2+4,3.3) {\small $1'$};
			\node[draw,shape=circle, inner sep=0.5] at (1,1) {\small $2$};
			\node[draw,shape=circle, inner sep=0.5] at (5,1) {\small $2'$};
			\node[draw,shape=circle, inner sep=0.5] at (3,1) {\small $3$};
			\node[draw,shape=circle, inner sep=0.5] at (7,1) {\small $3'$};
			\node[draw,shape=circle, inner sep=0.5] at (4,-1.3) {\small $4$};
			\node[draw,shape=circle, inner sep=0.5] at (8,-1.3) {\small $4'$};

		\end{tikzpicture}	
		\caption{$f=16$, $T(8\bbb\ddd^2,8\aaa^2\ccc^2,2\bbb^{4})$.} \label{ad2 b2c2}
	\end{figure}
	\subsubsection*{Case \{$\bbb\ddd^2,\bbb^4$\}}
    We will find all possible vertices $\bn$.	
	If $n_1=0$, then $n_3=n_1=0$. By $\bbb=\frac12,\ddd=\frac34$ and Parity Lemma, we get $\bn=(0\,4\,0\,0)$ or $(0\,1\,0\,2)$.	
	
	If $n_1=1$, then $n_3=n_1=1$. By Parity Lemma and Lemma \ref{anglesum}, we get $n_2=0$ and $n_4\ge3$, contradicting $\aaa+\ccc+3\ddd>2$.
	
	If $n_1\ge2$, by $n_1=n_3$ and $R(\aaa^2\ccc^2\cdots)<\aaa+\ccc,\bbb,\ddd$, we get $\bn=(2\,0\,2\,0)$.
	
	Therefore, $\text{AVC}\sub\{\bbb\ddd^2,\aaa^2\ccc^2,\bbb^4\}$, which has been discussed by Case $\{\bbb\ddd^2,\aaa^2\ccc^2\}$ and gives a unique tiling $T(8\bbb\ddd^2,8\aaa^2\ccc^2,2\bbb^{4})$ in Fig.\,\ref{ad2 b2c2}.	
	
	\subsubsection*{Case \{$\bbb\ddd^2,\bbb^3\ccc$\}}	
	The AAD of $\bbb\ddd^2$ gives $\ccc^2\cdots$, denoted by $\bn$, $n_3\ge2$.
	
	If $n_2=0$, by $5n_1=6n_3+n_4\ge12$, then $n_1\ge3$. By $3\aaa+2\ccc+\ddd>2$, we get $n_4=0$. Therefore, $5n_1=6n_3$. Then $n_1\ge6$ and $n_3\ge5$, contradicting $6\aaa+5\ccc>2$. 
	
	If $n_2\ge1$, by $\bbb^3\ccc$ and $5n_1+2n_2=6n_3+n_4\ge12$, we deduce that $n_2=1,2$ and $n_1\ge2$, contradicting $2\aaa+\bbb+2\ccc>2$.
	
	
	\subsubsection*{Case \{$\bbb\ddd^2,\bbb^2\ccc^2$\}}	
	The AAD of $\bbb\ddd^2$ gives $\aaa^2\cdots$, denoted by $\bn$, $n_1\ge2$.
	
	If $n_3=0$, then $n_4=n_1+2n_2\ge2$, contradicting $2\aaa+2\ddd>2$. 
	 
	If $n_3=1$, by $n_1+2n_2=2+n_4$ and $\text{deg}\,\bn\ge4$, then $n_4\ge1$. By $2\aaa+2\ddd>2$ and Parity Lemma, we get $n_4=1$ and $n_1\ge3$, contradicting $3\aaa+\ccc+\ddd>2$.
	 
	If $n_3\ge2$, by $n_1+2n_2=2n_3+n_4$, we get either $n_2\ge1$ or $n_1\ge4$, contradicting $2\aaa+\bbb+2\ccc>2$ and $4\aaa+2\ccc>2$. 
	 	
	
	\subsubsection*{Case \{$\bbb\ddd^2,\bbb\ccc^3$\}}
	The AAD of $\bbb\ddd^2$ gives $\aaa^2\cdots$, denoted by $\bn$, $n_1\ge2$. By $6n_2=n_1+2n_3+3n_4\ge2$, we get $n_2\ge1$.
	
	If $n_4=0$, by $n_1+2n_3=6n_2\ge6$ and Parity Lemma, we get $n_1\ge6$, or $n_1=4$ and $n_3\ge1$, or $n_1=2$ and $n_3\ge2$, contradicting $6\aaa+\bbb>2,4\aaa+\bbb+\ccc>2$ and $2\aaa+\bbb+2\ccc>2$.
	
	If $n_4\ge1$, by $\bbb\ddd^2$ and Parity Lemma, then $n_4=1$ and $n_1\ge3$,  contradicting $3\aaa+\bbb+\ddd>2$.
	
	
	\subsubsection*{Case \{$\bbb\ddd^2,\ccc^4$\}}	
	We will find all possible vertices $\bn$. 
	If $n_2=0$, by $n_1+n_4=2n_2$, then $n_1=n_4=0$. By $\ccc^4$, we get $\bn=(0\,0\,4\,0)$.
	
	If $n_2=1$, by $n_1+n_4=2$, we get $\bn=(2\,1\,n_3\,0),(1\,1\,n_3\,1)$ or $(0\,1\,n_3\,2)$. 
	By $\bbb\ddd^2$, $2\aaa+\bbb+2\ccc>2$ and Lemma \ref{anglesum}, we get $\bn=(2\,1\,1\,0)$ or $(0\,1\,0\,2)$. 
	
	If $n_2\ge2$, by $\bbb\ddd^2$ and $n_1+n_4=2n_2$, we get either $n_4=1,n_1\ge3$ or $n_4=0,n_1\ge4$, contradicting $3\aaa+2\bbb+\ddd>2$ and $4\aaa+2\bbb>2$. 
	
	Therefore, $\text{AVC}\sub\{\bbb\ddd^2,\aaa^2\bbb\ccc,\ccc^4\}$, which has been discussed by Case $\{\bbb\ddd^2,\aaa^2\bbb\ccc\}$ and gives a rational quadrilateral admitting two tilings.
\end{proof}

\begin{remark}
	Combined with the classification of rational $a^3b$-tilings in \cite{lw}, the only other two tilings with $\bbb\ddd^2$ as the unique degree $3$ vertex type come from Case $(1,4,2,2)/4$ of Table $1$ in \cite{lw}. Both tilings have the same vertices $\{8\bbb\ddd^2,8\aaa^2\bbb\ccc,2\ccc^4\}$ with $f=16$.

\end{remark}

\section{Conclusion}\label{conclusion}
Table \ref{Tab-summary-1} and \ref{Tab-summary-2} present Sakano-Akama's classification \cite{sa} for type $a^2 b^2$ and $a^4$ using our notations.

\begin{table*}[htp]                        
	\centering     

	~\\ 
	\begin{tabular}{c|c|c}	 
		
	 $f$&$(\aaa,\ddd,\ccc,\ddd)$ &Tilings\\
	\hline				
	\multirow{2}{*}{$24$}& \multirow{2}{*}{$(3,3,4,3)/6$}& $T(8\ccc^3,6\aaa^4, 12\ddd^4)$\\		
	\cline{3-3}
	&& $T(8\ccc^3,2\aaa^4,8\aaa^2\ddd^2,8\ddd^4)$\\	
	\cline{1-3}
    $60$&$(12,15,20,15)/30$& $ T(20\ccc^3,30\ddd^4,12\aaa^5)$\\
	\hline
	$\ge6$ (even)&$(2-2\ddd,\ddd,\frac4f,\ddd),\ddd\in(0,1)\setminus\{\frac{f-2}{f}\}$& $ T(f\aaa\ddd^2,2\ccc^{\frac{f}{2}})$\\		
	\hline	
	\end{tabular}
	\caption{All tilings of type $a^2b^2$. }\label{Tab-summary-1}        
\end{table*}

\begin{table*}[htp]                        
	\centering     

	~\\ 
	\begin{tabular}{c|c|c}	 
			
			$f$&$(\aaa,\bbb,\aaa,\bbb)$ &Tilings\\
			\hline			
			$6$&$(2,2,2,2)/3$& $ T(8\aaa^3),\,\bbb=\aaa$\\
			\cline{1-3}
			$12$&$(4,3,4,3)/6$& $ T(8\aaa^3,6\bbb^4)$\\
			\cline{1-3}
			$30$&$(10,6,10,6)/15$& $ T(20\aaa^3,12\bbb^5)$\\
			\hline
			\multirow{2}{*}{$\ge8$  (even)}&\multirow{2}{*}{$(4,f-2,4,f-2)/f$}& $ T(f\aaa\bbb^2,2\aaa^{\frac{f}{2}})$\\	
			\cline{3-3}		
			&&$f=4k+2$: $T((f-2)\aaa\bbb^2,4\aaa^{\frac{f+2}{4}}\bbb)$ \\			
			\hline
	\end{tabular}
	\caption{All tilings of type $a^4$.}\label{Tab-summary-2}        
\end{table*}

\begin{figure}[htp]
	\centering
	\begin{tikzpicture}

		\draw
		(-0.8,-0.8) -- (-0.8,0.8) -- (0.8,0.8) -- (0.8,-0.8)-- (-0.8,-0.8);					
		
		\node at (-0.5,0.5) {\small $\aaa$};
		\node at (0.5,0.5) {\small $\aaa$};
		\node at (-0.5,-0.5) {\small $\aaa$};
		\node at (0.5,-0.5) {\small $\aaa$};
		
		\begin{scope}[scale=0.4,xshift=-14cm,yshift=-1cm]	
			\draw (-2,0)--(0,5)--(2,0);			
			\draw[line width=1.5]
			(-2,0)--(0,-2)--(2,0);	    
			
			\node at (0,3.8) {\small $\ccc$};
			\node at (0,-1) {\small $\aaa$};
			\node at (-1.4,0.2) {\small $\ddd$};
			\node at (1.4,0.2) {\small $\ddd$};	
		\end{scope}
		\begin{scope}[scale=0.36,xshift=-8cm,yshift=0.8cm]	
			\draw (-2,0)--(0,4)--(2,0)			
			(-2,0)--(0,-4)--(2,0);	   
			
			\node at (0,2.8) {\small $\aaa$};
			\node at (0,-2.8) {\small $\aaa$};
			\node at (-1.4,0) {\small $\bbb$};
			\node at (1.4,0) {\small $\bbb$}; 	
		\end{scope}
		
		\begin{scope}[scale=0.4,xshift=-21cm,yshift=-1cm]	
			\draw (-2,-2)--(0,5)--(2,-2);			
			\draw[line width=1.5]
			(-2,-2)--(0,-2)--(2,-2);	    
			
			\node at (0,3.6) {\small $\ccc$};
			\node at (0,0.4-2) {\small $\aaa$};
			\node at (-1.4,0.4-2) {\small $\ddd$};
			\node at (1.4,0.4-2) {\small $\ddd$};	
			\fill (0,-2) circle (0.1);
		\end{scope}	
		
		\begin{scope}[scale=0.4,xshift=-28cm,yshift=-1cm]	
			\draw (-2,-2)--(0,5)--(2,-2);			
			\draw[line width=1.5]
			(-2,-2)--(0,0)--(2,-2);	    
			
			\node at (0,3.6) {\small $\ccc$};
			\node at (0,0.5) {\small $\aaa$};
			\node at (-1.25,-0.75) {\small $\ddd$};
			\node at (1.28,-0.75) {\small $\ddd$};	
		\end{scope}	
		
	\end{tikzpicture} 
	\caption{darts, triangles, kites, rhombi and squares.}
	\label{summary}
\end{figure}
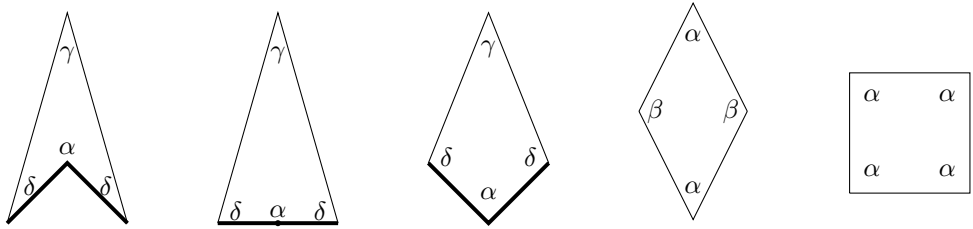

The infinite sequence in Table \ref{Tab-summary-1} are $1$-parameter families of $2$-layer earth map tilings represented by the dashed line $AA'$ inside the $2$-dimensional moduli in Fig.\,\ref{modular}. Moving from $A$ to $A'$, the quadrilateral has the shape of concave dart, degenerate triangle, and convex kite in Fig.\,\ref{summary}. Case $f=24$ and $f=60$ in Table \ref{Tab-summary-1} come from the standard quadrilateral subdivisions of the octahedron and icosahedron, together with a  modification. The standard quadrilateral subdivision divides each $n$-gon face into $n$ quadrilaterals by linking the center of this face to all middle points of its edges. Such subdivision of the tetrahedron gives an $a^4$-tiling with $f=12$ in Table \ref{Tab-summary-2}.

The infinite sequence in Table \ref{Tab-summary-2} is represented by the point $V$ in the moduli in Fig.\,\ref{modular}, and each quadrilateral becomes a rhombus in Fig.\,\ref{summary} which also gives the unique square tiling when $f=6$. For $f=4k+2\ge10$, these $2$-layer earth map tilings admit a basic flip modification. Case $f=6$, $f=12$ and $f=30$ in Table \ref{Tab-summary-2} come from the quadricentric subdivisions of the tetrahedron, octahedron and icosahedron. The quadricentric subdivision constructs each quadrilateral face using two vertices of an old edge and two centers of the two faces adjacent to this edge, so that this old edge becomes a diagonal of the new quadrilateral face. 

Together with \cite{sa} and our two papers \cite{lpwx,lw}, this paper completes the classification of  edge-to-edge spherical tilings by congruent quadrilaterals. We summarize them into four classes: 
\begin{enumerate}
	\item A sequence of $2$-parameter families of $2$-layer earth map tilings in Fig.\,\ref{abd-1} by any even $f\ge6$ tiles,  with the moduli in Fig.\,\ref{modular} showing the reductions from type $a^2 bc$ to $a^2 b^2$ (the dashed line $AA'$),  $a^3 b$ (the dotted curve $UVW$), and $a^4$ (the point $V$), among which  
	\begin{itemize}
		\item special $a^3b$-quadrilaterals with $\bbb$ being an integer multiple of $\ccc$ (i.e. some special points on the dotted curve $UVW$ counted by $\mathcal{Q}_1 (f)$ in Table \ref{Tab-1.3}) admit basic flip modifications in Table \ref{tab-3} and Fig.\,\ref{flip1}; 
		\item more special $a^3 b$-quadrilaterals with $4$ angles $(\frac 2f,\frac{4f-4}{3f},\frac{4}{f},\frac{2f-2}{3f})$ for $f=6k+4\ge10$ admit $3$ additional modifications in  \cite{lw} Table $2$ or  \cite{lw} Section $6$; 
		\item $a^4$-quadrilaterals with $4$ angles $(\frac{4}{f},\frac{f-2}{f},\frac{4}{f},\frac{f-2}{f} )$ for $f=4k+2\ge10$ admit a basic flip modification in Table \ref{Tab-summary-2}. 
	\end{itemize}		
	
	\item A $1$-parameter family of quadrilateral subdivisions of the octahedron in \cite{lpwx} Fig.\,$36$ by $24$ tiles, with the moduli $b\in (0,\frac{\pi}{4}]$ having the reductions from type $a^2 bc$ to $a^2 b^2$ ($b=\frac{\pi}{4}$) and $a^3 b$ ($b=\arctan(\sqrt{3}-1)\approx0.201 \pi$) in Fig.\,\ref{case a3 d4}, among which 
	\begin{itemize}
		\item case $b=\frac{\pi}{4}$ admits a  modification in Table \ref{Tab-summary-1};
		\item case $b=\arctan(\frac{3-\sqrt{5}}{2})\approx0.116\pi$ ($\bbb=\frac{2\pi}{3}=2\ccc$) admits a flip modification in the $3$rd picture of \cite{lpwx} Fig.\,$2$ and $29$.
	\end{itemize}   
	
	\item A sequence of $3$-layer earth map tilings by $f=8k$ tiles for $k\ge2$ (each has a unique $a^2 bc$-quadrilateral with $4$ angles $(\frac{f-8}{f},\frac8f,\frac{f+8}{2f},\frac12)$) in \cite{lpwx} Fig.\,$16$, among which there are exactly two flip modifications for each odd $k$ in \cite{lpwx} Fig.\,$26$ and $28$. The quadrilateral for $f=24 \, (k=3)$ coincides with case $b=\arctan(\frac{3-\sqrt{5}}{2})$ ($\bbb=\frac{2\pi}{3}=2\ccc$) in the previous class and admits $2+3=5$ tilings in total, as shown in \cite{lpwx} Fig.\,$2$.  	
	
	\item $11$ sporadic tilings: 
	\begin{itemize}
		\item $3$ from platonic solids in \cite{sa} Fig.\,$1$: Case $f=60$ of type $a^2b^2$ in Table \ref{Tab-summary-1} as the quadrilateral subdivision of the icosahedron, and 	  
		Case $f=12,30$ of type $a^4$ in Table \ref{Tab-summary-2} as the quadricentric subdivisions of the octahedron and icosahedron;	
		\item $7$ earth map tilings of type $a^3b$: Case $f=12,16,16,16$ in Table \ref{tab-1} and Fig.\,\ref{a2}, and Case $f=16,16,36$ in \cite{lw} Table $1$ and Fig.\,$2$;	
		\item Case $f=36$ of type $a^3b$ with $4$ angles  $(15,6,10,7)/18$ in \cite{lw} Table $1$ and the last picture of Fig.\,$2$. This is the only tiling, among all edge-to-edge triangular, quadrilateral, pentagonal tilings of the sphere,  which has no apparent relation with any platonic solids or earth map tilings.
	\end{itemize}		
\end{enumerate}

We also summarize all non-edge-to-edge triangular tilings induced from the above list: 
\begin{enumerate}
	\item All degenerate $2$-layer earth map tilings represented by the great arcs $EF$, $A'E$ and $A'F$ inside the $2$-dimensional moduli in Fig.\,\ref{modular}.
	\item Case $f=16,16$ with $4$ angles $(1,4,2,2)/4$ in \cite{lw} Table $1$ and Fig.\,$2,35$ and $36$. In fact, the triangle is half of the regular triangle face of the octahedron. One can cut each face to $2$ such triangles in $3$ different ways, most of which induce non-edge-to-edge tilings and admit $1$-parameter rotations along any full great circle boundary. 
	\item Special $a^3 b$-quadrilaterals with $4$ angles $(\frac 2f,\frac{4f-4}{3f},\frac{4}{f},\frac{2f-2}{3f})$ for any even $f\ge6$ satisfy $\ccc=2\aaa,\bbb=2\ddd$ and can be subdivided into $3$ congruent triangles as shown in the third picture of \cite{lw} Fig.\,$6$. Each tiling in the second case of \cite{lw} Table $2$ induces a non-edge-to-edge triangular tiling.
\end{enumerate}

These are new examples, comparing to early explorations of non-edge-to-edge triangular tilings in  \cite{D1,D2,DD1,DD2,DD3}. 

For monohedral spherical tilings, there remain three open problems: non-edge-to-edge tilings of the sphere by congruent triangles, quadrilaterals, or pentagons. They prove to be more challenging after a few early explorations, and will be our next study goal.


\begin{thebibliography}{1}
	
	
	\bibitem{ac}
	Y.~Akama, N.~Van Cleemput. 
	\newblock Spherical tilings by congruent quadrangles: Forbidden cases and substructures. 
	\newblock {\em Ars Math. Contemp.} 8 (2015), 297-318. 	
	
	\bibitem{awy}
	Y.~Akama, E.~Wang, M.~Yan.
	\newblock Tilings of sphere by congruent pentagons III: edge combination $a^5$. 
	\newblock {\em Adv. Math.}, 394 (2022), 107881. 
	
	\bibitem{coolsaet}
	K.~Coolsaet.
	\newblock Spherical tilings with three equal sides and rational angles. 
	\newblock {\em Ars Math. Contemp.}, 12:415–424, 2017.  	
	
	\bibitem{D1}
	Dawson, R. J. MacG., 
	\newblock {T}ilings of the Sphere with Isosceles Triangles.  
	\newblock \emph{Discrete Comput. Geom.}, 30 (2003) 467–487.

	\bibitem{D2}
	Dawson, R. J. MacG., 
	\newblock {A}n isosceles triangle that tiles the sphere in exactly three ways.  
	\newblock \emph{Discrete Comput. Geom.}, 30 (2003) 459–466.	
	
	\bibitem{DD1}
	Dawson, R. J. MacG. and Doyle, B., 
	\newblock {T}ilings of the sphere with right triangles I: the asymptotically right families. 
	\newblock \emph{Electronic Journal of Combinatorics}, 13(1) (2006) \#R48.
	
	\bibitem{DD2}
	Dawson, R. J. MacG. and Doyle, B., 
	\newblock {T}ilings of the sphere with right triangles II: the(1, 3, 2), (0, 2, n) subfamily.  
	\newblock \emph{Electronic Journal of Combinatorics}, 13(1) (2006) \#R49.
	
	\bibitem{DD3}
	Dawson, R. J. MacG. and Doyle, B., 
	\newblock {T}ilings of the sphere with right triangles III: the asymptotically obtuse families. 
	\newblock \emph{Electronic Journal of Combinatorics}, 13(1) (2007) \#R48.
	
	\bibitem{lpwx}
	Y.~Liao, P.~Qian, E.~Wang, Y.~Xu.
	\newblock Tilings of the sphere by congruent quadrilaterals I: edge combination $a^2bc$. 
	\newblock  \textit{preprint}, arXiv: 2110.10087, 2021.   
	
	
	
	\bibitem{lw}
	Y.~Liao, E.~Wang.
	\newblock Tilings of the sphere by congruent quadrilaterals II: edge combination $a^3b$ with rational angles. 
	\newblock  \textit{preprint}, arXiv: 2205.14936, 2022. 
	
	\bibitem{rao}
	M.~Rao.
	\newblock Exhaustive search of convex pentagons which tile the plane.
	\newblock  \textit{preprint}, arXiv: 1708.00274, 2017. 
	
	
	\bibitem{sa}
	Y.~Sakano, Y.~Akama. 
	\newblock Anisohedral spherical triangles and classification of spherical tilings by congruent kites, darts and rhombi. 
	\newblock {\em Hiroshima Mathematical Journal}, 2015, 45(3): 309-339.
	\bibitem{so}
	D.~M.~Y.~Sommerville.
	\newblock Division of space by congruent triangles and tetrahedra.
	\newblock {\em Proc. Royal Soc. Edinburgh}, 43:85--116, 1924.
	
	\bibitem{ua}
	Y.~Ueno, Y.~Agaoka. 
	\newblock Classification of tilings of the 2-dimensional sphere by congruent triangles.
	\newblock {\em Hiroshima Math. J.}, 32(3):463--540, 2002.
	
	\bibitem{ua2}
	Y.~Ueno, Y.~Agaoka. 
	\newblock Examples of spherical tilings by congruent quadrilaterals.
	\newblock {\em Math. Inform. Sci., Fac. Integrated Arts and Sci., Hiroshima Univ., Ser. IV}, 27:135--144, 2001.
	
	\bibitem{wy1}
	E.~Wang, M.~Yan.
	\newblock Tilings of the sphere by congruent pentagons I: edge combinations $a^2b^2c$ and $a^3bc$. 
	\newblock {\em Adv. Math.}, 394 (2022), 107866.
	
	\bibitem{wy2}
	E.~Wang, M.~Yan.
	\newblock Tilings of the sphere by congruent pentagons II: edge combination $a^3b^2$. 
	\newblock {\em Adv. Math.}, 394 (2022), 107867.
	
	\bibitem{wy3}
	E.~Wang, M.~Yan.
	\newblock Moduli of pentagonal subdivision tiling.
	\newblock {\em preprint}, arXiv: 1907.08776, 2019.
	
	
	
\end{thebibliography}
\end{document}